\documentclass{article}

\usepackage[latin1]{inputenc}
\usepackage[english]{babel}

\usepackage{amsmath}
\usepackage{amsthm}
\usepackage{amssymb}
\usepackage{mathtools}

\usepackage[notref, notcite, final]{showkeys}
\usepackage{color}
\usepackage{graphicx}
\usepackage[hidelinks]{hyperref}
\usepackage{cleveref}
\usepackage{comment}
\usepackage{doi}

\usepackage[shortlabels]{enumitem}

\newcommand{\id}{\mathcal{I}}

\newcommand{\nn}{\mathbb{N}}
\newcommand{\rr}{\mathbb{R}}
\newcommand{\cc}{\mathbb{C}}
\newcommand{\hh}{\mathbb{H}}
\renewcommand{\SS}{\mathbb{S}}

\renewcommand{\Re}{\mathrm{Re}}

\newcommand{\lhol}{\mathcal{SH}_L}
\newcommand{\rhol}{\mathcal{SH}_R}
\newcommand{\intrin}{\mathcal{SH}}

\newcommand{\SliceL}{\mathcal{SF}_L}
\newcommand{\SliceR}{\mathcal{SF}_R}
\newcommand{\SliceReal}{\mathcal{SF}}

\newcommand{\SCon}{\mathcal{SC}}

\newcommand{\I}{\mathsf{i}}
\newcommand{\J}{\mathsf{j}}
\newcommand{\K}{\mathsf{k}}

\newcommand{\sderiv}[1][]{\partial_{S#1}}

\newcommand{\ran}{\operatorname{ran}}
\newcommand{\dom}{\mathcal{D}}
\newcommand{\dist}{\operatorname{dist}}

\newcommand{\fdom}{\mathcal{D}}

\newcommand{\linspan}[1]{\operatorname{span}_{#1}}

\newcommand{\clos}[1]{\overline{#1}^{\mathsf{cl}}}

\newcommand{\hil}{\mathcal{H}} 
\newcommand{\sBorel}{\mathsf{B_S}}
\newcommand{\Borel}{\mathsf{B}}

\newcommand{\sMeas}{\mathcal{M}_s}
\newcommand{\bsMeas}{\mathcal{M}_s^{\infty}}
\newcommand{\bsIntrin}{\SliceReal_{0}^{\infty}}

\newcommand{\RRes}{\mathcal{R}}
\newcommand{\sgn}{\operatorname{sgn}}

\newcommand{\bw}{\mathbf{w}}
\newcommand{\bv}{\mathbf{v}}
\newcommand{\bu}{\mathbf{u}}
\newcommand{\bb}{\mathbf{b}}

\newcommand{\bof}{\mathbf{f}}
\newcommand{\bg}{\mathbf{g}}
\newcommand{\bh}{\mathbf{h}}
\newcommand{\bO}{\mathbf{0}}

\newcommand{\boundOP}{\mathcal{B}}
\newcommand{\closOP}{\mathcal{K}}
\newcommand{\Q}{\mathcal{Q}}

\newcommand{\diag}{\operatorname{diag}}

\numberwithin{equation}{section}

\theoremstyle{plain}
\newtheorem{theorem}{Theorem}[section]
\newtheorem{lemma}[theorem]{Lemma}
\newtheorem{proposition}[theorem]{Proposition}
\newtheorem{corollary}[theorem]{Corollary}

\theoremstyle{definition}
\newtheorem{definition}[theorem]{Definition}
\newtheorem{example}[theorem]{Example}
\newtheorem{notation}[theorem]{Notation}

\theoremstyle{remark}
\newtheorem{remark}[theorem]{Remark}
\crefname{enumi}{}{}

\title{Operator Theory on One-Sided Quaternionic Linear Spaces: Intrinsic $S$-Functional Calculus and Spectral Operators}

\author{Jonathan Gantner\\
Politecnico di Milano\\
Dipartimento di Matematica\\
Via E. Bonardi, 9\\
20133 Milano, Italy\\
jonathan.gantner@polimi.it}

\begin{document}

\maketitle

\begin{abstract}
Two major themes drive this article: identifying the minimal structure necessary to formulate quaternionic operator theory and revealing a deep  relation between complex and quaternionic operator theory. 

The theory for quaternionic right linear operators is usually formulated under the assumption that there exists not only a right- but also a left-multiplication on the considered Banach space $V$. This has technical reasons, as the space of bounded operators on $V$ is otherwise not a quaternionic linear space. A right linear operator is however only associated with the right multiplication on the space and in certain settings, for instance on quaternionic Hilbert spaces, the left multiplication is not defined a priori but must be chosen randomly. Spectral properties of an operator should hence be independent of the left multiplication on the space.

We show that results derived from functional calculi involving intrinsic slice functions can be formulated without the assumption of a left multiplication. We develop the $S$-functional calculus in this setting and a new approach to spectral integration of intrinsic slice functions. This approach has a clear interpretation in terms of the right linear structure on the space and allows to formulate the spectral theorem without using any randomly chosen structure. These techniques only apply to intrinsic slice functions, but this is an negligible restriction. Indeed, only these functions are compatible with the very basic intuition of a functional calculus that $f(T)$ should be defined by letting $f$ act on the spectral values of $T$.

Using the above tools we develop a theory of spectral operators and obtain results analogue to those of the complex theory. In particular we show the existence of a canonical decomposition of a spectral operator and discuss its behavior under the $S$-functional calculus. 

Finally, we show a beautiful relation with complex operator theory: if we embed the complex numbers into the quaternions and consider the quaternionic vector space as a complex one, then complex and quaternionic operator theory are consistent. Again it is the symmetry of intrinsic slice functions that guarantees that this compatibility is true for any possible imbedding of the complex numbers.

\end{abstract}

{\bf 2010 Mathematics Subject Classification.} 47A60, 47B40.\\
{\bf Key words and phrases.} quaternionic operators, spectral operators, $S$-functional calculus, spectral integration, $S$-spectrum. 

\section{Introduction}
Birkhoff and von Neumann  showed  1936  in their celebrated paper  \cite{Birkhoff:1936} titled {\em The Logic of Quantum Mechanics} that quantum mechanics can only be formulated in essentially two different settings, either on complex or on quaternionic Hilbert spaces. Since then mathematicians have tried to develop methods for studying quaternionic linear operators analogue to those for complex linear ones. First abstract results were obtained quite soon \cite{Teichmuller:1936,Viswanath:1971, Sharma:1987} but it was only ten to fifteen years ago that the fundamental principles of a rigorous theory of quaternionic linear operators were understood. In particular the definition of a proper notion of spectrum caused severe problems.

The spectrum $\sigma(A)$ of a closed operator $A$ on a complex Banach space generalizes the set of its eigenvalues. It is defined as the set of all complex numbers $\lambda$ such that the operator of the eigenvalue equation ${(\lambda\id - A)\bv} = \bO$ does not have a bounded inverse. If $\lambda\notin\sigma(A)$, the bounded inverse $R_{\lambda}(A) := (\lambda\id - A)^{-1}$ is called the resolvent of $A$ at $\lambda$ and $\lambda\mapsto R_{\lambda}(A)$ is an operator-valued holomorphic function on the resolvent set $\rho(A) := \cc\setminus\sigma(A)$ of $A$. Trying to generalize these concepts to the quaternionic setting by simply mimicking this procedure ends up in an unclear and ambiguous situation. Indeed, already the concept of  an eigenvalue is ambiguous. The quaternions $\hh$ are a number system that generalize the complex numbers, but their multiplication is not commutative. Hence, we can formulate two different eigenvalue problems for a quaternionic right linear operator $T$, i.e. an operator with the property that $T(\bu a + \bv) = T(\bu) a+ T(\bv)$ for any vectors $\bu$ and $\bv$ and any $a\in\hh$. We can consider right or left eigenvalues, which satisfy
\[
T\bv - \bv s = \bO \quad\text{resp.}\quad  T\bv - s\bv = \bO
\]
for some $\bv\neq \bO$. Since a right linear operator is only related to the right multiplication on the vector space,  right eigenvalues seem the natural notion of eigenvalue for such operator. Moreover, right eigenvalues are meaningful both in applications such as quaternionic quantum mechanics \cite{Adler:1995} and in the mathematical theory, for instance for showing the spectral theorem for matrices with quaternionic entries \cite{Farenick:2003}. Due to the non-commutativity of the quaternionic multiplication, right eigenvalues are however accompanied by strange phenomena. For a right eigenvector $\bv$ satisfying $T\bv = \bv s$ and $a\in\hh$ with $a\neq 0$, we namely have
\begin{equation}\label{RightEVRel}
T(\bv a) = (T\bv)a = (\bv s)a = (\bv a )(a^{-1}sa).
\end{equation}
Hence, $\bv a$ is not an eigenvector associated with $s$, but an eigenvector associated with $a^{-1}sa$. Right eigenvalues do therefore not appear individually, but in equivalent classes of the form
\begin{equation}\label{EVSym}
[s] = \{a^{-1}sa: a\in\hh, a\neq 0\}.
\end{equation}
Although right eigenvalues seem to constitute a meaningful notion of eigenvalues for quaternionic right linear operators, the operator of the right eigenvalue equation can moreover not be used to generalize this concept to some notion of spectrum as it is not quaternionic right linear because of \eqref{RightEVRel}. Left eigenvalues on the other hand are characterized by a right linear operator, but they do not seem useful neither in the mathematical theory nor in applications.

Further difficulties were caused by the fact that it was not clear on which notion of generalized holomorphicity one could build a theory of quaternionic linear operators. The fact that the resolvent is an operator-valued holomorphic function is essential in complex operator theory. The most successful notion of generalized holomorphicity in the quaternionic setting was the notion of Fueter regularity introduced by Fueter  in \cite{Fueter:1932,Fueter:1934}. The theory of Fueter regular functions was extremely successful in generalizing the theory of holomorphic functions, but it was not suitable for applications in operator theory because it does not contain power series in a quaternionic variable. 

Hence, little progress was made until these problems were resolved by the introduction of the theory of slice hyperholomorphic functions and the associated notion of $S$-spectrum. Slice hyperholomorphic functions are functions that are in some sense holomorphic and compatible with the symmetry of the set of eigenvalues in \eqref{EVSym}. As power series in the quaternionic variable are slice hyperholomorphic, the power series expansion of the classical resolvent 
\begin{equation}\label{SResSer}
S_L^{-1}(s,T) = \sum_{n=0}^{+\infty}T^n s^{-(n+1)}
\end{equation}
is slice hyperholomorphic in the scalar variable $s$. It is possible to find a closed form for this power series under the assumption that $T$ and $s$ do not commute, namely \begin{equation}\label{SRes1}
S_L^{-1}(s,T) = \Q_{s}(T)^{-1}(\overline{s}\id - T),
\end{equation}
 with
 \[
 \Q_{s}(T) = T^2 - 2\Re(s)T + |s|^2\id.
 \]
  This closed form allowed to define the (left) $S$-resolvent operator $S_L^{-1}(s,T)$ and in turn the $S$-resolvent set $\rho_{S}(T)$ and the $S$-spectrum $\sigma_{S}(T)$ of $T$, namely
\[
\rho_{S}(T) := \{s\in\hh: \Q_{s}(T)^{-1}\in\boundOP(V)\}\qquad\text{and}\qquad\sigma_{S}(T) := \hh\setminus\rho_{S}(T),
\]
where $\boundOP(V)$ denotes the set of all bounded quaternionic right linear operators on the quaternionic Banach space we consider. The $S$-spectrum is a proper generalization of the set of right eigenvalues, just as the spectrum generalizes the set of eigenvalues in complex operator theory. Indeed, $\Q_{s}(T)$ is the operator of a generalized eigenvalue equation that is associated not with the single eigenvalue $s$, but with the equivalence class~$[s]$. Any eigenvector associated with a right eigenvalue in $[s]$ belongs to the kernel of $\Q_{s}(T)$. Conversely,  any vector in $\ker\Q_{s}(T)$ is the linear combination of eigenvectors associated with different eigenvalues in $[s]$ as we show in \Cref{EVSplitLem}.

The introduction of these fundamental concepts was the starting point of a mathematically rigorous theory of quaternionic linear operators that expanded in several directions during the last decade. A first step  was the generalization of the Riesz-Dunford functional calculus for holomorphic functions. Its quaternionic counterpart is the $S$-functional calculus for slice hyperholomorphic functions, which is based on the $S$-spectrum and the $S$-resolvent \cite{Colombo:2010b}. This tool allowed to develop a theory of strongly continuous quaternionic groups and semigroups \cite{Colombo:2011b,Alpay:2014} and to generalize the $H^{\infty}$-functional calculus to the quaternionic setting \cite{Alpay:2016c}. The continuous functional calculus for normal operators on quaternionic Hilbert spaces was introduced in  \cite{Ghiloni:2013} and finally even the spectral theorem for quaternionic linear operators based on the $S$-spectrum was shown, for unitary operators in \cite{Alpay:2016a} and for unbounded normal operators in \cite{Alpay:2016} and later also in \cite{Ghiloni:2017}.

The theory is not only restricted to the quaternionic setting. Using a Clifford-algebra approach, Colombo, Sabadini and Struppa  defined the $S$-func\-tion\-al calculus  in \cite{Colombo:2008} also for $n$-tuples of non-commuting operators. Moreover, a theory of semigroups over real alternative $\ast$-algebras was introduced in \cite{Ghiloni:2016}. The $S$-resolvent is also crucial in slice hyperholomorphic Schur analysis, which was developed in \cite{Alpay:2015a,Alpay:2013} and the monograph \cite{Alpay:2017}. Finally, the point-$S$-spectrum coincides with the set of right eigenvalues such that the above concepts are suitable for applications in quaternionic quantum mechanics \cite{Adler:1995}. Fractional powers of quaternionic linear operators, which were introduced in \cite{fracpow}, might furthermore provide a tool for formulating new fractional diffusion problems \cite{HInftyComing,Marketing}.

Two overall themes drive the present article: the question of the minimal structure that is necessary to formulate the theory of quaternionic linear operators and the revelation of a beautiful and deep relation of quaternionic and complex operator theory. This relation shows that these two theories are compatible and  proves once more that the techniques associated to the $S$-spectrum are the correct approach towards operator theory in the non-commutative setting.

Quaternionic operator theory is usually formulated under the assumption that there exists a left multiplication on the considered Banach space. Indeed, this seemed a technical requirement since the formulas \eqref{SResSer} and \eqref{SRes1} for the $S$-resolvents contain the multiplication of operators with non-real scalars. For an operator $T$ and a scalar $a\in\hh$ the operators $aT$ and $Ta$ are supposed to act as $(aT)\bv = a(T\bv)$ and $(Ta)\bv = T(a\bv)$. This is only meaningful if a multiplication of vectors with scalars from the left is defined. Otherwise only the multiplication of operators with real scalars is defined, namely as $(aT)\bv = (Ta)\bv = T(\bv a)$. (This definition does however not yield a quaternionic linear operator if $a$ is not real!) 

A right linear operator is via the linearity-condition however only related to the right multiplication on the Banach space and hence its spectral properties should be independent of such a left multiplication. This could be considered a philosophical problem, but for instance on quaternionic Hilbert spaces only a right sided structure is defined a priori. In order to define the $S$-functional calculus for closed operators on quaternionic Hilbert spaces, one has to choose a random left multiplication in order to turn it into a two-sided space. The spectral properties of the operator should then be independent of this left-multiplication. Similarly, the proofs of the spectral theorem introduce left multiplications on the considered quaternionic Hilbert space (in \cite{Alpay:2016} it is only partially defined, in \cite{Ghiloni:2017} it is defined for all quaternions). This left multiplication is chosen to suite the considered operator in a certain sense, but it is only partially determined by the operator and then extended randomly. As the authors point out, the spectral measure associated with a normal operator is however independent of the chosen extension. 

The $S$-spectrum of a right linear operator is independent of the left multiplication as it is determined by the operator $\Q_{s}(T) = T^2 - 2s_0 T + |s|^2\id$, which only contains real scalars. In this article we show that the essential results in quaternionic operator theory can be formulated without the assumption of a left multiplication. Precisely, we show that the $S$-functional calculus can, for intrinsic slice hyperholomorphic functions, be expressed using only the right-multiplication on the vector space and we introduce an approach to spectral integration in the quaternionic setting that does not require the introduction of a left multiplication and has a clear interpretation in terms of the right multiplication on the considered vector space.

Any quaternion $x$ can be written as $x = x_0 + \I_{x} x_1$, where $x_0$ is the real part of $x$,  $x_1$ is the 
modulus of the imaginary part of $x$ and $\I_{x} = \underline{x}/|x|$ is the normalized imaginary part of $x$, which satisfies with $\I_{x}^{2} = -1$.  A function is called (left) slice function if it is of the form $f(x) = \alpha(x_0,x_1) + \I_{x}\beta(x_0,x_1)$ such that $\alpha(x_0,-x_1) = \alpha(x_0,x_1)$ and $\beta(x_0,-x_1) = -\beta(x_0,x_1)$. If $\alpha$ and $\beta$ satisfy the Cauchy-Riemann-equations, then $f$ is called left slice hyperholomorphic and if $\alpha$ and $\beta$ are real-valued then it is called intrinsic. Given a bounded operator $T$ and a left slice hyperholomorphic function $f$ that is defined on a suitable set $U$ with $\sigma_S(T)\subset U$, the operator $f(T)$ can be defined via the $S$-functional calculus as
\[
f(T) := \frac{1}{2\pi}\int_{\partial(U\cap\cc_{\I})}S_L^{-1}(s,T)\,ds_{\I}\,f(s)
\]
where $\I \in\SS := \{ x\in\hh: \Re(x) = 0, |x| = 1\} = \{x\in\hh: x^2 = -1\}$ is arbitrary, where $\cc_{\I}:= \{x_0 + \I x_1:x_0,x_1\in\rr\}$ and where $ds_{\I}:= ds (-\I)$. This corresponds to writing $f(x)$ in terms of the slice hyperholomorphic Cauchy formula and then formally replacing the scalar variable $x$ by the operator $T$, which clearly follows the idea of the Riesz-Dunford-functional calculus. We show that the symmetry $f(\overline{x}) = \overline{f(x)}$ satisfied by intrinsic slice functions allows us to rewrite this integral for intrinsic $f$ as
\begin{equation}\label{ZuMsti1}
f(T)\bv := \frac{1}{2\pi} \int_{\partial(U\cap\cc_{\I})} \RRes_{z}(T;\bv) f(z)\, dz\frac{1}{2\pi \I}\qquad \forall \bv\in V
\end{equation}
with
\begin{equation}\label{ZuMsti}
\RRes_{z}(T;\bv):= \Q_{z}(T)^{-1}(\bv z - T\bv),
\end{equation}
if $f$ is intrinsic. Observe that the expression \eqref{ZuMsti1} does not contain any multiplication of $\bv$ with non-real scalars from the left. It is hence independent from any left multiplication on the space $V$ and we can even use it to define the $S$-functional calculus on Banach spaces that are only right sided. 

The symmetry of the set of eigenvalues and of the $S$-spectrum explained in \eqref{RightEVRel} and \eqref{EVSym} is essential in order to properly define spectral integrals in the quaternionic setting. The equivalence classes of eigenvalues in \eqref{EVSym} are 2-spheres of the form $[s] = \{s_0 + \I s_1 : \I \in\SS\}$ and a set $U$ compatible with this symmetry (i.e. such that $[s]\subset U$ for all $s\in U$) is called axially symmetric. Observe that intrinsic slice functions of the form $f(s) = \alpha(s_0,s_1) + \I_{s}\beta(s_0,s_1)$ are compatible with this symmetry as $[f(s)] = f([s])$.

 A quaternionic spectral measure is a set function that associates a projection onto a right-linear subspace of $V$ to any $\Delta$ in the axially symmetric Borel sets $\sBorel(\hh)$ of $\hh$, for instance the generalized eigenspace $\ker\Q_{s}(T)$ of $T$ to the sphere $[s]$. It does this respecting the usual axioms, in particular $E(\Delta_1\cap\Delta_2) = E(\Delta_1)E(\Delta_2)$ and $E\left(\bigcup_{n\in\nn}\Delta_n\right) = \sum_{n=1}^{+\infty}E(\Delta_n)$ for $\Delta_n\cap\Delta_m = \emptyset$. For any bounded real-valued $\sBorel(\hh)$-measurable function $f$, we can define the spectral integral of $f$ with respect to $E$ as usual as
\begin{equation}\label{ZJUZJU}
\int_{\hh}f(s)\,dE(s) = \lim_{n\to+\infty}\int_{\hh}f_{n}(s)\,dE(s) = \lim_{n\to+\infty} \sum_{\ell = 1}^{N_n} a_{n,\ell}E(\Delta_n),
\end{equation}
where $f_{n}(s) = \sum_{\ell = 1}^{N_n} a_{n,\ell}\chi_{\Delta_n}(s)$ is a sequence of real-valued $\sBorel(\hh)$-measurable simple functions that converges uniformly to $f$. Observe that this integral is defined on spaces that are only right-sided: the coefficients  $a_{n,\ell}$ are real so that no multiplication of operators with non-real scalars appears. 

The class of quaternion-valued $\sBorel(\hh)$-$\sBorel(\hh)$-measurable functions is exactly the class of intrinsic slice functions on $\hh$ that are measurable with respect to the usual Borel sets. In order to integrate such functions, we introduce the concept of a spectral system $(E,J)$ consisting of a spectral measure $E$ and a spectral orientation $J$, which is a bounded operator satisfying $\ker J = \ran E(\rr)$ and $- J^2 =  E(\hh\setminus\rr)$. On $\ran J = \ran E(\hh\setminus\rr)$, the operator $J$ is essentially a right linear multiplication with values in $\SS$ from the right. For a bounded measurable intrinsic slice function $f = \alpha + \I \beta$, we then define
\[
 \int_{\hh} f(s)\,dE_{J}(s)  = \int_{\hh}\alpha(s_0,s_1)\,dE(s) + J \int_{\hh}\beta(s_0,s_1)\,dE(s),
\]
where $ \int_{\hh}\alpha(s_0,s_1)\,dE(s)$ and $\int_{\hh}\beta(s_0,s_1)\,dE(s)$ are defined in the sense of \eqref{ZJUZJU}. The idea of our approach to spectral integration is the following: while the spectral measure tells us how to associate subspaces of $V$ to sets of spectral spheres of the form $[s] = s_0 + \SS s_1$, the spectral orientation tells us how to multiply the different spectral values $s_{\I} = s_0 + \I s_1, \I\in\SS$, onto the vectors in these  subspaces.

The spectral theorem can  be formulated with our approach to spectral integration. In contrast to \cite{Alpay:2016} and \cite{Ghiloni:2017} our approach  does however not require to randomly introduce any undetermined structure: the operator fully determines its  spectral system. The proof of the spectral theorem by Alpay, Colombo and Kimsey in \cite{Alpay:2016} translates easily into the language of spectral systems. Actually this requires just a change of notation: they only use spectral integrals of intrinsic slice functions and for such functions the different approaches to spectral integration are equivalent if things are interpreted correctly.

The price that we have to pay for not introducing any undetermined structure and having a clear and intuitive interpretation of spectral integrals in terms of the right linear structure on the space is  that we can only integrate intrinsic slice functions. This class is much smaller than the class that in particular the approach using iqPVM in  \cite{Ghiloni:2017} allows to integrate. In \Cref{Conclusions} we however argue that extending the class of integrable functions beyond the class of intrinsic slice functions is not meaningful in the quaternionic setting. The intuition that underlies spectral integration is lost if this is done.

Although we are able to develop the essential results of quaternionic operator in this paper without assuming the existence of a left multiplication, the situation might be different for certain results that assume that the operator has commuting components. For instance in the characterization of the $S$-spectrum using the kernel of the $\mathcal{SC}$-functional calculus in \cite{Colombo:2012a} or the definition of $\mathcal{F}$-functional calculus in \cite{Colombo:2010c} make this assumption. For these results the left multiplication is essential, as it determines the components of the operator. Via the commutativity condition it is hence deeply related with the assumptions on the operator and so these results can not be independent of the left multiplication.
 
Using the tools described above, we develop then a theory of bounded spectral operators on quaternionic right linear Banach spaces, analogue to the complex theory in \cite{Dunford:1958}. A bounded complex linear spectral operator is a bounded complex linear operator that has a spectral resolution, i.e. there exists a spectral measure $E$ on the Borel sets $\Borel(\cc)$ of $\cc$  such that $E(\Delta) T = TE(\Delta)$ and $\sigma(T_{\Delta}) \subset \clos{\Delta}$ with $T_{\Delta} = T|_{\ran E(\Delta)}$ for any $\Delta\in\Borel(\cc)$, where $\clos{\Delta}$ denotes the closure of the set $\Delta$. As explained above, spectral measures are replaced by spectral systems in the quaternionic setting. A bounded quaternionic spectral operator is hence a right linear operator $T$ that has a spectral decomposition, i.e. a spectral system $(E,J)$ that commutes with $T$ such that $\sigma_{S}(T_{\Delta}) \subset\clos{\Delta}$ with $T_{\Delta} = T|_{\ran E(\Delta)}$ for all $\Delta\in\sBorel(\hh)$ and such that 
$ (s_0 \id - J s_1 - T)|_{\ran E(\hh\setminus\rr)}$ has a bounded inverse on $\ran E(\hh\setminus\rr)$. For such operator we generalize several known results from the complex setting. In particular any quaternionic spectral operator has a canonical decomposition $ T = S+ N$ into a scalar operator $S = \int_{\hh} s\,dE_J(s)$ and a quasi-nilpotent operator $N = T-S$. Furthermore, the operator $f(T)$ obtained via the $S$-functional calculus can be expressed as
\[
f(T) = \sum_{n=0}^{+\infty}N^n\frac{1}{n!} \int_{\hh}(\partial_S^nf)(s)\,dE_J(s),
\]
where $\partial_S^nf$ denotes the $n$-th slice derivate of $f$. Hence, $f(T)$ depends only on the values of $f$ on the $S$-spectrum  $\sigma_S(T)$ and not on its values on an entire neighborhood of $\sigma_S(T)$. Moreover, $f(T)$ is again a spectral operator, the spectral decomposition of which can be constructed from the spectral decomposition $(E,J)$ of $T$.

As the second theme of this article we mentioned a deep relation between complex and quaternionic operator theory. Let $V$ be a quaternionic right Banach space. If we choose $\I\in\SS$ and restrict the (right) scalar multiplication on $V$ to the complex plane $\cc_{\I}$, then we obtain a complex Banach space over $\cc_{\I}$, which we shall denote by $V_{\I}$. Any quaternionic linear operator on $V$ is then also a complex linear operator on $V_{\I}$. We show that the spectral properties of $T$ as a quaternionic linear operator on $V$ and the spectral properties of $T$ as a $\cc_{\I}$-linear operator on $V_{\I}$ are compatible. The spectrum $\sigma_{\cc_{\I}}(T)$ of $T$ as an operator on $V_{\I}$ equals $\sigma_{\cc_{\I}}(T) = \sigma_{S}(T) \cap\cc_{\I}$. Moreover, the $\cc_{\I}$-linear resolvent of $T$ at $z$ applied to $\bv$ is exactly the operator $\RRes_{z}(T;\bv)$ in \eqref{ZuMsti}. Hence the $S$-functional calculus for intrinsic slice hyperholomorphic functions on $V$ coincides because of \eqref{ZuMsti1} with the Riesz-Dunford functional calculus of $T$ on $V_{\I}$. Similarly, we can show that a spectral system $(E,J)$ generates a spectral measure  on  $V_{\I}$ over the complex plane $\cc_{\I}$. Integrating an intrinsic slice function $f$ with respect to $(E,J)$ is then equivalent to integrating the restriction of $f$ to $\cc_{\I}$ with respect to $E_{\I}$. Finally, a bounded quaternionic right linear operator $T$ on $V$ is a quaternionic spectral operator if and only if it is a complex spectral operator on $V_{\I}$. The spectral resolution of $T$ as a $\cc_{\I}$-complex spectral operator is exactly the spectral measure on $V_{\I}$ generated by its quaternionic spectral decomposition $(E,J)$. This correspondence of quaternionic and complex operator theory is due to the axial symmetry of the $S$-spectrum and compatibility of intrinsic slice functions with axial symmetry. The same reasons also guarantee that the choice of the imaginary unit $\I$ is irrelevant: this correspondence holds for any imaginary unit.

The paper is structured as follows: in \Cref{PrelimSect} we recall the most important facts of quaternionic operator theory and the related theory of slice hyperholomorphic functions. We do this in considerable detail (in particular for the spectral theorem) since this material is nonstandard and essential for the comparison of the different approaches to spectral integration. \Cref{IntCalcSect} contains the formulation of the $S$-functional calculus for intrinsic slice hyperholomorphic functions in terms of the right linear structure and the relation between the complex and the quaternionic notions of spectrum and resolvent. In \Cref{SpecIntSect} we develop our approach to spectral integration and compare it to the existing ones in \cite{Alpay:2016,Ghiloni:2017}. In \Cref{SpecOpSect} we introduce quaternionic spectral operators, we show that their spectral decompositions are unique and we show the relation to the complex theory of spectral operators. In \Cref{SpecOpSect2} we show that any quaternionic spectral operator has a canonical decomposition into a scalar and a radical part and we discuss the behavior of spectral operators and their spectral decompositions under the intrinsic $S$-functional calculus. \Cref{Conclusions} finally contains several general remarks on quaternionic operator theory. In particular we show that intrinsic slice functions are sufficient to reveal all the spectral information about a quaternionic linear operator.

\section{Preliminaries}\label{PrelimSect}
The skew-field of quaternions consists of the real vector space 
\[
\hh:=\left\{x = \xi_0 + \sum_{\ell =1}^3\xi_\ell e_\ell : \xi_\ell \in\rr\right\},
\]
 which is endowed with an associative product with unity $1$ such that
$e_\ell^2 = - 1$ and $e_\ell e_\kappa = -e_\kappa e_\ell$ for $\ell,\kappa\in\{1,2,3\}$ with $\ell \neq \kappa$.
The real part of a quaternion $x = \xi_0 + \sum_{\ell=1}^3\xi_\ell e_\ell$ is defined as $\Re(x) := \xi_0$, its imaginary part as $\underline{x} := \sum_{\ell=1}^3\xi_\ell e_\ell$ and its conjugate as $\overline{x} := \Re(x) - \underline{x}$. Furthermore, the modulus of $x$ is given by $|x|^2 = \overline{x}x = x\overline{x} = \sum_{\ell=0}^3 \xi_\ell^2$ and hence its inverse is $x^{-1} = \overline{x}|x|^{-2}$.

We denote the sphere of all normalized purely imaginary quaternions by
\[\SS := \{ x\in\hh: \Re(x) = 0, |x| = 1 \}.\]
If  $\I\in\SS$ then $\I^2 = -1$ and $\I$ therefore called an imaginary unit. The subspace $\cc_{\I} := \{x_0 + \I x_1: x_0,x_1\in\rr\}$ is then an isomorphic copy of the field of complex numbers. We moreover introduce the notation $\cc_{\I}^{+} = \{ x_0 + \I x_1: x_0\in\rr, x_1 >0\}$  for the open upper half plane,  $\cc_{\I}^{-} = \{ x_0 + \I x_1: x_0\in\rr, x_1 <0\}$  for the open lower half plane, and $\cc_{\I}^{\geq} = \{x_0 + \I x_1: x_0\in\rr, x_1 \geq 0\}$ for the closed upper halfplane in $\cc_{\I}$.

If $\I,\J \in\SS$ with $\I\perp \J$ and we set $\K=\I\J = -\J\I$, then $1$, $\I$, $\J$ and $\K$ form an
 orthonormal basis of $\hh$ as a real vector space and $1$ and $\J$ form a basis of $\hh$ as a left or right vector space over the complex plane $\cc_{\I}$, that is
 \[ \hh = \cc_{\I} + \cc_{\I} \J \quad\text{and}\quad \hh = \cc_{\I} + \J\cc_{\I}.\]
 Any quaternion $x$ belongs to a complex plane $\cc_{\I}$: if we set
 \[\I_x := \begin{cases}\frac{1}{|\underline{x}|}\underline{x},& \text{if  }\underline{x} \neq 0 \\ \text{any }\I\in\SS, \quad&\text{if }\underline{x}  = 0,\end{cases}\]
 then $x = x_0 + \I_x x_1$ with $x_0 =\Re(x)$ and $x_1 = |\underline{x}|$. (In several situations we shall also set $\I_{x} = 0$ if $x\in\rr$, which is sometimes more convenient.) The set
 \[
 [x] := \{x_0 + \I x_1: \I\in\SS\},
 \]
is a 2-sphere, that reduces to a single point if $x$ is real. Quaternions that belong to the same sphere can be transformed into each other by multiplication with a third quaternion.
\begin{corollary}\label{TurnCor}
Let $x\in\hh$. A quaternion $y\in\hh$ belongs of $[x]$ if and only if there exists $h\in\hh\setminus\{0\}$ such that $x = h^{-1}yh$.
\end{corollary}

\subsection{Slice hyperholomorphic functions} 
The theory of complex linear operators is based on the theory of holomorphic functions. Similarly, the theory of quaternionic linear operators is based on the theory of slice hyperholomorphic functions. We briefly recall the main properties of this class of functions. The corresponding proofs  can be found in \cite{Colombo:2011}.

\begin{definition}
A set $U\subset\hh$ is called
\begin{enumerate}[label = (\roman*)]
\item axially symmetric if $[x]\subset U$ for any $x\in U$ and
\item a slice domain if $U$ is open, $U\cap\rr\neq 0$ and $U\cap\cc_{\I}$ is a domain for any ${\I}\in\SS$.
\end{enumerate}
\end{definition}
In order to avoid confusion, we shall throughout this paper denote the closure of a set $U$ by $\clos{U}$ and its conjugation by $ \overline{U} = \{ \overline{x}:x\in U\}$.
\begin{definition}\label{sHolDef}
Let $U\subset\hh$ be axially symmetric. A function $f: U \to \hh$ is called left slice function, if it is of the form
\begin{equation}
\label{lHolDef}f(x) = \alpha(x_0,x_1) + \I_x \beta(x_0, x_1) \quad \forall x\in U,
\end{equation}
where $\alpha$ and $\beta$ are functions that take values in $\hh$ and satisfy the compatibility condition
\begin{equation}\label{CCond}
\alpha(x_0,x_1) = \alpha(x_0,-x_1)\quad\text{and}\quad \beta(x_0,x_1) = -\beta(x_0,-x_1).
\end{equation}
If in addition $U$ is open and $\alpha$ and $\beta$ satisfy the Cauchy-Riemann-differential equations
\begin{equation}\label{CR}
\frac{\partial}{\partial x_0} \alpha(x_0,x_1) = \frac{\partial}{\partial x_1} \beta(x_0,x_1)\quad\text{and}\quad \frac{\partial}{\partial x_1}\alpha(x_0,x_1) = - \frac{\partial}{\partial x_0} \beta(x_0,x_1),
\end{equation}
then $f$ is called left slice hyperholomorphic.
A function $f: U \to \hh$ is called  right slice function, if  it is of the form
\begin{equation}\label{rHolDef}
f(x) = \alpha(x_0,x_1) +\beta(x_0,x_1)  \I_x \quad \forall x \in U,
\end{equation}
with functions $\alpha$ and $\beta$ satisfying \eqref{CCond}. If in addition $U$ is open and $\alpha$ and $\beta $ satisfy \eqref{CR}, then $f$ is called right slice hyperholomorphic. 

Finally, a left slice function $f = \alpha + \I\beta$ is called  intrinsic slice function if $\alpha$ and $\beta$ take values in $\rr$.

We denote the set of all left slice functions on $U$ by $\SliceL(U)$, the set of all right slice functions on $U$ by $\SliceR(U)$ and the set of all intrinsic slice functions on $U$ by $\SliceReal(U)$. Moreover, we denote by $\lhol(U)$ (or $\rhol(U)$ resp. $\intrin(U)$) the set of all functions that are left slice hyperholomorphic (or right slice hyperholomorphic resp. intrinsic slice hyperholomorphic) on an open axially symmetric set $U'$ with $U\subset U'$.
\end{definition}
\begin{remark}
Observe that any quaternion $x$ can be represented using two different  imaginary units, namely $x = x_0 + \I_x x_1 = x_0 + (-\I_x)(-x_1)$. If $x\in\rr$, then we can even choose any imaginary unit we want in this representation. The compatibility condition \eqref{CCond} assures that the choice of this imaginary unit is irrelevant.
\end{remark}

Any intrinsic slice function is both a left and a right slice function because $\I$ and $\beta$ commute. The converse is however not true: the constant function $f \equiv c \in \hh\notin \rr$ is obviously a left and a right slice function, but not intrinsic. Intrinsic slice functions can be characterized in several ways.
\begin{corollary}\label{RealSliceChar}
If $f\in\SliceL(U)$ or $f\in\SliceR(U)$, then the following statements are equivalent.
\begin{enumerate}[(i)]
\item The function $f$ is an intrinsic slice function.
\item \label{COnj} We have $f(\overline{x}) = \overline{f(x)}$ for any $x\in U$. 
\item We have $f(U\cap\cc_{\I}) \subset \cc_{\I}$ for all $\I\in\SS$.
\end{enumerate}
\end{corollary}

The importance of this subclass is due to the fact that the pointwise multiplication and the composition with intrinsic slice functions preserve the slice structure. This is not true for arbitrary slice functions. Moreover, if the function is also slice hyperholomorphic, then slice hyperholomorphicity is preserved too. 
\begin{corollary}Let $U\subset\hh$ be axially symmetric.
\begin{enumerate}[(i)]
\item If $f\in\SliceReal(U)$ and $g\in\SliceL(U)$, then $fg\in\SliceL(U)$. If $f\in\SliceR(U)$ and $g\in\SliceReal(U)$, then $fg\in\SliceR(U)$.
\item If $f\in\intrin(U)$ and $g\in\lhol(U)$, then $fg\in\lhol(U)$. If $f\in\rhol(U)$ and $g\in\intrin(U)$, then $fg\in\rhol(U)$.
\item If $g\in\SliceReal(U)$ and $f\in\SliceL(g(U))$, then $f\circ g\in \SliceL(U)$. If $g\in\SliceReal(U)$ and $f\in\SliceR(g(U))$, then $f\circ g\in \SliceR(U)$.
\item If $g\in\intrin(U)$ and $f\in\lhol(g(U))$, then $f\circ g\in \lhol(U)$. If $g\in\intrin(U)$ and $f\in\rhol(g(U))$, then $f\circ g\in \rhol(U)$.
 \end{enumerate}
\end{corollary}

The values of a slice  function are uniquely determined by its values on an arbitrary complex plane $\cc_{\I}$. \begin{theorem}[Representation Formula]\label{RepFo}
Let $U\subset \hh$ be axially symmetric  and let $\I\in\SS$. For any $x = x_0 + \I_x x_1\in U$ set $x_{\I} := x_0 + \I x_1$. If $f\in\SliceL(U)$, then
\begin{equation}\label{RepFoL}
f(x) = \frac12(1-\I_x\I)f(x_{\I}) + \frac12(1+\I_x\I)f(\overline{x_{\I}}) \quad\text{for all $x\in U.$}
\end{equation}
If $f\in\rhol(U)$, then
\[f(x) = f(x_{\I})(1-\I\I_x)\frac12 + f(\overline{x_{\I}})(1+\I\I_x x)\frac12 \quad\text{for all $x\in U$.}\]
\end{theorem}

As a  consequence, any quaternion-valued function on a suitable subset of a complex plane possesses a unique slice extension.
\begin{corollary}\label{extLem}
Let $\I\in\SS$ and let $f:O\to\hh$, where $O$ is a set in $\cc_{\I}$ that is symmetric with respect to the real axis.
We define the axially symmetric hull of  $O$ as $[O]: = \bigcup_{z\in O}[z]$.
\begin{enumerate}[label = (\roman*)]
\item There exists a unique function $f_L \in \SliceL([O]) $ such that $f_L |_{O\cap\cc_{\I}} = f$. Similarly, there exists a unique function $f_R \in \SliceR(O) $ such that $f_R |_{O\cap\cc_{\I}} = f$.
\item If $f$ satisfies $\frac{1}{2}\left(\frac{\partial}{\partial x_0} f(x) + \I\frac{\partial}{\partial x_1} f(x)\right) = 0$, then  $f_L$ is left slice hyperholomorphic.
\item If $f$ satisfies $ \frac{1}{2}\left(\frac{\partial}{\partial x_0} f(x) +\frac{\partial}{\partial x_1} f(x) \I \right)   = 0$, then $f_{R}$ is right slice hyperholomorphic.
\end{enumerate}
\end{corollary}
\begin{remark}\label{ExtensionOfRealSlicefFs}
If $O \subset \cc_{\I}^{\geq}$ and $f: O \mapsto \cc_{\I}$ is such that $f(O\cap\rr)\subset\rr$, then there exists a unique intrinsic slice extension $\tilde{f}\in\SliceReal([O])$ of $f$ to $[O]$. Indeed, we can use the relation \cref{COnj} in  \Cref{RealSliceChar} to extend $f$ to $O \cup \overline{O}$ by setting $f(\overline{z}) = \overline{f(z)}$. The set  $O\cup\overline{O}$ is symmetric with respect to the real axis and hence \Cref{extLem} implies the existence of a left slice extension $f_L$ of $f$. If we write $f(z) = \alpha(z_0,z_1) + \I \beta(z_0,z_1)$ with $\alpha(z), \beta(z)\in \rr$ for $z\in O$, this slice extension is because of \eqref{RepFoL}  given by 
\[
f_L(x) = \frac12 \left(f(x_{\I}) + \overline{f(x_{\I})}\right) + \I_x(-\I)\frac12 \left(f(x_{\I}) - \overline{f(x_{\I})}\right) = \alpha(x_0,x_1) + \I_{x}\beta(x_0,x_1)
\]
Hence, it is intrinsic. Any intrinsic slice function is therefore already entirely determined by its values on a complex halfplane $\cc_{\I}^{\geq}$. 
\end{remark}
Let us now turn our attention to the generalizations of classical function theoretic results. Important examples of slice hyperholomorphic functions are power series with quaternionic coefficients: series of the form $\sum_{n=0}^{+\infty}x^na_n$ are left slice hyperholomorphic and series of the form $\sum_{n=0}^{\infty} a_nx^n$ are right slice hyperholomorphic on their domain of convergence. A power series is intrinsic if and only if its coefficients are real. Conversely, any slice hyperholomorphic function can be expanded into a power series at any real point. 

\begin{definition}
Let $U\subset \hh$ be an axially symmetric open set. For any $f\in\lhol(U)$, the function
\[
\sderiv f(x) = \lim_{\cc_{\I_{}}\ni s\to x} (s - x)^{-1} (f(s)- f(x)),
\]
where $\lim_{\cc_{\I_x}\ni s \to x} g(s)$ denotes the limit of $g$ as $s$ tends to $x$ in $\cc_{\I_{x}}$, is called the slice derivative of $f$. Similarly, if $f\in\rhol(U)$, then the function 
\[
\sderiv f(x) = \lim_{\cc_{\I_{x}}\ni s\to x} (f(s)- f(x)) (s - x)^{-1}  
\]
is called the slice derivative of $f$. 
\end{definition}
\begin{remark}
The slice derivative of a left or right slice hyperholomorphic function is again left resp. right slice hyperholomorphic. Moreover, it coincides with the partial derivative $\frac{\partial}{\partial x_0} f(x)$ of $f$ with respect to the real part $x_0$ of $x$. It is therefore also at  points $x\in\rr$ well defined and independent of the choice of $\I_{x}$.
\end{remark}

\begin{theorem}\label{Taylor}
If $f$ is left slice hyperholomorphic on the ball $B_r(a)$ of radius $r>0$ centered at $a\in\rr$, then
\[f(x) = \sum_{n=0}^{+\infty} (x-a)^n \frac{1}{n!}\left( \sderiv^n f\right)(a)\quad\text{for $x\in B_r(a)$.}\]
If $f$ is right slice hyperholomorphic on $B_r(a)$, then
\[f(x) = \sum_{n=0}^{+\infty}\frac{1}{n!}\left( \sderiv^n f\right)(a)(x-a)^n \quad\text{for $x\in B_r(\alpha)$.}\]
\end{theorem}

If we restrict a slice hyperholomorphic function to one slice $\cc_{\I}$, then we obtain a vector-valued function that is holomorphic in the usual sense.
\begin{lemma}\label{HolLem}
Let $U\subset\hh$ be an axially symmetric open set. If $f\in\lhol(U)$, then for any $\I\in\SS$ the restriction $f_{\I} := f|_{U\cap\cc_{\I}}$ is left holomorphic, i. e.
\begin{equation}\label{SplitEQL}
 \frac{1}{2}\left(\frac{\partial}{\partial x_0} f_{\I}(x) + \I\frac{\partial}{\partial x_1} f_{\I}(x)\right) = 0,\qquad \forall x = x_0 + \I x_1 \in U\cap\cc_{\I}.
 \end{equation}
If $f\in\rhol(U)$, then for any $\I\in\SS$ the restriction $f_{\I} := f|_{U\cap\cc_{\I}}$ is right holomorphic, i. e.
\begin{equation}\label{SplitEQR}
\frac{1}{2}\left(\frac{\partial}{\partial x_0} f_{\I}(x) + \frac{\partial}{\partial x_1} f_{\I}(x)\I\right) = 0,\qquad \forall x = x_0 + \I x_1 \in U\cap\cc_{\I}.
\end{equation}
\end{lemma}

Finally, slice hy\-per\-ho\-lo\-mor\-phic functions satisfy an adapted version of the Cauchy integral theorem and an integral formula of Cauchy-type with a modified kernel.
\begin{definition}
We define the  left slice hyperholomorphic Cauchy kernel as
\[S_L^{-1}(s,x) = (x^2-2x_0x + |s|^2)^{-1}(\overline{s} - x)\quad\text{for }x\notin[s]\]
and the right slice hyperholomorphic Cauchy kernel as
\[S_R^{-1}(s,x) = (\overline{s} - x)(x^2-2s_0x + |s|^2)^{-1}\quad\text{for }x\notin[s].\]
\end{definition}

\begin{corollary}
The left slice hyperholomorphic Cauchy-kernel $S_L^{-1}(s,x)$ is left slice hyperholomorphic in the variable $x$ and right slice hy\-per\-ho\-lo\-mor\-phic in the variable $s$ on its domain of definition. Moreover, we have $S_R^{-1}(s,x) = - S_L^{-1}(x,s)$.
\end{corollary}
\begin{remark}
If $x$ and $s$ belong to the same complex plane, they commute and the slice hyperholomorphic Cauchy-kernels reduce to the classical one, i.e.
\[ (s-x)^{-1} = S_L^{-1}(s,x) = S_R^{-1}(s,x).\]
\end{remark}
\begin{theorem}[Cauchy's integral theorem]
Let $U\subset\hh$ be an axially symmetric open set, let $\I\in\SS$ and let $D_{\I}$ be a bounded open subset of $O\cap\cc_{\I}$ with $\clos{D_{\I}}\subset O\cap\cc_{\I}$ such that its boundary consists of a finite number of continuously differentiable Jordan curves. For any $f\in\rhol(U)$ and $g\in\lhol(U)$, it is
\[\int_{\partial D_{\I}}f(s)\,ds_{\I}\,g(s) = 0,\]
where $ds_{\I} = -{\I}\ ds$.
\end{theorem}
\begin{definition}
An axially symmetric open set $U\subset \hh$ is called a slice Cauchy domain if $U\cap\cc_{\I}$ is a Cauchy domain for any ${\I}\in\SS$, that is
\begin{enumerate}[(i)]
\item $U\cap\cc_{\I}$ is open
\item $U\cap\cc_{\I}$ has a finite number of components (i.e. maximal connected subsets), the closures of any two of which are disjoint
\item the boundary of $U\cap\cc_{\I}$ consists of a finite positive number of closed piecewise continuously differentiable Jordan curves.
\end{enumerate}
\end{definition}

\begin{remark}
A slice Cauchy domain is either bounded or has exactly one unbounded component. If it is unbounded, then the unbounded component contains a neighborhood of infinity.
\end{remark}

\begin{theorem}[Cauchy's integral formula]\label{Cauchy}
Let $U\subset\hh$ be a slice Cauchy domain, let $\I\in\SS$ and set $ds_{\I} = -{\I}\, ds$. If $f$ is left slice hyperholomorphic on an open set that contains $\clos{U}$, then
\[f(x) = \frac{1}{2\pi}\int_{\partial(U\cap\cc_{\I})} S_L^{-1}(s,x)\,ds_{\I}\, f(s)\quad\text{for all }x\in U.\]
If $f$ is right slice hyperholomorphic on an open set that contains $\clos{U}$, then
\[f(x) = \frac{1}{2\pi}\int_{\partial(U\cap\cc_{\I})}f(s)\, ds_{\I}\, S_R^{-1}(s,x)\quad\text{for all }x\in U.\]
\end{theorem}

%
%

\begin{remark}\label{OpVal}
The results presented in this section can  be extended to functions with values in a two-sided quaternionic Banach space. Problems concerning vector-valued functions can be reduced to scalar problems by applying elements of the dual space, analogue to what is done in the complex setting \cite{Alpay:2015a}.
\end{remark}
\subsection{The S-functional calculus}
The natural extension of the Riesz-Dunford-functional calculus for complex linear operators  to quaternionic linear operators is the so-called $S$-functional calculus. It is based on the theory of slice hyperholomorphic functions and follows the idea of the classical case: to formally replace the scalar variable $x$ in the Cauchy formula by an operator. Unless stated differently, the proofs of the results in this subsection can be found in \cite{Colombo:2011}.

Let us start with a precise definition of the various structures of vector, Banach and Hilbert spaces in the quaternionic setting.
\begin{definition}
A quaternionic right vector space is an additive group $(V,+)$ endowed with a scalar multiplication from the right such that for all $\bu, \bv\in V$ and all $a,b\in\hh$
\begin{align}\label{RVecS} 
(\bu + \bv)a& = \bu a + \bv a &   \bu (a+b) &= \bu a+ \bu  b & \bv (ab) &= (\bv a)b & \bv 1 & = \bv.
\end{align}
A quaternionic left vector space is an additive group $(V,+)$ endowed with a scalar multiplication from the left such that for all $\bu,\bv\in V$ and all $a,b\in\hh$
\begin{align}\label{LVecS}
a(\bu+\bv) &= a\bu + a\bv & (a+b)\bv &= a\bv + b\bv & (ab)\bv &= a(b\bv) & 1\bv &= \bv.
\end{align}
Finally, a two sided quaternionic vector space is an additive group $(V,+)$ together with a scalar multiplication from the right and a scalar multiplication from the left that satisfy \eqref{RVecS} resp. \eqref{LVecS} such that in addition $a\bv = \bv a$ for all $\bv\in V$ and all $a\in\rr$.
\end{definition}
\begin{remark}
Starting from a real vector space $V_{\rr}$, we can easily construct a two-sided quaternionic vector space by setting 
\[ V_{\rr} \otimes \hh = \left\{ \sum_{\ell=0}^{3} \bv_{\ell}\otimes e_{\ell}: \bv_{\ell}\in V_{\rr}\right\},\]
where we denote $e_{0}=1$ for neatness. Together with the componentwise addition $V_{\rr}\otimes\hh$ forms an additive group, which is a two-sided quaternionic vector space, if we endow it with the right and left scalar multiplication
\[ a\bv = \sum_{\ell,\kappa=0}^3 (a_{\ell}\bv_{\kappa})\otimes (e_{\ell}e_{\kappa})\qquad\text{and}\qquad \bv a = \sum_{\ell,\kappa=0}^3 (a_{\ell}\bv_{\kappa})\otimes (e_{\kappa}e_{\ell})\]
for $a = \sum_{\ell=0}^3 a_{\ell}e_{\ell}\in\hh$ and $\bv = \sum_{\kappa=0}^{3} \bv_{\kappa}\otimes e_{\kappa}\in V_{\rr}\otimes\hh$. Usually one omits the symbol $\otimes$ and simply writes $\bv = \sum_{\ell=0}^3 \bv_{\ell}e_{\ell}$.

It was shown in \cite{Ng:2007} that any two-sided quaternionic vector-space $\bv$ is essentially of this form. Indeed, we can set
\begin{equation}\label{VR}
V_{\rr}  = \{\bv\in V: a\bv = \bv a\ \forall a\in\hh\}
\end{equation}
and find that $V$ is isomorphic to $V_{\rr}\otimes \hh$. If we set $\Re(\bv) := \frac{1}{4} \sum_{\ell = 0}^{3} \overline{ e_{\ell}} \bv e_{\ell}$, then $\Re(\bv) \in V_{\rr}$ and $\bv =  \sum_{\ell=0}^3 \Re(\overline{e_{\ell}} \bv )e_{\ell}$.
\end{remark}
\begin{remark}\label{RealVS}
A quaternionic right or left vector space also carries the structure of a real vector space: if we simply restrict the quaternionic scalar multiplication to $\rr$, then we obtain a real vector space. Similarly, if we choose some $\I\in\SS$ and identify $\cc_{\I}$ with the field of complex numbers, then $V$ also carries the structure of a complex vector space over $\cc_{\I}$, which we obtain by restricting the quaternionic scalar multiplication to $\cc_{\I}$.

If we consider a two-sided quaternionic vector space, then the left and the right scalar multiplication coincide for real numbers so that we can restrict them to $\rr$ in order to obtain again a real vector space. This is however not true for multiplication with scalars in one complex plane $\cc_{\I}$. In general $a\bv \neq \bv a$ for $a\in\cc_{\I}$ and $\bv\in V$. Hence, we can only restrict either the left or the right multiplication to $\cc_{\I}$ in order to consider $V$ as a complex vector space over $\cc_{\I}$, but not both simultaneously.
\end{remark}
\begin{definition}
A quaternionic right (left or two-sided) Banach space is a quaternionic right (left or two-sided) vector space endowed in a norm $\|\cdot\|$ in the sense of real vector spaces (cf. \Cref{RealVS}) that is compatible with the quaternionic right (left resp. two-sided) scalar multiplication, i.e. such that $\|\bv a\| = \| \bv \| |a|$ ( or $\|a\bv\| = | a|\|\bv\|$ resp. $\| a\bv \| = |a|\|\bv\| = \| \bv a\|$) for all $a\in\hh$ and all $\bv\in V$. 
\end{definition}
\begin{remark}\label{RealBS}
Similar to \Cref{RealVS}, we obtain a real or complex Banach space if we restrict the left or right scalar multiplication on a quaternionic Banach space to $\rr$ resp. $\cc_{\I}$ for some $\I\in\SS$.
\end{remark}
\begin{definition}
A quaternionic right Hilbert space $\hil$ is a quaternionic right vector space equipped with a scalar product $\langle\cdot,\cdot\rangle: \hil\times\hil \to \hh$ that satisfies $\langle \bu,\bv a + \bw\rangle = \langle \bu,\bv\rangle a + \langle \bu, \bw\rangle$, $\langle \bu, \bv\rangle = \overline{\langle\bv,\bu\rangle}$ and $\langle \bu,\bu \rangle \geq 0$ for all $\bu,\bv,\bw\in\hil$ and all $a\in\hh$ such that $\hil$ is complete with the norm $\|\bv\| = \sqrt{\langle \bu, \bu\rangle}$. 
\end{definition}
\begin{remark}\label{realHil}
Also a quaternionic Hilbert space carries natural real and complex Hilbert space structures: if we restrict the right scalar multiplication on $\hil$ to $\rr$ and define $\langle \bu, \bv\rangle_{\rr} := \Re \langle \bu,\bv\rangle$, then $(\hil,\langle\cdot,\cdot\rangle_{\rr})$ is a real Hilbert space. Similarly, if we choose $\I\in\SS$, then we can set $\langle \bu,\bv\rangle_{\I}:= \{\langle \bu,\bv\rangle\}_{\I}$, where $\{a\}_{\I}$ denotes the $\cc_{\I}$-part of a quaternion $a$, i.e. $\{a\}_{\I} := \eta_1 $ if $a = \eta_1 + \eta_2 \J$ with $\eta_1,\eta_2\in\cc_{\I}$ and $\J\in\SS$ with $\I\perp \J$. If we restrict the right scalar multiplication on $\hil$ to $\cc_{\I}$, then $(\hil,\langle.,.\rangle_{\I})$ is a complex Hilbert space over $\cc_{\I}$. Observe that for $\J\in\SS$ with $\J\perp\I$ and any $\bv$ in $\hil$, the vectors $\bv$ and $\bv\J$ are orthogonal in this structure.
\end{remark}
\begin{remark}
The fundamental concepts of complex Hilbert spaces such as orthogonality, orthonormal bases, the Riesz representation theorem etc. can be defined as in the complex case. For more detailed definitions we refer to \cite{Ghiloni:2013}.
\end{remark}
\begin{notation}
Since we are working with different number systems and vector space structures in this article, we introduce for a set of vectors $ \mathbf{B} := (\bb_{\ell})_{\ell\in\varLambda}$ the quaternionic right-linear span of $\mathbf{B}$ 
\[
\linspan{\hh} \mathbf{B} := \left\{ \sum_{\ell\in I} \bb_{\ell}a_{\ell}: a_{\ell}\in\hh, I\subset\varLambda \ \text{finite}\right\}
\]
and the $\cc_{\I}$-linear span of $\mathbf{B}$
\[
\linspan{\cc_{\I}} \mathbf{B} := \left\{ \sum_{\ell\in I} \bb_{\ell}a_{\ell}: a_{\ell}\in\cc_{\I}, I\subset\varLambda \ \text{finite}\right\}.
\]
\end{notation}
\begin{definition}
An operator $T:  V_{R}\to  V_{R}$ on a right-sided Banach space $ V_{R}$ is called right linear if $T(\bu a + \bv) = T(\bu) a + T(\bv)$ for all $\bu,\bv\in  V_{R}$ and $a\in\hh$. The operator $T$ is bounded if $\| T\| := \sup_{\|\bv \| = 1}\|T \bv\|$ is finite and closed if the graph of $T$ is closed in $ V_{R}\times V_{R}$. We denote the set of all bounded right linear operators on $ V_{R}$ by $\boundOP( V_{R})$ and the set of all closed operators on $V_{R}$ by $\closOP( V_{R})$. 
\end{definition}
\begin{remark}
The notions of adjoints of self-adjoint, unitary and normal operators etc. are in the quaternionic setting defined as in the complex one. We refer to \cite{Ghiloni:2013} for precise definitions.
\end{remark}

\begin{remark}\label{B(V)Struct}
If $ V_{R}$ is even a two-sided Banach space, then $\boundOP( V_{R})$ is again a two-sided Banach space with the scalar multiplications 
\begin{equation}\label{OPScalMul}
(aT)(\bv) = a(T(\bv))\quad \text{and}\quad (Ta)(\bv ) = T(a\bv )\quad\forall a\in\hh, \bv\in V_{R}.
\end{equation}
 If however $ V_{R}$ is just one-sided, then these formulas are not meaningful and $ V_{R}$ is only a real Banach space on which no multiplication with quaternionic scalars is defined.
 \end{remark}
Throughout the paper, let $ V_{R}$ denote a right-sided and let $V$ denote a two-sided quaternionic Banach space.
For $T\in\closOP( V_{R})$, we define
\[ \Q_{s}(T) := T^2 - 2 s_0 T + |s|^2\id, \qquad \text{for $s\in\hh$}. \]
\begin{definition}
Let $T\in\closOP( V_{R})$. We define the $S$-resolvent set of  $T$ as
\[
\rho_S(T):= \left\{ s\in\hh: \Q_{s}(T)^{-1}\in\boundOP(V_R)\right\}
\]
and the $S$-spectrum of $T$ as
\[\sigma_S(T):=\hh\setminus\rho_S(T).\]
If $s\in\rho_S(T)$, then $\Q_{s}(T)^{-1}$ is called the pseudo-resolvent (sometimes also the  spherical resolvent) of $T$ at $s$. Furthermore, we define the extended $S$-spectrum $\sigma_{SX}(T)$ as
\begin{equation*}
\sigma_{SX}(T) := \begin{cases} \sigma_{S}(T) & \text{if  $T$ is bounded,}\\
\sigma_{S}(T)\cup\{\infty\} & \text{if  $T$ is unbounded.}
\end{cases}
\end{equation*}
\end{definition}
The $S$-spectrum generalizes the set of right eigenvalues and has properties that are similar to those of the spectrum of a complex linear operator.
\begin{theorem}\label{SpecProp}
Let $T\in\closOP( V_{R})$.
\begin{enumerate}[label = (\roman*)]
\item The $S$-spectrum $\sigma_{S}(T)$ of $T$ is axially symmetric. It contains the set of right eigenvalues of $T$ and if $V_{R}$ has finite dimension, then it equals the set of right eigenvalues.
\item The $S$-spectrum $\sigma_{S}(T)$ is a closed subset of $\hh$ and the extended $S$-spectrum $\sigma_{SX}(T)$ is a closed and hence compact subset of $\hh_{\infty}:= \hh \cup\{\infty\}$.
\item If $T$ is bounded, then $\sigma_{S}(T)$ is nonempty and bounded by the norm of $T$, i.e. $\sigma_{S}(T)\subset \clos{B_{\|T\|}(0)}$. 
\end{enumerate}
\end{theorem}
\begin{remark}\label{Usually2Sided}
In quaternionic operator theory, one usually works on a two-sided Banach space $V$ because the $S$-resolvent operators (cf. \Cref{SResDef}) can only be defined if  a quaternionic multiplication on $\boundOP(V)$ is available. This is not the case if $V$ is only one-sided, cf. \Cref{B(V)Struct}.

The operator $\Q_{s}(T)$ on the other hand involves only real scalars and hence the $S$-spectrum $\sigma_S(T)$ of an operator $T$  can also be defined if $T$ acts on a one-sided Banach space. Although the properties of $\sigma_S(T)$ given in \Cref{SpecProp} are for the above reason usually proved under the assumption that $T$ acts on a two-sided Banach space (see for instance \cite{Colombo:2011}), they remain true for an operator on a one-sided Banach space. Indeed, their proofs do not make any use of a quaternionic multiplication on $\boundOP(V)$. One only has to pay attention when showing that $\sigma_{S}(T)$ is bounded by $\|T\|$. This is usually shown using a power series expansion of the $S$-resolvent that involves quaternionic coefficients, cf. \cite{Colombo:2011}. However,  \cite[Section~4]{Colombo:2016a} introduced a new series expansion for the pseudo-resolvent $\Q_{s}(T)^{-1}$ converging for $|s|>\|T\|$, namely 
\begin{equation}\label{ALPT}
\Q_{s}(T)^{-1} = \sum_{n=0}^{+\infty}T^na_n\quad\text{with}\quad a_n= |s|^{-2n-2}\sum_{k=0}^{n}\overline{s}^{k}s^{n-k}.
\end{equation}
 Since $\overline{a_n} = a_{n}$ this series contains only real coefficients and can hence also be defined on a one-sided Banach space. The computations in the proof of \cite[Lemma~4.2]{Colombo:2016a} that showed that $(\sum_{n=0}^{+\infty}T_na_n)\Q_{s}(T)^{-1} = \id$ require a quaternionic multiplication on $\boundOP(V)$ and can not be performed if $V$ is only one-sided. However, if we write $(\sum_{n=0}^{+\infty}T^na_n)\bv = \sum_{n=0}^{+\infty}T^n\bv a_n$ for any $\bv \in V$, then we can show that $\Q_{s}(T)^{-1}(\sum_{n=0}^{+\infty}T^na_n)\bv = \Q_{s}(T)^{-1}(\sum_{n=0}^{+\infty}T^n\bv a_n) = \bv$ by essentially the same computations, which now only make use of the right multiplication on $V$. The series expansion \eqref{ALPT} can therefore serve for showing that $\sigma_{S}(T)\subset \clos{B_{\|T\|}(0)}$ also for operators on one-sided Banach spaces such that all properties in \Cref{SpecProp} are also true in this setting. This is in particular important because we will mainly work on right-sided and not on two-sided quaternionic Banach spaces in this article.
\end{remark}
\begin{definition}\label{SResDef}
Let $T\in\closOP(V)$ be a closed operator on a two-sided quaternionic Banach space. For $s\in\rho_S(T)$, the left $S$-resolvent operator is defined as
\begin{equation}\label{SresolvoperatorL}
S_L^{-1}(s,T):= \Q_s(T)^{-1}\overline{s} -T\Q_s(T)^{-1}
\end{equation}
and the right $S$-resolvent operator is defined as
\begin{equation}\label{SresolvoperatorR}
S_R^{-1}(s,T):=-(T-\id \overline{s})\Q_s(T)^{-1}.
\end{equation}
\end{definition}
\begin{remark}\label{RkResExtension} One clearly obtains the right $S$-resolvent operator by formally replacing the variable $x$ in the right slice hyperholomorphic Cauchy kernel by the operator $T$. The same procedure yields
\begin{equation}\label{LResShort}
S_L^{-1}(s,T)\bv = -\Q_s(T)^{-1}(T-\overline{s}\id)\bv,\quad\text{for }\bv\in\dom(T)
\end{equation}
for the left $S$-resolvent operator. This operator is only defined  on  the domain $\dom(T)$ of $T$ and not on the entire space $V$. However,  $\Q_s(T)^{-1}T\bv = T\Q_{s}(T)^{-1}\bv$ for $\bv\in\dom(T)$ and commuting $T$ and $Q_s(T)$ in \eqref{LResShort} yields \eqref{SresolvoperatorL}. For arbitrary $s\in\hh$, the operator $\Q_{s}(T) = T^2 - 2s_0T + |s|^2\id$ maps $\dom(T^2)$ to $V$. Hence, the pseudo-resolvent $\Q_s(T)^{-1}$ maps $V$ to $\dom(T^2)\subset \dom(T)$ if $s\in\rho_S(T)$. Since $T$ is closed and $\Q_s(T)^{-1}$ is bounded, equation \eqref{SresolvoperatorL} then defines a continuous and therefore bounded right linear operator on the entire space $V$. Hence, the left resolvent $S_L^{-1}(s,T)$ is the closed extension of the operator \eqref{LResShort} to~$V$. In particular, if $T$ is bounded, then $S_L^{-1}(s,T)$ can directly be defined by \eqref{LResShort}.

If one considers left linear operators, then one must modify the definition
of the right $S$-resolvent operator for the same reasons.
\end{remark}
\begin{remark}
The $S$-resolvent operators reduce to the classical resolvent if $T$ and $s$ commute, that is
\[S_L^{-1}(s,T) = S_R^{-1}(s,T) = (s\id - T)^{-1}.\]
This is in particular the case if $s$ is real.
\end{remark}

The proof of the following crucial lemma  can be found for instance in \cite{Colombo:2011} for the case of bounded operators. For the case of unbounded operators, several additional technical difficulties have to be overcome. The respective proof can be found in \cite{fracpow}.
\begin{lemma}\label{ResHol2342}
Let $T\in\closOP(V)$. The map $s\mapsto S_L^{-1}(s,T)$ is a right slice-hyperholomorphic function on $\rho_S(T)$ with values in the two-sided  quaternionic Banach space $\boundOP(V)$. The map $s\mapsto S_R^{-1}(s,T)$ is a left slice-hyperholomorphic function on $\rho_S(T)$ with values in the two-sided quaternionic Banach space $\boundOP(V)$. 
\end{lemma}
The $S$-resolvent equation plays the role of the classical resolvent equation in the quaternionic setting and has been first introduced in \cite{Alpay:2015}.
\begin{theorem}[$S$-resolvent equation]Let $T\in\closOP(V)$. For  $s,x \in  \rho_S(T)$ with $s\notin[x]$, it is
\begin{equation}\label{SresEQ1}
\begin{split}
S_R^{-1}(s,T)S_L^{-1}(x,T)=&\big[[S_R^{-1}(s,T)-S_L^{-1}(x,T)]x
\\
&
-\overline{s}[S_R^{-1}(s,T)-S_L^{-1}(x,T)]\big](x^2-2s_0x+|s|^2)^{-1}.
\end{split}
\end{equation}
\end{theorem}
If we follow the idea of the Riesz-Dunford functional calculus and formally replace $x$ by $T$ in the slice hyperholomorphic Cauchy-formula, then we obtain its natural generalization to the quaternionic setting \cite{Colombo:2010b}.
\begin{definition}[$S$-functional calculus for bounded operators]\label{SCalcBd}
Let $T\in\boundOP(V)$, choose $\I\in\SS$ and set $ds_{\I }= -\I\, ds$. For $f\in\lhol(\sigma_S(T))$, we choose a bounded slice Cauchy domain $U$ with  $\sigma_{S}(T)\subset U$ and  $\clos{U}\subset \fdom(f)$ and define
\begin{equation}\label{EQLCalcBd}
f(T) := \frac{1}{2\pi}\int_{\partial(U\cap\cc_{\I})} S_L^{-1}(s,T)\, ds_{\I}\, f(s).
\end{equation}
For $f\in\rhol(\sigma_S(T))$, we choose again a bounded slice Cauchy domain $U$ with $\sigma_{S}(T)\subset U$ and $\clos{U}\subset \fdom(f)$ and define
\begin{equation}\label{EQRCalcBd}
f(T) := \frac{1}{2\pi}\int_{\partial(U\cap\cc_{\I})} f(s)\, ds_{\I}\, S_R^{-1}(s,T).
\end{equation}
These integrals are independent of the choices of the slice domain $U$ and the imaginary unit $\I\in\SS$.
\end{definition}
In order to extend these definitions from bounded to closed operators one can, as in the complex case, transform the closed operator into a bounded one and then apply the $S$-functional calculus for bounded operators, cf. \cite{Colombo:2011}. The obtained operator can then be represented by an integral of Cauchy type. Since this procedure works in the quaternionic setting only if the $S$-resolvent set of the operator contains a real point, we shall instead directly start from a Cauchy integral, cf. \cite{Gantner:2017}.

We say that a function is left (or right) slice hyperholomorphic at infinity, if $f$ is left (or right) slice hyperholomorphic on its domain $\fdom(f)$ and if there exists $r>0$ such that $\hh\setminus B_{r}(0)$ is contained $\fdom(f)$ of $f$ and  $f(\infty) = \lim_{x\to\infty}f(x)$ exists.
\begin{definition}[$S$-functional calculus for closed operators]\label{SCalcCl}
Let $T\in\closOP(V)$ with $\rho_S(T) \neq 0$, choose $\I\in\SS$ and set $ds_{\I} = -\I\, ds$. For $f\in\lhol(\sigma_{S}(T)\cup\{\infty\})$, we choose an unbounded slice Cauchy domain $U$ with $\sigma_{S}(T)\subset U$ and $\clos{U}\subset \fdom(f)$ define
\begin{equation}\label{EQLCalcCl}
f(T) := f(\infty)\id + \frac{1}{2\pi}\int_{\partial(U\cap\cc_{\I})} S_L^{-1}(s,T)\, ds_{\I}\, f(s).
\end{equation}
For $f\in\rhol(\sigma_S(T)\cup\{\infty\})$, we choose again a bounded slice Cauchy domain $U$ with $\sigma_{S}(T)\subset U$ and $\clos{U}\subset \fdom(f)$ and define
\begin{equation}\label{EQRCalcCl}
f(T) := f(\infty)\id + \frac{1}{2\pi}\int_{\partial(U\cap\cc_{\I})} f(s)\, ds_{\I}\, S_R^{-1}(s,T).
\end{equation}
These integrals are independent of the choices of the slice domain $U$ and the imaginary unit $\I\in\SS$.
\end{definition}
Observe that we can include the case of bounded operators into the one of unbounded operators, if we extend $f$ by setting $f(s) =0$ on a neighborhood of $\infty$ that is far from the $S$-spectrum of the operator. Since $\sigma_{SX}(T) = \sigma_{S}(T)$ if $T$ is bounded and $\sigma_{SX}(T) = \sigma_{S}(T)\cup\{\infty\}$ if $T$ is unbounded, we shall,  for neatness, denote the classes of admissible functions by $\lhol(\sigma_{SX}(T))$ and $\rhol(\sigma_{SX}(T))$ in the following. Most of the following properties of the $S$-functional calculus can be found in \cite{Colombo:2011}. The product rule  and the existence of Riesz-projectors were shown in \cite{Alpay:2015} for bounded and in \cite{Gantner:2017} for unbounded operators and the compatibility of the $S$-functional calculus for closed operators with intrinsic polynomials can be found in \cite{Gantner:2017}.
\begin{lemma}\label{SCProp}
Let $T\in\closOP(V)$ with $\rho_S(T)\neq \emptyset$. The $S$-functional calculus has the following properties:
\begin{enumerate}[label = (\roman*)]
\item If $f\in\lhol(\sigma_{SX}(T))$ or $f\in\rhol(\sigma_{SX}(T))$, then $f(T)$ is a bounded operator. 
\item  If $f\in\intrin(\sigma_{SX}(T))$, then the left slice hyperholomorphic and the right slice hyperholomorphic approach are consistent, that is \eqref{EQLCalcBd} and \eqref{EQRCalcBd}  resp. \eqref{EQLCalcCl} and \eqref{EQRCalcCl} yield the same operator.
\item If $f,g\in\lhol(\sigma_{SX}(T))$ and $a\in\hh$, then $(fa+g)(T) = f(T)a+g(T)$.  If $f,g\in\rhol(\sigma_{SX}(T))$ and $a\in\hh$, then $(af+g)(T) = af(T)+g(T)$.
\item If $f\in\intrin(\sigma_{SX}(T))$ and $g\in\lhol(\sigma_{SX}(T))$ or if $f\in\rhol(\sigma_{SX}(T))$ and $g\in\intrin(\sigma_{SX}(T))$, then $(fg)(T) = f(T)g(T)$.
\item If $g\in\intrin(\sigma_{SX}(T))$, then $ \sigma_S(g(T)) = g(\sigma_{SX}(T))$ and $f(g(T)) = (f\circ g)(T)$ if $f\in \lhol(g(\sigma_S(T)))$ or $f\in\rhol(g(\sigma_S(T)))$.
\item If $P$ is an intrinsic polynomial of order $n\geq 1$, then $P$ is not slice hyperholomorphic at infinity and does hence not belong to $\intrin(\sigma_{SX}(T))$ if $T$ is unbounded. However, even if $T$ is unbounded and $f\in\intrin(\sigma_{SX}(T))$ has a zero of order greater or equal to $n$ at infinity (i.e. $\lim_{p\to\infty}P(p)f(p) = 0$), then $P(T)f(T) = (Pf)(T)$, where $P(T)$ is defined as usual as $T^{0} = \id$ with $\dom(T^0) = V$ and $T^{n+1}\bv = T (T^n\bv)$ for $\bv\in\dom(T^n)$ such that $T^n\bv \in \dom(T)$.
\item If $\sigma$ is an open and closed subset of  $\sigma_{SX}(T)$, let $\chi_{\sigma}$ be equal to $1$ on an axially symmetric neighborhood of $\sigma$ in $\hh_{\infty}$ and equal to $0$ on an axially symmetric neighborhood of $\sigma_{SX}(T)\setminus \sigma$ in $\hh_{\infty}$. Then $\chi_{\sigma}\in\intrin(\sigma_{SX}(T))$ and $\chi_{\sigma}(T)$ is a projection onto a right linear subspace of $V$ that is invariant under $T$. Moreover, if we denote the restriction of $T$ to the range of $\chi_{\sigma}(T)$ by $T_{\sigma}$, then $\sigma_{SX}(T_{\sigma}) = \sigma$.
\end{enumerate}
\end{lemma}
Finally,  \cite{Colombo:2016a} showed the Taylor series expansion of the $S$-functional calculus in the operator.
\begin{theorem}\label{SCTaylor}
Let $T,N\in\boundOP(V)$ such that $TN = NT $ and such that $\sigma_{S}(T) \subset \clos{B_{\varepsilon}(0)}$ for some $\varepsilon >0$. Then
\[
\sigma_{S}(T+N) \subset \clos{B_{\varepsilon}(\sigma_{S}(T))} = \{ s\in\hh: \dist(s,\sigma_{S}(T)\leq \varepsilon\}
\]
with $\dist(s,\sigma_{S}(T)) = \inf_{x\in\sigma_{S}(T)}|s-x|$. For any $f\in \lhol\left( \clos{B_{\varepsilon}(\sigma_{S}(T))} \right)$ and any $g\in \rhol\left( \clos{B_{\varepsilon}(\sigma_{S}(T)} \right)$, we have
\[
f(T+N) = \sum_{n=0}^{+\infty}N^n \frac{1}{n!}\left(\sderiv^nf\right)(T)\qquad \text{and}\qquad g(T+N) = \sum_{n=0}^{+\infty} \frac{1}{n!}\left(\sderiv^ng\right)(T)N^n.
\]
\end{theorem}

\subsection{The spectral theorem for normal operators}\label{SpectThmSect}
We conclude this section by recalling the spectral theorem for bounded normal operators on a quaternionic Hilbert space $\hil$. This theorem was shown in \cite{Alpay:2016} and later also in \cite{Ghiloni:2017}. The two articles use different strategies for proving this theorem and they consider  two different approaches towards spectral integration in the quaternionic setting.

 We recall that $\hil$ is only a right Banach space. Therefore the space $\boundOP(\hil)$ of all bounded right linear operators on $\hil$ is only a real Banach space, cf. \Cref{B(V)Struct}. 
\begin{definition}\label{SpecMeasHil}
A spectral measure $E$ over $\cc_{\I}^{\geq}$ on the quaternionic Hilbert space $\hil$ is a set function defined on the Borel sets $\Borel\big(\cc_{\I}^{\geq}\big)$ of $\cc_{\I}^{\geq}$ the values of which are orthogonal projections on $\hil$ such that
\begin{enumerate}[(i)]
\item  $E\big(\cc_{\I}^{\geq}\big) = \id$
\item $E(\Delta_{1}\cap\Delta_{2}) = E(\Delta_{1}) E(\Delta_2)$ for all $\Delta_1,\Delta_2\in\Borel\big(\cc_{\I}^{\geq}\big)$
\item  $E(\cup_{n\in\nn} \Delta_{n})\bv = \sum_{n\in\nn} E(\Delta_n)\bv$ for all $\bv\in\hil$ and any sequence  $(\Delta_n)_{n\in\nn}$ of pairwise disjoint sets in~$\Borel\big(\cc_{\I}^{\geq}\big)$. 
\end{enumerate}
\end{definition}

In the complex setting, integrals with respect to a spectral measure $E$ defined on the Borel sets $\Borel(\cc)$ of $\cc$ are defined via approximation by simple functions. For a simple function $f(s) = \sum_{\ell=1}^{n} a_{\ell}\chi_{\Delta_\ell}(s)$, where $\chi_{\Delta_{\ell}}$ denotes the characteristic function of the set $\Delta_{\ell}\in\Borel(\cc)$, one defines 
\begin{equation}\label{SINt}
\int_{\cc}f(z)\,dE(z) = \sum_{\ell=1}^{n}a_{\ell}E(\Delta_{\ell}),
\end{equation}
 and for arbitrary bounded and measurable $f$ one defines  
 \begin{equation}\label{SINt2}
 \int_{\cc}f(z)\,dE(z) = \lim_{n\to+\infty}\int_{\cc} f_{n}(z)\,dE(z),
 \end{equation}
  where $f_n$ is a suitable sequence of simple functions converging uniformly to $f$. In order to proceed similarly in the quaternionic setting, we hence need additional structure (i.e. a left multiplication) on $\hil$. Otherwise only real-valued functions can be integrated  because \eqref{SINt} is only defined for real coefficients $a_{\ell}$. 
  
 \begin{definition}\label{LMult}
 A left multiplication on $\hil$ is a real algebra-homomorphism $\mathcal{L}: \hh\to\boundOP(\hh), a\mapsto L_a$ such that $L_a = a \id$ for any $a\in\rr$ and $(L_a)^* = L_{\overline{a}}$ for any $a\in\hh$. 
 \end{definition}
 If it is clear which left-multiplication we consider, we will also write $a\bv$ instead of  $L_a\bv$. Since $\mathcal{L}$ is a real algebra-homomorphism and $a\bv = L_a\bv = \bv a$ for any $a\in\rr$, such a left multiplication turns $\hil$ into a two-sided quaternionic vector space. Since moreover $L_{\overline{a}} = (L_a)^*$, we have $\langle \bu, a \bv\rangle = \langle \overline{a}\bu,\bv\rangle$ and in turn 
 \[
 \| a\bv\|^2 = \langle a\bv, a\bv\rangle = \langle \bv, \overline{a}a\bv\rangle = \langle \bv,\bv\rangle |a|^2 = |a|^2\|\bv\|^2.
 \]
Hence, any left multiplication as in \Cref{LMult} turns $\hil$ into a two-sided quaternionic Banach space. 

If $\mathbf{B} = (\bb_{\ell})_{ \ell\in \varLambda}$ is an orthonormal basis of $\hil$, then the left multiplication induced by $\mathbf{B}$ is the map $\mathcal{L}_{\mathbf{B}}$ given by $a\mapsto L_{a} : = \sum_{\ell\in\varLambda} \bb_{\ell} a \langle\bb_{\ell}, \bv\rangle$ for $a\in\hh$. If on the other hand $\hil$ is endowed with a left multiplication, then the real component space $\hil_{\rr}$ defined \eqref{VR} is, endowed with the scalar product on $\hil$, a real Hilbert space. Indeed, since $a \bv = \bv a $ for $a\in\hh$ and $\bv\in\hil_{\rr}$, we have for $\bu,\bv\in\hil_{\rr}$ that
\[
a \langle \bu, \bv \rangle = \langle \bu \overline{a}, \bv\rangle = \langle \overline{a}\bu,\bv\rangle = \langle\bu, a\bv\rangle= \langle \bu, \bv a \rangle = \langle \bu,\bv\rangle a \qquad \forall a\in\hh,
\]
which implies that $\langle \bu,\bv\rangle$ is real. 
Any orthonormal basis $(\bb)_{\ell\in\varLambda}$ of this space is also an orthonormal basis  of the quaternionic Hilbert space $\hil$ and  induces the left multiplication given on $\hil$ as $a\bv = a \sum_{\ell\in\varLambda}\bb_{\ell}\langle \bb_{\ell},\bv\rangle = \sum_{\ell\in\varLambda}\bb_{\ell}a \langle \bb_{\ell},\bv\rangle$ for any $a\in\hh$. (This fact was also shown in \cite[Theorem~4.3]{Ghiloni:2017} with different arguments).

 For any imaginary unit $\I\in\SS$, the mapping $\bv \mapsto \I\bv = L_{\I}\bv = \sum_{\ell\in\varLambda} \bb_{\ell}\I \langle \bb_{\ell},\bv\rangle$ is a unitary anti-selfadjoint operator on $\hil$. Conversely, if $\mathsf{J}$ is a unitary and anti-selfadjoint operator and $\I\in\SS$ is any imaginary unit, then there exists an orthonormal basis $(\bb_{\I,\ell})_{ \ell\in \varLambda}$ of $\hil$ such that $\mathsf{J}\bv =  \sum_{\ell\in\varLambda} \bb_{\I,\ell}\I \langle {\bb}_{\I,\ell},\bv\rangle$. Hence, applying $\mathsf{J}$ can be considered as a multiplication with $\I$ on the left. Note however, that $\I$ is not determined by $\mathsf{J}$. We can interpret applying $\mathsf{J}$ as a multiplication with any $\I\in\SS$, but obviously we cannot make this identification for two different imaginary units simultaneously. If we choose a different imaginary unit $\J\in\SS$, then we do moreover not have $\mathsf{J}\bv =  \sum_{\ell\in\varLambda} {\bb}_{\I,\ell}\J \langle {\bb}_{\I,\ell},\bv\rangle$. Instead, $\mathsf{J}\bv =  \sum_{\ell\in\varLambda} {\bb}_{\J,\ell}\J \langle {\bb}_{\J,\ell},\bv\rangle$, where $\bb_{\J,\ell} = \bb_{\I,\ell}h$ with $h\in\hh$ such that $|h|=1$ and $\J = h^{-1} \I h$, cf. \eqref{EVSym}.

Any left multiplication on $\hil$  is fully determined by its multiplications $L_{\I}, L_{\J}$ with two different imaginary units $\I,\J\in\SS$ with $\I\perp\J$. Indeed, setting $\K = \I\J$, any quaternion $a\in\hh$ can then be written as $a = a_0 + \I a_1 + \J a_2 + \K a_3$ and we have $a \bv = \bv a_0 + L_{\I} \bv a_1 + L_{\J}\bv a_2 + L_{\I}L_{\J}\bv a_3$. Observe that $L_{\I}$ and $L_{\J}$ are two  unitary anti-selfadjoint operators that anti-commute. Conversely, if $\mathsf{I}$ and $\mathsf{J}$ are two such operators, then we can choose $\I,\J\in\SS$ with $\I\perp\J$ and define a left multiplication on $\hil$ by setting $L_{\I}:=\mathsf{I}$ and $L_{\J}:=\mathsf{J}$. However,  $\I$ and $\J$ are  arbitrary. Hence, any couple consisting of two anti-commuting unitary and anti-selfadjoint operators can be used to generate an infinite number of distinct left multiplications on $\hil$.

Assume now that we are given a spectral measure $E$ on $\hil$ over $\cc_{\I}^{\geq}$ and a unitary anti-selfadjoint operator $\mathsf{J}$ that commutes with $E$. We can  interpret the application of $\mathsf{J}$ as a multiplication from the left with an arbitrary but fixed imaginary unit $\I\in\SS$. Since $E$ and $J$ commute, the imaginary unit $\I$ does  also commute with $E$ so that \eqref{SINt} is meaningful not only for real coefficients $a_{\ell}$, but even for coefficients in $\cc_{\I}$. We can hence define spectral integrals for $\cc_{\I}$-valued functions via the usual procedure in \eqref{SINt} and \eqref{SINt2}. Observe that if $f(z) = \alpha(z) + \I\beta(z)$ with $\alpha(z),\beta(z)\in\rr$ and $f_n(z) = \sum_{\ell=1}^{N_{n}}(\alpha_{n,\ell} + \I\beta_{n,\ell}) \chi_{\Delta_{n}}(z)$ with $\alpha_{n,\ell},\beta_{n,\ell}\in\rr$ for $n\in\nn$ is a sequence of measurable simple functions tending uniformly to $f$, then
\begin{align*}
\int_{\cc_{\I}^{\geq}}f(z)\,dE(z) = &\lim_{n\to+\infty}\sum_{\ell=1}^{N_{n}}(\alpha_{n,\ell} + \I\beta_{n,\ell}) E(\Delta_n)\\
=& \lim_{n\to+\infty}\sum_{\ell=1}^{N_{n}}\alpha_{n,\ell} E(\Delta_{n}) + \mathsf{J}  \lim_{n\to+\infty}\sum_{\ell=1}^{N_{n}}\beta_{n,\ell} E(\Delta_{n}) \\
=& \int_{\cc_{\I}^{\geq}}\alpha(z)\,dE(z) + \mathsf{J} \int_{\cc_{\I}^{\geq}}\beta(z)\,dE(z).
\end{align*}
The characterization of this spectral integral given in the following lemma was shown in \cite[Lemma~5.3]{Alpay:2016}.
\begin{lemma}\label{FabLem}
Let $E$ be a quaternionic spectral measure on $\hil$ over $\cc_{\I}^{\geq}$ and let $J$ be a unitary and anti-selfadjoint operator $J$ that commutes with $E$. If we interpret $J$ as the multiplication with $\I$ from the left and $f$ is a bounded measurable $\cc_{\I}$-valued function on $\cc_{\I}^{\geq}$ with $f(z) = \alpha(z) + \I\beta(z)$ where $\alpha(z),\beta(z)\in\rr$, then
\[
\left\langle \bu, \int_{\cc_{\I}^{\geq}}f(z)\,dE(z)\bv\right\rangle = \int_{\cc_{\I}^{\geq}}\alpha(z)\,d\langle \bu, E(z)\bv\rangle +  \int_{\cc_{\I}^{\geq}}\beta(z)\,d\langle \bu, E(z)\mathsf{J}\bv\rangle
\]
for all $\bu,\bv\in\hil$, where the quaternion-valued measure  $\langle \bu, E \bv\rangle$  on $\big(\cc_{\I}^{\geq},\Borel\big(\cc_{\I}^{\geq}\big)\big)$ is given by $\Delta\mapsto \langle \bu, E(\Delta) \bv\rangle$.
\end{lemma}

If $T$ is a bounded normal operator on $\hil$ (i.e. a bounded operator that commutes with its adjoint $T^*$), then we have
\begin{equation}\label{SplitTeasy}
T = A + \mathsf{J}_0B = \frac{1}{2} (T + T^*) + \mathsf{J}_0 \frac{1}{2}\left |T -T^*\right|,
\end{equation}
where $A := \frac{1}{2} (T + T^*)$ is self adjoint. The operator $C:= \frac{1}{2}(T - T^*)$ is anti-selfadjoint. Hence the polar decomposition theorem \cite[Theorem~2.20]{Ghiloni:2013} implies the existence of an anti-selfadjoint partially unitary operator $\mathsf{J}_0$ such that $C = \mathsf{J}_0|C|$, where $|C| = \sqrt{C^*C}$ is the unique positive square root of the positive operator $C^*C$. The operator $J_0$ is partially unitary with  $\ker \mathsf{J}_0 = \ker C$, i.e  the restriction of $J_0$ to  $\ran J_{0} = \clos{\ran {C}} = \ker\mathsf{J}_{0}^{\perp}$ is a unitary operator. Moreover, $A$, $B$, and $\mathsf{J}_0$ commute mutually and also with any operator that commutes with $T$ and $T^*$. Finally, we have $T^* = A - \mathsf{J}_0B$.  Teichm\"{u}ller already showed these facts in 1936 in \cite{Teichmuller:1936}, but proofs in English can also be found  in \cite{Ghiloni:2013}. As shown in \cite{Ghiloni:2013}, the operator $\mathsf{J}_0$ can furthermore be extended to a unitary and anti-selfadjoint operator $\mathsf{J}$ on all of $\hil$ that commutes with $T$ and $T^*$ and such that 
\begin{equation}\label{SplitT}
T = A + \mathsf{J} B.
\end{equation}
We however stress that unlike $\mathsf{J}_0$ the operator $\mathsf{J}$ is not unique.

Alpay, Colombo and Kimsey used these facts to show the spectral theorem for normal quaternionic linear operators in \cite{Alpay:2016}. We shall only recall the result for bounded operators \cite[Theorem~4.7]{Alpay:2016}, the proof of which is based on the continuous functional calculus for quaternionic normal operators introduced in \cite{Ghiloni:2013}. Let $p(z_0,z_1) = \sum_{0\leq |\ell|\leq n}a_{\ell}z_0^{\ell_{1}}z_{1}^{\ell_{2}}$ with the multi-index $\ell=(\ell_1,\ell_2)$ be a polynomial in the variables $z_0$ and $z_1$ with real coefficients $a_{\ell}$. The function $p(s) := p(s_0,s_1)$ for $s = s_0 + \I_s s_{1}\in\hh$ with $s_{1}\geq 0$ is then a real-valued and hence intrinsic slice function. For a normal operator $T$ decomposed as in \eqref{SplitT}, we can then define 
\begin{equation}\label{FAB}
p(T) := p(A,B) = \sum_{0\leq |\ell|\leq n}a_{\ell}A^{\ell_{1}}B^{\ell_{2}}.
\end{equation}
The operator $\mathsf{J}$ is essentially a multiplication with imaginary units $\I\in\SS$. Hence, we can define $f(T) = p_1(A,B) + Jp_2(A,B)$ for any intrinsic slice function function $f(s) = p_1(s_0,s_1) + \I p_2(s_0,s_1)$ with real-valued polynomials $p_1$ and $p_2$. We have $\|f(T)\| = \sup_{z\in\sigma_{S}(T)}|f(z)|$. If $f$ is any continuous slice function on $\sigma_{S}(T)$, then the Weierstrass approximation theorem implies the existence of a sequence of intrinsic slice functions $f_n = p_{n,1} + \I p_{n,2}$ with real-valued polynomials $p_{n,1},p_{n,2}$ such that $f_n$ tends to $f$ uniformly on $\sigma_{S}(T)$. Hence we can define $f(T) = \lim_{n\to\infty}f_n(T)$, where the sequence $f_n(T)$ converges in $\boundOP(\hil)$. We then have $\sigma_{S}(f(T)) = f(\sigma_{S}(T))$.

Alpay, Colombo and Kimsey follow in \cite{Alpay:2016} a well-known strategy from the complex case in order to show the spectral theorem for normal quaternionic linear operators. We choose $\I\in\SS$. By \Cref{ExtensionOfRealSlicefFs}, the mapping $f \mapsto f_{\I} = f|_{\cc_{\I}^{\geq}}$ determines a bijective relation between the set $\SCon(\sigma_{S}(T),\rr)$ of all real valued continuous slice functions on $\sigma_{S}(T)$ and the set $C(\Omega_{\I},\rr)$ of all real-valued continuous functions on $\Omega_{\I} := \sigma_{S}(T)\cap\cc_{\I}^{\geq}$. It also defines a bijective relation between the set $\SCon(\sigma_{S}(T))$ of all continuous intrinsic slice functions on $\sigma_{S}(T)$ and the set $C_{0}(\Omega_{\I},\cc_{\I})$ of all continuous $\cc_{\I}$-valued functions $f_{\I}$ on $\Omega_{\I}$ with $f_{\I}(\rr\cap\Omega_{\I})\subset\rr$. For any $\bv\in\hil$, the mapping $f_{\I}\mapsto \langle \bv, f(T) \bv\rangle$ is a continuous and positive $\rr$-linear functional on $C(\Omega_{\I},\rr)$. The Riesz representation theorem hence implies the existence of a positive Borel measure $\mu_{\bv,\bv}$ on $\Omega_{\I}$ such that 
\[
\langle \bv, f(T) \bv \rangle = \int_{\Omega_{\I}} f_{\I}(z)\, d\mu_{\bv,\bv}(z)\qquad \forall f_{\I}\in C(\Omega_{\I},\rr).
\]
Using the polarization identity, the authors deduced for any $\bu,\bv\in\hil$ the existence of a quaternion-valued Borel measure $\mu_{\bu,\bv}$ such that
\[
\langle \bu, f(T)\bv\rangle = \int_{\Omega_{\I}} f_{\I}(z)\,d\mu_{\bu,\bv}(z)\qquad \forall f_{\I}\in C(\Omega_{\I},\rr).
\]
For each $\Delta\in\Borel(\Omega_{\I})$ and each $\bu\in\hil$, the map $\bv \mapsto \mu_{\bu,\bv}(\Delta)$ is then a continuous quaternionic right linear functional on $\hil$. Hence there exists a unique $\bw \in\hil$ such that $\mu_{\bu,\bv}(\Delta) = \langle \bw , \bv\rangle $ for all $\bv\in\hil$. We define $E_{\I}(\Delta)\bu:= \bw$ such that
\[
\langle E_{\I}(\Delta)\bu, \bv\rangle = \langle \bw, \bv \rangle = \mu_{\bu,\bv}(\Delta).
\]
The operator $E_{\I}(\Delta)$ is then an orthogonal projection on $\hil$ and the mapping $\Delta\mapsto E_{\I}(\Delta)$ turns out to be a quaternionic spectral measure over $\cc_{\I}^{\geq}$ on $\hil$. Finally, one arrives at the spectral theorem for bounded normal operators on quaternionic Hilbert spaces \cite[Theorem~4.7]{Alpay:2016}.

\begin{theorem}\label{SpecThm1}
Let $T\in\boundOP(\hil)$ be normal, let $T = A+\mathsf{J}B$  as in \eqref{SplitT}, consider $\mathsf{J}$ as a multiplication with $\I\in\SS$ and set $\Omega_{\I}:=\sigma_{S}(T)\cap\cc_{\I}^\geq$ . Then there exists a unique quaternionic spectral measure $E_{\I}$ over $\cc_{\I}^{\geq}$ such that
\begin{equation}\label{JJJ}
f(T) = \int_{\Omega_{\I}} f_{\I}(z)\,dE_{\I}(z)\qquad\forall f\in\SCon(\sigma_{S}(T)),
\end{equation}
where the spectral integral is intended as above. cf. also \Cref{FabLem}.
Moreover, if $\J\in\SS$ and we define $\varphi_{\I,\J}(z_0 + \I z_1) = z_0 + \J z_1$, then $E_{\I}(\Delta) = E_{\J}(\varphi_{\I,\J}(\Delta))$ for any $\Delta\in\Borel(\Omega_{\I})$. 
\end{theorem}
\begin{remark}\label{RemSpecThm}
We want to point out that the construction of the spectral measure $E$ only used real-valued functions on $\Omega_{\I}$. These functions are restrictions of real-valued slice functions.

Moreover, the spectral measure was actually constructed using functions of $A$ and $B$ in \eqref{FAB}. Hence, it depends only on $A$ and $B$ but not on $J$. The spectral measures of $T = A + \mathsf{J}B$ and $T^*= A- \mathsf{J}B$ therefore coincide. As we see, in the quaternionic setting, a normal operator is not fully determined by its spectral measure. Essential information is also contained in the operator $\mathsf{J}$, which will force us to introduce the notion of a spectral system in \Cref{SpecSys}.
\end{remark}

In \cite{Ghiloni:2017} the authors go one step further in their theory of spectral integration and introduce {\em intertwining quaternionic projection-valued measures}, which allow them to define spectral integrals for $\hh$-valued and not only for $\cc_{\I}$-valued functions.
\begin{definition}
Let $\I\in\SS$. An intertwining quaternionic projection-valued measure (for short iqPVM) over $\cc_{\I}^{\geq}$ on $\hil$ is a couple $\mathcal{E} = (E, \mathcal{L})$ consisting of a quaternionic spectral measure $E$ over $\cc_{\I}$ and a left multiplication $\mathcal{L}: \hh\to \boundOP(\hil), a\mapsto L_{a}$ that commutes with $E$, i.e. such that $E(\Delta)L_{a} = L_{a} E(\Delta)$ for any $\Delta\in\Borel\big(\cc_{\I}^{\geq}\big)$ and any $a\in\hh$.
\end{definition}
With respect to an iqPVM $\mathcal{E}$ one can define spectral integrals of functions with arbitrary values in the quaternions since \eqref{SINt} is meaningful for simple functions with arbitrary quaternionic coefficients $a_{\ell}$. Ghiloni, Moretti and Perotti  show the following result \cite[Theorem~4.1]{Ghiloni:2017}
\begin{theorem}\label{SpecThm2}
Let $T\in\boundOP(\hil)$ be normal and let $\I\in\SS$. There exists an iqPVM $\mathcal{E} = (E, \mathcal{L})$ over $\cc_{\I}^{\geq}$ on $\hil$ such that
\begin{equation}\label{IJI}
T = \int_{\cc_{\I}^{\geq}} z\,d\mathcal{E}(z).
\end{equation}
The spectral measure $E$ is uniquely determined by $T$ and the left multiplication is uniquely determined for $a\in\cc_{\I}$ on $\ker(T-T^*)^{\perp}$. Precisely, we have for any other left multiplication $\mathcal{L}'$  such that $\mathcal{E}'=(E,\mathcal{L}')$ is an iqPVM satisfying \eqref{IJI} that $L_{a}\bv = L_{a}'\bv$ for any $a\in\cc_{\I}$ and any $\bv\in\ker(T-T^*)^{\perp}$. (Even more specifically, we have $\I\bv  = \mathsf{J}_{0}\bv$ for any $\bv\in \ker(T-T^*)^{\perp} = \ran \mathsf{J}_{0}$.)
\end{theorem}

The proof of the above result in \cite{Ghiloni:2017} follows a completely different strategy than the proof of the spectral theorem for bounded operators in \cite{Alpay:2016}. Similar to \cite{Sharma:1987, Viswanath:1971}, it reduces the quaternionic problem to a complex problem in order to apply the well-known results from complex operator theory. Instead of working with the symplectic image, the authors however apply the classical results to the restriction of $T$ to a suitably chosen complex component space.

We decompose $T = A + \mathsf{J}B$ as in \eqref{SplitT} and choose $\I,\J\in\SS$ with $\I\perp\J$. Since $\mathsf{J}$ is unitary and anti-selfadjoint, there exists an orthonormal basis $(\bb_{\ell})_{\ell\in\varLambda}$ of $\hil$ such that $\mathsf{J}\bb_{\ell} = \bb_{\ell}\I$.  We can write any $\bv = \sum_{\ell\in\varLambda}\bb_{\ell}a_{\ell}\in\hil$ with $a_{\ell}\in\hh$ as $\bv = \bv_{1} + \bv_{2}\J = \sum_{\ell\in\varLambda}\bb_{\ell}a_{\ell,1} + \sum_{\ell\in\varLambda}\bb_{\ell}a_{\ell,2}\J$ with $a_{\ell,1},a_{\ell,2}\in\cc_{\I}$ such that $a_{\ell} = a_{\ell,1}+a_{\ell,2}\J$. Setting $\hil_{\mathsf{J},\I}^{+}:=\clos{\linspan{\cc_{\I}}\{\bb_{\ell}:\ell\in\varLambda \}}$ and $\hil_{\mathsf{J},\I}^{-}:=\hil_{\mathsf{J},\I}^{+}\J$, we find that
 $\hil = \hil_{\mathsf{J,}\I}^{+}\oplus \hil_{\mathsf{J},\I}^{-}$.  The spaces $\hil_{\mathsf{J},\I}^{+}$ and $\hil_{\mathsf{J},\I}^{-}$ are  $\cc_{\I}$-complex Hilbert spaces if we endow them with the right scalar multiplication on $\hil$ restricted to $\cc_{\I}$ and with the scalar product on $\hil$. Obviously $\hil_{\mathsf{J},\I}^{+}$ consists of all eigenvectors of $\mathsf{J}$ associated with the eigenvalue $\I$ and since $\I\J = \J(-\I)$ the space $\hil_{\mathsf{J},\I}^{-}$ consists of all eigenvectors of $\mathsf{J}$ associated with the eigenvalue $-\I$. For more detailed arguments we refer to \cite{Ghiloni:2013}, but we stress that these facts motivate \Cref{JVSplitThm} of the present paper.

For $\bv\in\hil_{\mathsf{J},\I}^{+}$, we have $\mathsf{J}(T\bv) = T( \mathsf{J}\bv) = (T\bv)\I$ and so $T\bv\in\hil_{\mathsf{J},\I}^{+}$. Hence, $T$ leaves $\hil_{\mathsf{J},\I}^{+}$ invariant. The restriction $T_+$ of $T$ to $\hil_{\mathsf{J},\I}^{+}$  defines a bounded normal operator on $\hil_{\mathsf{J},\I}^{+}$ with $\sigma(T_{+}) = \sigma_{S}(T)\cap\cc_{\I}^{\geq}$. Applying the spectral theorem for normal complex linear operators, one obtains a spectral measure $E_{+}$ on $\hil_{\J,\I}^{+}$, the support of which is $\sigma(T_+)$, such that $T_{+} = \int_{\sigma(T_+)} z\, dE_{+}(z)$. The quaternionic linear extension $E$ of $E_{+}$ to all of $\hil$ is then a quaternionic spectral measure over $\cc_{\I}^{\geq}$ on $\hil$. This extension is obtained by writing $\bv\in\hil$ as $\bv = \bv_1 + \bv_2\J$ with $\bv_1,\bv_2\in\hil_{\mathsf{J},\I}^{+}$ and setting $E(\Delta)\bv = E_{+}(\Delta)\bv_{1} + (E(\Delta)\bv_{2})\J$ for $\Delta\in\Borel\big(\cc_{\I}^{\geq}\big)$. If we furthermore choose a suitable orthonormal basis $\mathbf{B} := (\bb_{\ell})_{\ell\in\varLambda}$ of the $\cc_{\I}$-complex Hilbert space $\hil_{\mathsf{J},\I}^{+}$, then we find that $(\bb_{\ell})_{\ell\in\varLambda}$ is also an orthonormal basis of the quaternionic Hilbert space $\hil$ and that $\mathcal{E} = (E,\mathcal{L}_{\mathbf{B}})$, where $\mathcal{L}_{\mathbf{B}}$ is the left multiplication on $\hil$ induced by $\mathbf{B}$, is an iqPVM such that \eqref{IJI} holds true. Again this is a rough summary of the strategy and we refer to \cite{Ghiloni:2017} for the technical details.

Since $\mathbf{B}$ is an orthonormal basis of $\hil_{\mathsf{J},\I}^{+}$, it consists of eigenvectors of $\mathsf{J}$ with respect to $\I$. We hence find for $\bv\in\hil_{\mathsf{J},\I}^{+}$ that $\I\bv = L_{\I}\bv = \sum_{\ell\in\varLambda}\bb_{\ell}\I \langle\bb_{\ell},\bv\rangle  = \sum_{\ell\in\varLambda}\mathsf{J}\bb_{\ell} \langle\bb_{\ell},\bv\rangle = \mathsf{J}\sum_{\ell\in\varLambda}\bb_{\ell} \langle\bb_{\ell},\bv\rangle = \mathsf{J}\bv$ and in turn also for arbitrary $\bv = \bv_{1} + \bv_{2}\J\in\hil$ with $\bv_{1},\bv_{2}\in\hil_{\mathsf{J},\I}^{+}$ that  $\I\bv = \I\bv_{1} + (\I\bv_{2})\J = \mathsf{J}\bv_{1} +( \mathsf{J}\bv_{2})\J = \mathsf{J}\bv$. Hence also in the approach using iqPVM the application of $\mathsf{J}$ is interpreted as  multiplication with a randomly chosen imaginary unit $\I$ on the left. Moreover, this multiplication is only on $\ker(T^*-T)^{\perp}$ determined by $T$. Indeed, on this space $\mathsf{J}$ is coincides with $\mathsf{J}_{0}$, whereas the extension of $\mathsf{J}_0$ to a unitary anit-sefladjoint operator on all of $\hil$ is arbitrary.  The left-multiplication $L_{\J}$ for $\J$ with $\J\perp \I$, which together with $L_{\I}=\mathsf{J}$ fully determines $\mathcal{L}$ in $\mathcal{E}$, is however completely arbitrary and not at all determined by $T$. We stress these facts because they shall be important in the discussion in \Cref{DiffSpInt}.

\section{\texorpdfstring{Intrinsic $S$-functional calculus on one-sided Banach spaces}{Intrinsic S-functional calculus on one-sided Banach spaces}}\label{IntCalcSect}
Although the $S$-resolvents are only defined on two-sided quaternionic Banach spaces, spectral properties of a right linear operator should be independent of the left multiplication on the space as we pointed out in the introduction. In this section we therefore derive---at least for an intrinsic function $f$---an integral representation of the operator $f(T)$ obtained via the $S$-functional calculus that is independent of the left multiplication as it only uses the multiplication of vectors with scalars from the right. We then use this integral representation to define the $S$-functional calculus for intrinsic functions on a quaternionic right Banach space.

We start by investigating in detail the relations between some of the fundamental notions of operator theory defined by the quaternionic structure and those defined by the various complex structures on a quaternionic right Banach space $V_R$. As pointed out in \Cref{RealBS}, any quaternionic right Banach space $V_R$ can in a natural way be considered as a complex Banach space over any of the complex planes $\cc_{\I}$ by restricting the multiplication with quaternionic scalars from the right to $\cc_{\I}$. In order to deal with the different structures on $V_{R}$, we introduce the following notation.
\begin{notation}
Let  $V_R$ be a quaternionic right Banach space. For $\I\in\SS$, we denote the space $V_{R}$ considered as a complex Banach space over the complex field $\cc_{\I}$ by $V_{R,\I}$. If $T$ is a quaternionic right linear operator on $V_R$, then $\rho_{\cc_{\I}}(T)$ and $\sigma_{\cc_{\I}}(T)$ shall denote its resolvent set and spectrum as a $\cc_{\I}$-complex linear operator on $V_{R,\I}$. If $A$ is a $\cc_{\I}$-complex linear, but not quaternionic linear operator on $V_{R,\I}$, then we denote its spectrum as usual by $\sigma(A)$.

If we want to distinguish between the identity operator on $V_{R}$ and the identity operator on $V_{R,\I}$ we denote them by $\id_{V_{R}}$ and $\id_{V_{R,\I}}$. We point out that the operator $\lambda \id_{V_{R,\I}}$ for $\lambda\in\cc_{\I}$ acts as $\lambda\id_{V_{R,\I}}\bv = \bv \lambda$ because the multiplication with scalars on $V_{R,\I}$ is defined as the quaternionic right scalar multiplication on $V_{R}$ restricted to $\cc_{\I}$. 
\end{notation}
We shall now specify an observation made in \cite[Theorem~5.9]{Ghiloni:}, where the authors showed that $\rho_S(T)$ is the axially symmetric hull of $\rho_{\cc_{\I}}(T)\cap\overline{\rho_{\cc_{\I}}(T)}$. 
\begin{theorem}\label{ABC}
Let $T\in\closOP(V_R)$ and choose $\I\in\SS$. The spectrum $\sigma_{\cc_{\I}}(T)$ of $T$ considered as a closed complex linear operator on $V_{R,\I}$ equals $\sigma_S(T)\cap\cc_{\I}$, i.e.
\begin{equation}\label{CiRes}
 \sigma_{\cc_{_{\I}}}(T) = \sigma_S(T)\cap\cc_{\I}.
\end{equation}
For any $\lambda$ in the resolvent set $\rho_{\cc_{\I}}(T)$ of $T$ as a complex linear operator on $V_{R,\I}$, the $\cc_{\I}$-linear resolvent of $T$ is given by $R_{\lambda}(T) = \left(\overline{\lambda}\id_{V_{R,\I}} - T\right)\Q_{\lambda}(T)^{-1}$, i.e.
\begin{equation}\label{JuTZu}
 R_{\lambda}(T)\bv := \Q_{\lambda}(T)^{-1}\bv\overline{\lambda}-T\Q_{\lambda}(T)^{-1}\bv.
\end{equation}
For any $\J\in\SS$ with $\I\perp\J$, we moreover have
\begin{equation}\label{RLambdaCon}
R_{\overline{\lambda}}(T)\bv = -[R_{\lambda}(T)(\bv\J)]\J.
\end{equation}
\end{theorem}
\begin{proof}
Let $\lambda \in\rho_S(T)\cap\cc_{\I}$. The resolvent $(\lambda\id_{V_{R,\I}} - T)^{-1}$ of $T$ as a $\cc_{\I}$-linear operator on $V_{R,\I}$ is then given by \eqref{JuTZu}. Indeed, since $T$ and $Q_{\lambda}(T)^{-1}$ commute, we have for $\bv\in \dom(T)$ that
\begin{align*}
 =& (\overline{\lambda}\id_{V_{R,\I}} - T)\Q_{\lambda}(T)^{-1}(\bv\lambda - T\bv)\\
=& (\overline{\lambda}\id_{V_{R,\I}} - T)\left(\Q_{\lambda}(T)^{-1}\bv\lambda - T\Q_{\lambda}(T)^{-1}\bv \right)\\
=& \Q_{\lambda}(T)^{-1}\bv\lambda\overline{\lambda} - T\Q_{\lambda}(T)^{-1}\bv\overline{\lambda}  -T\Q_{\lambda}(T)^{-1}\bv\lambda + T^2\Q_{\lambda}(T)^{-1}\bv \\
=& (|\lambda|^2\id_{V_{R,\I}} - 2\lambda_0 T + T^2)\Q_{\lambda}(T)^{-1}\bv  = \bv.
\end{align*}
Similarly, for $\bv\in V_{R, \I} = V_{R}$, we have
\begin{align*}
&(\lambda\id_{V_{R,\I}} - T)R_{\lambda}(T)\bv \\
=& (\lambda\id_{V_{R,\I}} - T)\left(\Q_{\lambda}(T)^{-1}\bv\overline{\lambda} - T\Q_{\lambda}(T)^{-1}\bv\right)\\
=& \Q_{\lambda}(T)^{-1}\bv\overline{\lambda}\lambda - T\Q_{\lambda}(T)^{-1}\bv\lambda - T\Q_{\lambda}(T)^{-1}\bv\overline{\lambda} + T^2\Q_{\lambda}(T)^{-1}\bv\\
=& (|\lambda|^2\id_{V_{R,\I}} - 2\lambda_0 T + T^2)\Q_{\lambda}(T)^{-1}\bv  = \bv.
\end{align*}
Since $\Q_{\lambda}(T)^{-1}$ maps $V_{R,\I}$ to $\dom(T^2)\subset\dom(T)$, the operator $R_{\lambda}(T)= (\lambda\id_{V_{R,\I}} - T)\Q_{\lambda}(T)^{-1}$ is bounded and so  $\lambda$ belongs to the resolvent set $\rho_{\cc_{\I}}(T)$ of $T$ considered as a $\cc_{\I}$-linear operator on $V_{R,\I}$. Hence, $\rho_S(T)\cap\cc_{\I} \subset \rho_{\cc_{\I}}(T)$ and in turn $\sigma_{\cc_{\I}}(T) \subset \sigma_S(T)\cap\cc_{\I}$. Together with the axial symmetry of the $S$-spectrum, this further implies 
\begin{equation}\label{AAaa}
\sigma_{\cc_{\I}}(T) \cup \overline{\sigma_{\cc_{\I}}(T)} \subset ( \sigma_S(T)\cap\cc_{\I} ) \cup (\overline{\sigma_S(T)\cap\cc_{\I}}) = \sigma_{S}(T)\cap\cc_{\I},
\end{equation}
where $\overline{A} = \{\overline{z}: z\in A\}$. 

If $\lambda$ and $\overline{\lambda}$ both belong to $\rho_{\cc_{\I}}(T)$, then $[\lambda]\subset \rho_S(T)$ because
\begin{align*}
&(\lambda \id_{V_{R,\I}} - T) ( \overline{\lambda} \id_{V_{R,\I}}-T)\bv \\
=& (\bv\overline{\lambda})\lambda - (T\bv)\lambda - T(\bv\overline{\lambda}) +T^2\bv \\
= & (T^2 - 2\lambda_0 T + |\lambda |^2) \bv
\end{align*}
and hence $\Q_{\lambda}(T)^{-1} = R_{\lambda}(T)R_{\overline{\lambda}}(T)\in\boundOP(V_R)$. Thus $\rho_S(T) \cap\cc_{\I} \supset \rho_{\cc_{\I}}(T) \cap \overline{\rho_{\cc_{\I}}(T)}$ and in turn 
\begin{equation}\label{BBbb}
\sigma_{S}(T)\cap\cc_{\I} \subset \sigma_{\cc_{\I}}(T)\cup\overline{\sigma_{\cc_{\I}}(T)}.
\end{equation}
The two relations \eqref{AAaa} and \eqref{BBbb} yield together  
\begin{equation}\label{JaTsUJ}
\sigma_S(T) \cap\cc_{\I} = \sigma_{\cc_{\I}}(T)\cup\overline{\sigma_{\cc_{\I}}(T)}.
\end{equation}
What remains to show is that $\rho_{\cc_{\I}}(T)$ and $\sigma_{\cc_{\I}}(T)$ are symmetric with respect to the real axis, which then implies 
\begin{equation}\label{Jukza} 
\sigma_S(T) \cap\cc_{\I} = \sigma_{\cc_{\I}}(T)\cup\overline{\sigma_{\cc_{\I}}(T)} = \sigma_{\cc_{\I}}(T).
\end{equation}

 Let $\lambda\in\rho_{\cc_{\I}}(T)$ and choose $\J\in\SS$ with $\I\perp \J$. We show that  $R_{\overline{\lambda}}(T)$ equals the mapping $A \bv := -\left[R_{\lambda}(T)(\bv \J)\right]\J$. As $\lambda \J = \J\overline{\lambda}$ and $\J\lambda = \overline{\lambda}\J$, we have for $\bv\in\dom(T)$ that
\begin{align*}
&A \left(\overline{\lambda}\id_{V_{R,\I}}-T\right) \bv = A\left(\bv\overline{\lambda} - T\bv\right) \\
=& - \left[R_{\lambda}(T)\left((\bv \overline{\lambda})\J - (T\bv)\J\right)\right]\J \\
=& - \left[R_{\lambda}(T)((\bv\J) \lambda - T(\bv\J))\right]\J \\
=& -\left[ R_{\lambda}(T) (\lambda \id_{V_{R,\I}} - T) (\bv\J)\right] \J = -\bv\J\J = \bv.
\end{align*}
Similarly, for arbitrary $\bv\in V_{R,\I}= V_{R}$, we have
\begin{align*}
&\left(\overline{\lambda}\id_{V_{R,\I}}-T\right) A \bv = \left(A\bv\right)\overline{\lambda} - T\left(A\bv\right)\\
= &-\left[R_{\lambda}(T)(\bv\J)\right]\J\overline{\lambda} + T\left(\left[R_{\lambda}(T)(\bv\J)\right]\J\right) \\
= &-\left[ R_{\lambda}(T)(\bv\J)\lambda - T(R_{\lambda}(T)(\bv\J))\right] \J \\
= & - \left[(\lambda\id_{V_{R,\I}} - T )R_{\lambda}(T) (\bv\J)\right]\J  = -\bv\J\J = \bv.
\end{align*}
 Hence, if $\lambda\in\rho_{\cc_{\I}}(T)$, then $R_{\overline{\lambda}}(T) = -\left[R_{\lambda}(T)(\bv \J)\right]\J$ such that in particular $\overline{\lambda}\in\rho_{\cc_{\I}}(T)$. Consequently $\rho_{\cc_{\I}}(T)$ and in turn also $\sigma_{\cc_{\I}}(T)$ are symmetric with respect to the real axis such that \eqref{Jukza} holds true.
 
 \end{proof}
 \begin{remark}
 The relation \eqref{JaTsUJ} has already been observed in \cite[Theorem~5.9]{Ghiloni:}. For the sake of completeness, we repeated the respective arguments here. Also the relation $ R_{\lambda}(T) R_{\overline{\lambda}}(T) = \Q_{\lambda}(T)^{-1}$, which is a consequence of \eqref{JuTZu}, was understood by the authors. The novelty in the above theorem is hence  the fact that for a quaternionic linear operator $T$ automatically $\sigma_{\cc_{\I}}(T) = \overline{\sigma_{\cc_{\I}}(T)}$ due to \eqref{RLambdaCon}. For unitary operators, this symmetry was  already understood by Sharma and Coulson, as the author of the present paper learned only recently. In \cite{Sharma:1987} they showed the following: if $T$ is a unitary quaternionic linear operator on a quaternionic Hilbert space $\hil$, then the spectrum $\sigma(T_s)$ of the operator $T_s$ induced by $T$ on the symplectic image $\hil_{s}$ satisfies $\sigma(T_s) = \overline{\sigma(T_s)}$. Their strategy for showing this was however different. Since the $S$-spectrum and the associated pseudo-resolvent had not yet been developed, they showed this symmetry for the set of approximate eigenvalues of $T_{s}$  using \eqref{RightEVRel}. Since with $T$ also $T_s$ is unitary, $\sigma(T_{s})$ only consists of approximate eigenvalues and the statement follows.  
 \end{remark}
 
The above result allows us to write the operator $f(T)$ obtained via the $S$-functional calculus for an intrinsic slice hyperholomorphic function without making use of a quaternionic multiplication on $\boundOP(V_R)$. 
\begin{definition}\label{CRes}
Let $T\in\closOP(V_{R})$. We define the $V_{R}$-valued function
\[\RRes_{s}(T;\bv) = \Q_{s}(T)^{-1}\bv\overline{s}-T\Q_{s}(T)^{-1}\bv\qquad \forall \bv\in V_{R}, \ s\in\rho_S(T) .\]
\end{definition}
\begin{remark}
By \Cref{ABC}, the mapping $\bv \mapsto \RRes_{s}(T;\bv)$ coincides with resolvent of $T$ at $s$ applied to $\bv$ if $T$ is considered a $\cc_{\I_{s}}$-linear operator on $V_{R,\I_{s}}$.
\end{remark}

\begin{corollary}
Let $T\in\closOP(V_{R})$. For any $\bv\in V_R$, the function $\bof(s):= \RRes_{s}(T;\bv)$ is right slice hyperholomorphic on $\rho_S(T)$.
\end{corollary}
\begin{proof}
Obviously $\bof(s) = \boldsymbol{\alpha}(s_0,s_1) + \boldsymbol{\beta}(s_0,s_1)\I_s$ with $\boldsymbol{\alpha}(s_0,s_1)  = \Q_{s}(T)^{-1}\bv s_0-T\Q_{s}(T)^{-1}\bv$ and $\boldsymbol{\beta}(s_0,s_1) =-\Q_{s}(T)^{-1}\bv s_1$, which satisfy the compatibility condition \eqref{CCond}. The property that $\boldsymbol{\alpha}$ and $\boldsymbol{\beta}$ satisfy the Cauchy-Riemann-equations \eqref{CR} is equivalent to $\bof|_{\rho_S(T)\cap \cc_{\I}}$ being (right) holomorphic for any $\I\in\SS$. But this is true because $\bof |_{\rho_S(T)\cap\cc_{\I}}$ coincides with the resolvent of $T$ as an operator on $V_{R,\I}$ applied to $\bv$ by \Cref{ABC}, which is a holomorphic function on $\rho_{\cc_{\I}}(T) = \rho_{S}(T)\cap\cc_{\I}$.

\end{proof}

\begin{theorem}\label{IntCalcThm}
Let $T\in\closOP(V)$ be a closed operator on a two-sided quaternionic Banach space $V$ with $\rho_S(T)\neq \emptyset$ and let $f\in\intrin(\sigma_{SX}(T)))$.
For any $\I\in\SS$ and any unbounded slice Cauchy domain $U$ with $\sigma_{S}(T)  \subset U$ and $\clos{U}\subset \fdom(F)$, the operator $f(T)$ obtained via the $S$-functional calculus satisfies
\begin{equation}
\label{IntSCalc} 
f(T)\bv = \bv f(\infty) + \int_{\partial(U\cap\cc_{\I})}\RRes_{z}(T;\bv)f(z)\,dz \frac{-\I}{2\pi } \quad\forall \bv\in V.
\end{equation}
\end{theorem}
\begin{proof}
Let $U$ be a slice Cauchy domain such that $\sigma_S(T)\subset U$ and $\clos{U}\subset \dom(f)$. We then have for any $\I\in\SS$ and any $\bv\in V$ that
\begin{equation}\label{XX}
 f(T)\bv = f(\infty)\bv + \frac{1}{2\pi}\int_{\partial(U\cap\cc_{\I})} f(s)\,ds_{\I}\,S_R^{-1}(s,T)\bv.
\end{equation}
The boundary $\partial(U\cap\cc_{\I})$ of $U$ in $\cc_{\I}$ consists of a finite number of Jordan curves and is symmetric with respect to the real axis. We divide the set of these curves into three subsets: those taking only values  in $\cc_{\I}^{+}$, those taking only values in $\cc_{\I}^{-}$ and those taking values in both $\cc_{\I}^{+}$ and $\cc_{\I}^{-}$.

Let $\gamma_{\ell}: [0,1]\to\cc_{\I}^+$, $\ell = 1,\ldots N_0$ be those Jordan curves that belong to $\partial(U\cap\cc_{\I})$ and only take values in $\cc_{\I}^+$. Due to the symmetry of $\partial(U\cap\cc_{\I})$,  the paths $-\overline{\gamma_{\ell}}$ given by $(-\overline{\gamma_{\ell}})(t) := \overline{\gamma_{\ell}(1-t)}$ for $\ell = 1,\ldots N_0$ are exactly those Jordan curves that belong to $\partial(U\cap\cc_{\I})$ and only take values in $\cc_{\I}^{-}$. Finally, let us consider those Jordan curves in $\partial(U\cap\cc_{\I})$ that take values  in both $\cc_{\I}^+$ and $\cc_{\I}^-$. Each of these curves $\Gamma_{\ell},~\ell = N_0+1,\ldots,N$ can be decomposed as $\Gamma_{\ell} = \gamma_{\ell} \cup -\overline{\gamma_{\ell}}$, where $\gamma_{\ell}:[0,1]\to \cc_{\I}^+$ is the part of $\Gamma_{\ell}$ that lies in the upper halfplane $\cc_{\I}^+$ and $-\overline{\gamma_{\ell}}$, which is again defined as $(-\overline{\gamma_{\ell}})(t) = \overline{\gamma_{\ell}(1-t)}$, is the part of $\Gamma_{\ell}$ that lies in the lower halfplane $\cc_{\I}^-$. Hence, we have 
\[\partial(U\cap\cc_{\I}) = \bigcup_{\ell=1}^{N} \gamma_{\ell}\cup-\overline{\gamma_{\ell}}.\]
For technically more detailed arguments, we refer to \cite[Section~3]{Gantner:2017}, where also the integral representation \eqref{KuKA} has been shown. We include the following computations nevertheless, in order to keep the paper self-contained. 

With the above notation and because of $(-\overline{\gamma_{\ell}})'(t) = -\overline{\gamma_{\ell}'(1-t)}$, we have
\begin{align*}
&\int_{\partial(U\cap\cc_{\I})}f(s)\,ds_{\I}\,S_R^{-1}(s,T) \bv\\
=&\sum_{\ell=1}^N\int_{\gamma_{\ell}} f(s)\,ds_{\I}\,S_R^{-1}(s,T)\bv + \sum_{\ell=1}^N\int_{-\overline{\gamma_{\ell}}}f(s)\,ds_{\I}\,S_R^{-1}(s,T)\bv\\
=&\sum_{\ell = 1}^N\int_{0}^{1}f(\gamma_\ell(t)) (-\I) \gamma_\ell'(t)\left(\overline{\gamma_{\ell}(t)}-T\right)\Q_{\gamma_{\ell}(t)}(T)^{-1}\bv\,dt\\
&+\sum_{\ell = 1}^N\int_{0}^{1}f\big(\overline{\gamma_\ell(1-t)}\big) \I \overline{\gamma_\ell'(1-t)}(\gamma_{\ell}(1-t)-T)\Q_{\overline{\gamma_{\ell}(1-t)}}(T)^{-1}\bv\,dt.
\end{align*}
Observe that $\Q_{\overline{s}}(T)^{-1} = \Q_{s}(T)^{-1}$ for $s\in\rho_S(T)$. Since $f$ is intrinsic, we have $f(\overline{x}) = \overline{f(x)}$. Moreover $f(\gamma_{\ell}(1-t))\in \cc_{\I}$ as $\gamma_{\ell}(1-t)\in\cc_{\I}$ such that $f\big(\overline{\gamma_\ell(1-t)}\big)$,  $\I$ and  $\overline{\gamma_\ell'(1-t)}$ commute mutually and in turn 
\[
f\big(\overline{\gamma_\ell(1-t)}\big) \I \overline{\gamma_\ell'(1-t)} = \overline{f\left(\gamma_\ell(1-t)\right) (-\I) \gamma_\ell'(1-t)}.
\]
 After a change of variables in the integrals of the second sum, we therefore get
\begin{align}
&\notag \int_{\partial(U\cap\cc_{\I})}f(s)\,ds_{\I}\,S_R^{-1}(s,T)\bv \\
=&\notag \sum_{\ell = 1}^N\int_{0}^{1}f(\gamma_\ell(t)) (-\I) \gamma_\ell'(t)\left(\overline{\gamma_{\ell}(t)}-T\right)\Q_{\gamma_{\ell}(t)}(T)^{-1}\bv\,dt\\
&+\notag \sum_{\ell = 1}^N\int_{0}^{1}\overline{f(\gamma_\ell(t))(-\I)\gamma_\ell'(t)}(\gamma_{\ell}(t)-T)\Q_{\gamma(t)}(T)^{-1}\bv\,dt\displaybreak[1]\\
\begin{split}\label{KuKA}
=& \sum_{\ell=1}^N\int_{0}^{1}2\Re\left(f(\gamma_{\ell}(t))(-\I)\gamma_{\ell}'(t)\overline{\gamma_{\ell}(t)}\bv
\right)\Q_{\gamma_{\ell}(t)}(T)^{-1}\bv\,dt\\
&-\sum_{\ell=1}^N\int_{0}^{1}2\Re\Big(f(\gamma_{\ell}(t))(-\I)\gamma_{\ell}'(t)\Big)T\Q_{\gamma_{\ell}(t)}(T)^{-1}\bv\,dt.
\end{split}
\end{align}
Since $\Q_{\gamma_{\ell}(t)}(T)^{-1}\bv$ and $TQ_{\gamma_{\ell}(t)}(T)^{-1}\bv$  commute with real numbers, we furthermore have
\begin{align*}
& \int_{\partial(U\cap\cc_{\I})}f(s)\,ds_{\I}\,S_R^{-1}(s,T) \bv \\
=& \sum_{\ell=1}^N\int_{0}^{1}\Q_{\gamma_{\ell}(t)}(T)^{-1}\bv\,2\Re\left(f(\gamma_{\ell}(t))(-\I)\gamma_{\ell}'(t)\overline{\gamma_{\ell}(t)}
\right)\,dt\displaybreak[1]\\
&-\sum_{\ell=1}^N\int_{0}^{1}T\Q_{\gamma_{\ell}(t)}(T)^{-1}\bv\,2\Re\Big(f(\gamma_{\ell}(t))(-\I)\gamma_{\ell}'(t)\Big)\,dt \displaybreak[1]\\
=&\sum_{\ell=1}^{N}\int_{0}^{1}\left(\Q_{\gamma_{\ell}(t)}(T)^{-1}\bv\overline{\gamma_{\ell}(t)}-T\Q_{\gamma_{\ell}(t)}(T)^{-1}\bv\right)f(\gamma_{\ell}(t))\gamma_{\ell}'(t)\,dt (-\I)\\
&-\sum_{\ell=1}^{N}\int_{0}^{1}\left(\Q_{\gamma_{\ell}(t)}(T)^{-1}\bv\gamma_{\ell}(t) - T\Q_{\gamma_{\ell}(t)}(T)^{-1}\bv\right) \overline{f(\gamma_{\ell}(t))}\!\overline{\gamma_{\ell}'(t)}\,dt (-\I).
\end{align*}
Recalling again that  $f(\overline{x}) = \overline{f(x)}$ because $f$ is intrinsic, that $\Q_{\overline{s}}(T)^{-1} = \Q_{s}(T)^{-1}$ for $s\in\rho_S(T)$ and that $(-\overline{\gamma_{\ell}})(t) = -\overline{\gamma_{\ell}'(1-t)}$, we thus find
\begin{align*}
& \int_{\partial(U\cap\cc_{\I})}f(s)\,ds_{\I}\,S_R^{-1}(s,T) \bv \\
=&\sum_{\ell=1}^{N}\int_{\gamma_{\ell}}\left(\Q_{z}(T)^{-1}\bv\overline{z}-T\Q_{z}(T)^{-1}\bv\right)f(z)\,dz (-\I)\\
&+\sum_{\ell=1}^{N}\int_{-\overline{\gamma_{\ell}}}\left(\Q_{z}(T)^{-1}\bv\overline{z} - T\Q_{z}(T)^{-1}\bv\right) f(z)\,dt (-\I)\\
=&\int_{\partial(U\cap\cc_{\I})}\left(\Q_{z}(T)^{-1}\bv\overline{z}-T\Q_{z}(T)^{-1}\bv\right)f(z)\,dz (-\I)\\
=&\int_{\partial(U\cap\cc_{\I})}\RRes_{z}(T;\bv)f(z)\,dz (-\I).
\end{align*}
Finally, observe that $f(\infty) = \lim_{s\to\infty} f(s) \in\rr$ because as an intrinsic function $f$ takes only real values on the real line. Since any vector commutes with real numbers, we can hence rewrite \eqref{XX} as
\[ 
f(T)\bv = \bv f(\infty) + \int_{\partial(U\cap\cc_{\I})}\RRes_{z}(T;\bv)f(z)\,dz \frac{(-\I)}{2\pi }.
\]

\end{proof}
The above theorem shows once again that the $S$-functional calculus is the proper generalization of the holomorphic Riesz-Dunford functional calculus. Indeed, combining it with \Cref{ABC} reveals another deep relation between these two techniques.
\begin{corollary}
Let $T\in\closOP(V)$, let $f\in\intrin(\sigma_{SX}(T))$ and let $\I\in\SS$. Applying the $S$-functional calculus for $T$ to $f$ or considering $T$ as a $\cc_{\I}$-linear operator on $V_{\I}$ and applying the holomorphic Riesz-Dunford functional calculus to $f_{\I} := f|_{\fdom(f)\cap\cc_{\I}}$ yields the same operator. Both techniques are compatible.
\end{corollary}
Another important observation is the independence of the $S$-functional calculus of the left multiplication if $f$ is intrinsic.
\begin{corollary}\label{IndCol}
Let $V_{R}$ be a quaternionic right Banach space and let $T\in\closOP(V_R)$. Let $f\in\intrin(\sigma_{SX}(T))$ and assume that it is possible to endow $V_R$ with two different left scalar multiplications that turn it into a two-sided quaternionic Banach space. We denote these two-sided Banach spaces by $V_1$ and $V_2$ and we denote the operators obtained by applying the $S$-functional calculus of $T$ on $V_{1}$ resp. $V_2$ to $f$  by  $[f(T)]_1$ and $[f(T)]_2$. We then have
\[ [f(T)]_1\bv = [f(T)]_2\bv\quad\forall \bv\in V_{R} = V_1 = V_2.\] 
\end{corollary}
\begin{proof}
The operators $[f(T)]_1$ and $[f(T)]_2$ can both be represented by \eqref{IntSCalc}. This formula does however not involve any multiplication with scalars from the left such that we obtain the statement.

\end{proof}
\begin{remark}
We point out again that the representation \eqref{IntSCalc} and also \Cref{IndCol} only hold for intrinsic functions. Indeed, the symmetry $f(\overline{s}) = \overline{f(s)}$, which is satisfied only by intrinsic functions, is crucial in the proof of \Cref{IntCalcThm}. For left or right slice-hyperholomorphic functions that are not intrinsic, the operator $f(T)$ will in general depend on the left multiplication. Any left slice hyperholomorphic function $f$ can for instance be written as 
\[
f(s) = f_0(s) + f_1(s)e_{1} + f_2(s)e_{2} + f_3(s)e_{3}
\]
 with intrinsic slice hyperholomorphic components $f_{\ell},\ell = 0,\ldots,3$, where $e_{\ell},\ell = 1,2,3$ is the generating basis of the quaternions. The operators $f_{\ell}(T),\ell = 0,\ldots,3$ are then determined by the right multiplication on the space, but 
\[
f(T) = f_0(T) + f_1(T)e_{1} + f_2(T)e_{2} + f_3(T)e_{3}
\]
  obviously depends on how the imaginary units $e_{\ell},\ell = 1,2,3$ are multiplied onto vectors from the left.
\end{remark}

We observe once again, that the formula \eqref{IntSCalc} does only contain multiplications of operators with real numbers and multiplications of vectors with quaternionic scalars from the right---operations that are also available on a right Banach space---but it does not contain any operation that requires a two-sided Banach space, i.e. neither multiplications of vectors with scalars from the left nor multiplications of operators with non-real scalars. We can hence use this formula to define the $S$-functional calculus for operators on right Banach spaces.
\begin{definition}\label{NewSCalc}
Let $V_R$ be a quaternionic right Banach space and let $T\in\closOP(V_R)$. For $f\in\intrin(\sigma_{SX}(T))$, we  define
\[f(T)\bv := \bv f(\infty) + \int_{\partial(U\cap\cc_{\I})} \RRes_{s}(T;\bv)f(z)\,dz \frac{-\I}{2\pi},\]
where $\I\in\SS$ is an arbitrary imaginary unit and $U$ is an unbounded slice Cauchy domain with $\sigma_S(T)\subset U$ and $\clos{U}\subset \fdom(f)$. 
\end{definition}
The operator is well-defined. Cauchy's integral theorem guarantees the independence of the choice of $U$ and if we perform the computations in the proof of \Cref{IntCalcThm} in the inverse order, we arrive at the representation \eqref{KuKA} for $f(T)$. This formula does not contain any imaginary units but only real scalars so that $f(T)$ is also independent of the choice of $\I\in\SS$. In particular this representation also guarantees that the obtained operator is again quaternionic right linear. Moreover, because of \Cref{ABC}, the operator $f(T)$ coincides with the operator that we obtain if we consider $T$ as a closed $\cc_{\I}$-linear operator on $V_{R,\I}$ and apply the Riesz-Dunford functional calculus to $f_{\I} = f|_{\fdom(f)\cap\cc_{\I}}$. The following lemma shows that the intrinsic $S$-functional calculus has properties analogue to those in \Cref{SCProp} and \Cref{SCTaylor}.

\begin{lemma}\label{NewSCalcProp}
Let $T\in\closOP(V_R)$.
\begin{enumerate}[(i)]
\item \label{Eins}We have $(af + g)(T) = af(T) + g(T)$ and $(fg)(T) = f(T)g(T)$ for all $f,g\in\intrin(\sigma_{SX}(T))$ and all $a\in\rr$.
\item \label{Zwei} For any  $f(T)\in\intrin(\sigma_S(T))$, the operator $f(T)$ commutes with  $T$ and moreover also with any bounded operator $A\in\boundOP(V_R)$ that commutes with $T$. 
\item \label{Drei} The spectral mapping theorem holds
\begin{equation}\label{YYyy}
f(\sigma_{SX}(T)) = \sigma_{S}(f(T)) =  \sigma_{SX}(f(T))
\end{equation}
for all $f\in\intrin(\sigma_{SX}(T))$. Moreover, if $g\in\intrin_{\sigma_{SX}(T)}$, then $(g\circ f)(T) = g(f(T))$. 
\item \label{Vier}  If $\Delta\subset \sigma_{SX}(T)$ is a spectral set (i.e. open and closed in $\sigma_{SX}(T)$), then let $\chi_{\Delta}\in\intrin(\sigma_{SX}(T))$ be equal to $1$ on a neighborhood of $\Delta$ and equal to  $0$ on a neighborhood of $\sigma_{SX}(T)\setminus \Delta$. The operator $\chi_{\Delta}(T)$ is a continuous projection that commutes with $T$ and the right linear subspace $V_{\Delta}:=\ran\chi_{\Delta}(T)$ of $V_{R}$  is invariant under $T$. Finally, if we denote $T_{\Delta} := T|_{V_{\Delta}}$, then $\sigma_{S}(T) = \Delta$ and
\begin{equation}\label{Subsp}
 f(T_{\Delta}) = f(T)|_{\Delta}.
\end{equation} 
\item \label{Funf} Assume that $T$ is bounded and assume that $N\in\boundOP(V_R)$ is another bounded operator that commutes with $T$. If $\sigma_{S}(N)\subset \clos{B_{\varepsilon}(0)}$ for some $\varepsilon >0$, then 
\[
\sigma_{S}(T+N) \subset \clos{B_{\varepsilon}(\sigma_S(T))},
\]
where 
\[
B_{\varepsilon}(\sigma_S(T)) = \{s\in\hh: \dist(s,\sigma_S(T)) <\varepsilon\}.
\]
 For any $f\in\intrin\Big(\clos{B_{\varepsilon}(\sigma_{S}(T))}\Big)$, we furthermore have
\[
f(T) = \sum_{n=0}^{+\infty}N^n \frac{1}{n!}(\sderiv^n f)(T) = \sum_{n=0}^{+\infty}\left( \frac{1}{n!}(\sderiv^n f)(T)\right)N^n,
\]
where this series converges in the operator norm.
\end{enumerate}
\end{lemma}
\begin{proof}For neatness let us denote by $f_{\I}[T]$ the operator obtained from the Riesz-Dunford functional calculus considering $T$ as a $\cc_{\I}$-linear operator on $V_{R,\I}$ with $\I\in\SS$. 
The class of intrinsic functions is closed under addition, pointwise multiplication and multiplication with real numbers. Since moreover $f(T) = f_{\I}[T]$, the properties in \cref{Eins} and \cref{Zwei} follow immediately from the properties of the Riesz-Dunford-functional calculus as
\[ 
(af + g)(T) = (af+g)_{\I}[T] = a f_{\I}[T] + g_{\I}[T] = af(T) + g(T) 
\]
as well as 
\[
(fg)(T) = (fg)_{\I}[T] = f_{\I}[T] g_{\I}[T] =f(T)g(T)
\]
and
\[
f(T)T = f_{\I}(T) T \subset Tf_{\I}(T) = Tf(T).
\] 
Any operator  $A\in\boundOP(V_{R})$ that commutes with $T$ is also a bounded $\cc_{\I}$-linear operator on $V_{R,\I}$ that commutes with $T$ and hence the consistency with the Riesz-Dunford functional calculus yields again
\[
f(T) A  = f_{\I}[T]A = Af_{\I}[T] = A f(T).
\]

Since $f$ is intrinsic, we have $f(\dom(f)\cap\cc_{\I})\subset \cc_{\I}$ for all $\I\in\SS$ and $f([s]) = [f(s)]$. The spectral mapping theorem in \cref{Drei} holds as because of \Cref{ABC}
\[
 \sigma_{S}(f(T))\cap\cc_{\I} = \sigma_{\cc_{\I}}(f(T)) = \sigma(f_{\I}[T]) = f_{\I}(\sigma_{\cc_{\I}X}(T)) = f(\sigma_{SX}(T))\cap\cc_{\I}.
 \]
Taking the axially symmetric hull yields \eqref{YYyy}. If $g\in\intrin(f(\sigma_{S}(T)))$, we thus also have $g_{\I}\in\sigma(f_{\I}[T])$ and  
\[
 (g\circ f)(T) = (g\circ f)_{\I} [T] = (g_{\I}\circ f_{\I})[T]= g_{\I}[f_{\I}[T]]=  g(f(T)). 
\]

If $\Delta$ is a spectral set, then $\chi_{\Delta}$ is a projection that commutes with $T$ by what we have shown above as $\chi_{\Delta}(T)^2 = \chi_{\Delta}^2(T) = \chi_{\Delta}(T)$. Moreover, $\Delta_{\I} := \Delta\cap\cc_{\I}$ is a spectral set of $T$ as a $\cc_i$-linear operator on $V_{R,\I}$ and $(\chi_{\Delta})_{\I} = \chi_{\Delta_{\I}}$. Hence, $\chi_{\Delta}(T) = \chi_{\Delta_{\I}}[T]$ and so $V_{\Delta} = \ran\chi_{\Delta}(T) = \ran \chi_{\Delta_{\I}}[T]$. The properties of the Riesz-Dunford functional calculus  imply
\[
 \sigma_{S}(T_{\Delta})\cap\cc_{\I} = \sigma_{\cc_{\I}}(T_{\Delta}) = \Delta_{\I} = \Delta\cap\cc_{\I}.
 \]
 Taking the axially symmetric hull yields $\sigma_{S}(T_{\Delta}) = \Delta$. Furthermore
\[ 
f(T_{\Delta}) = f_{\I}[T_{\Delta}] = f_{\I}[T]|_{\ran \chi_{\Delta_{\I}}[T]} = f_{\I}[T]|_{\ran \chi_{\Delta}(T)} = f(T)|_{V_{\Delta}}.
\]

Let finally $T$ be  bounded and let $N\in\boundOP(V_R)$ with $\sigma_{S}(N)\subset \clos{B_{\varepsilon}(0)}$ commute with $T$. By \Cref{ABC} we have $\sigma_{\cc_{\I}}(N) \subset \clos{B_{\varepsilon}(0)\cap\cc_{\I}}$ and hence the properties of the Riesz-Dunford-functional calculus (precisely Theorem~10 in \cite[Chapter VII.6]{Dunford:1958}) imply
\begin{gather*}
\sigma_{S}(T + N) \cap\cc_{\I} = \sigma_{\cc_{\I}}(T+N) \\
\subset \left\{ z\in\cc_{\I}: \dist\left(z, \sigma_{\cc_{\I}}(T)\right) \leq \varepsilon\right\} = \clos{B_{\varepsilon}(\sigma_{S}(T))}\cap\cc_{\I}
\end{gather*}
and 
\[
f_{\I}[T+N] = \sum_{n=0}^{+\infty} N^n \frac{1}{n!}\left(f_{\I}^{(n)}\right)[T],
\]
where this series converges in the operator norm. We therefore find
\[
\sigma_{S}(T+N) = \left[\sigma_{S}(T + N) \cap\cc_{\I}\right] =\left[ \clos{B_{\varepsilon}(\sigma_{S}(T))}\cap\cc_{\I} \right] = \clos{B_{\varepsilon}(\sigma_{S}(T))}
\]
and, using the fact that $(\sderiv f)_{\I} = (f_{\I})'$, we also find that
\[
f(T+N) = f_{\I}[T+N] = \sum_{n=0}^{+\infty} N^n \frac{1}{n!} f_{\I}^{(n)}[T] =  \sum_{n=0}^{+\infty} N^n \frac{1}{n!}(\sderiv^n f)(T).
\]

\end{proof}
\begin{remark}
We point out that \eqref{Subsp} can not be shown as a property of the $S$-functional calculus with the usual approach: since invariant subspaces of right linear operators are only right linear subspace, it is in general not possible to define $f(T_{\Delta})$. This would require $V_{\Delta}$ to be a two-sided subspace of the Banach space. This fact might cause technical difficulties as it happened for instance in the proof of the spectral mapping theorem in \cite{HInftyComing}. The approach in \Cref{NewSCalc} provides a work-around for this problem.
\end{remark}

\section{Spectral integration in the quaternionic setting}\label{SpecIntSect}
Before we pass to the study of quaternionic spectral operators, let us discuss spectral integration in the quaternionic setting. The basic idea of spectral integration is well known: it generates a multiplication operator in a way that generalizes the multiplication with eigenvalues in the discrete case. If for instance $\lambda \in \sigma(A)$ of some normal operator $A:\cc^{n}\to\cc^{n}$, then we can define $E(\{\lambda\})$ to be the orthogonal projection of $\cc^{n}$ onto the eigenspace associated to $\lambda$ and we find $A = \sum_{\lambda \in\sigma(A)}\lambda \,E(\{\lambda\})$. Setting $E(\Delta) = \sum_{\lambda \in \Delta}E(\{\lambda\})$ one obtains a discrete measure on $\cc$, the values of which are orthogonal projections on $\cc^n$, and $A$ is the integral of the identity function with respect to this measure. Changing the notation accordingly we have
\begin{equation}\label{TIntC}
A= \sum_{\lambda \in \sigma(A)} \lambda \,E(\{\lambda\})\qquad \Longrightarrow\qquad A = \int_{\sigma(A)} \lambda\, dE(\lambda).
\end{equation}
Via functional calculus it is possible to define functions of an operator. The fundamental intuition of a functional calculus is that $f(A)$ should be defined by the action of $f$ on the spectral values (resp. the eigenvalues) of $A$. For our normal operator $A$ on $\cc^n$ the operator $f(A)$ is the operator with the following property: if $\bv\in\cc^n$ is an eigenvector of $A$ with respect to $\lambda$, then $\bv$ is an eigenvector of $f(A)$ with respect to $f(\lambda)$, just as it happens for instance naturally for powers and polynomials of $A$. Using the above notation we thus have
\begin{equation}\label{f(T)IntC}
f(A) = \sum_{\lambda\in\sigma(A)}f(\lambda)\,E(\{\lambda\})\qquad\Longrightarrow\qquad f(A) = \int_{\sigma(A)}f(\lambda)\,dE(\lambda).
\end{equation}

In infinite dimensional Hilbert spaces, the spectrum of a normal operator might be not discrete such that the expressions on the left-hand side of \eqref{TIntC} and \eqref{f(T)IntC} do not make sense. If $E$ however is a suitable projection-valued measure, then it is possible to give the expression \eqref{f(T)IntC} a meaning by following the usual way of defining integrals via the approximation of $f$ by simple functions. The spectral theorem then shows that for every normal operator $T$ there exists a spectral measure such that \eqref{TIntC} holds true.

If we want to introduce similar concepts in the quaternionic setting, we are---even in the discrete case---faced with several unexpected phenomena.
\begin{enumerate}[{(P1)}]
\item\label{p1} The space of bounded linear operators on a quaternionic Banach space $V_{R}$ is only a real Banach space and not a quaternionic one as pointed out in \Cref{B(V)Struct}. Hence, the expressions in \eqref{TIntC} and \eqref{f(T)IntC} are a priori only defined if $\lambda$ resp. $f(\lambda)$ are real. Otherwise one needs to give meaning to the multiplication of the operator  $E(\{\lambda\})$ with nonreal scalars.
\item \label{p2} The multiplication with a (non-real) scalar on the right is not linear such that $a E(\{\lambda\})$ can for $a\in\hh$ not be defined as $(a E(\{\lambda\})a)\bv = (E(\{\lambda\})\bv )a $. Moreover, the set of eigenvectors associated to a specific eigenvalue does not constitute a linear subspace of $V_{R}$: if for instance $T\bv = \bv s$ with $s = s_0 + \I_s s_1$ and $\J\in\SS$ with $\I_s\perp \J$, then $T(\bv\J) = (T\bv)\J = (\bv s)\J = (\bv\J)\overline{s}\neq (\bv\J)s$. 
\item \label{p3} Finally, the set of eigenvalues is in general not discrete: if $s \in\hh$ is an eigenvalue of $T$ with $T\bv = \bv s$ for some $\bv\neq \bO$ and $s_{\I} = s_0 + \I s_1\in[s]$, then there exists $h \in\hh\setminus\{0\}$ such that $s_{\I} = h^{-1}sh$ and so
\begin{equation}\label{EVCo}
 T(\bv h) = T(\bv)h = \bv sh = (\bv h)h^{-1}sh = (\bv h)s_{\I}.
 \end{equation}
Thus, any $s_{\I}\in[s]$ is also an eigenvalue of $T$.  
\end{enumerate}

As a first consequence of \Cref{p2,p3} the notion of eigenvalue and eigen\-space have to be adapted: linear subspaces are in the quaternionic setting not associated with individual eigenvalues $s$ but with spheres~$[s]$ of equivalent eigenvalues.
\begin{definition}
Let $T\in\closOP(V_{R})$ and let $s\in\hh\setminus\rr$. We say that $[s]$ is an eigensphere of $T$ if there exists a vector $\bv\in V_{R}\setminus\{\bO\}$ such that 
\begin{equation}\label{SEigenVal}
(T^2 - 2s_0T + |s|^2\id)\bv = \Q_{s}(T)\bv= \bO.
\end{equation}
The eigenspace of $T$ associated with $[s]$ consists of all those vectors that satisfy \eqref{SEigenVal}.
\end{definition}
\begin{remark}
For real values, things remain as we know them from the classical case: a quaternion $s\in\rr$ is an eigenvalues of $T$ if $T\bv - \bv s = 0$ for some $\bv\neq\bO$. The quaternionic right linear subspace $\ker ( T- s\id)$ is then called the eigenspace of $T$ associated with $s$. 
\end{remark}
Any eigenvector $\bv$ that satisfies $T(\bv) = \bv s_{\I}$ with $s_{\I} = s_0 + \I s_1$ for some $\I\in\SS$ belongs to the eigenspace associated with the eigensphere $[s]$. Note however that the eigenspace associated with an eigensphere $[s]$ does not only consist of eigenvectors. It contains also linear combinations of eigenvectors associated with different eigenvalues in $[s]$ as the next lemma shows.
\begin{lemma}\label{EVSplitLem}
Let $T\in\closOP(V_R)$ and let $[s]$ be an eigensphere of $T$ and let $\I\in\SS$.  A vector $\bv$ belongs to the eigenspace associated with  $[s]$ if and only if $ \bv = \bv_1 + \bv_2$ such that $T\bv_1 = \bv_1 s_{\I}$ and  $T \bv_2 = \bv_2\overline{s_{\I}}$ where $s_{\I} = s_0 + \I s_1$.
\end{lemma}
\begin{proof}
Observe that
\begin{equation}\label{as}
\Q_{s}(T)\bv = T^2\bv - T\bv2s_0  +  \bv|s|^2 = T(T\bv - \bv\overline{s_{\I}}) - (T\bv - \bv\overline{s_{\I}})s_{\I} 
\end{equation}
and
\begin{equation}\label{bs}
\Q_{s}(T)\bv = T^2\bv -  T\bv 2s_0 +  \bv |s|^2 = T(T\bv - \bv s_{\I}) - (T\bv-\bv s_{\I})\overline{s_{\I}}.
\end{equation}
Hence, $\Q_{s}(T)\bv = \bO$ for any eigenvector associated with $s_{\I}$ or $\overline{s_{\I}}$ and in turn also for any $\bv$ that is the sum of two such vectors. 

If conversely $\bv$ satisfies \eqref{SEigenVal}, then we deduce from \eqref{as}, that $T\bv-\bv\overline{s_{\I}}$ is a right eigenvector associated with $s_{\I}$ and that $T\bv - \bv s_{\I}$ is a right eigenvalues of $T$ associated with $\overline{s_{\I}}$. Since $s_{\I}$  and $\I$ commute, the vectors $\bv_1 = (T\bv - \bv\overline{s_{\I}})\frac{-\I}{2s_1}$ and $\bv_{2} = (T\bv - \bv s_{\I})\frac{\I}{2s_1}$ are right eigenvectors associated with $s$ resp. $\overline{s_{\I}}$, too. Hence we found the desired decomposition as
\begin{align*}
\bv_1 + \bv_2 = (T\bv - \bv\overline{s_{\I}})\frac{-\I}{2s_1}+ (T\bv - \bv s_{\I})\frac{\I}{2s_1} = \bv (\overline{s_{\I}} - s_{\I}) \frac{\I}{2s_1} = \bv .
\end{align*}

\end{proof}
\begin{remark}
If $\J\in\SS$ with $\J\perp\I$, then $\tilde{\bv}_2 := \bv_2(-\J)$ is an eigenvector of $T$ associated with $s$. Hence we can write $\bv$ also as $\bv = \bv_{1} + \tilde{\bv}_{2}\J$, where $\bv_1,\tilde{\bv}_2$ are both eigenvectors associated with $s_{\I}$. 
\end{remark}

Since the eigenspaces of quaternionic linear operators are not associated with individual eigenvalues but instead with eigenspheres, quaternionic spectral measures must not be defined on arbitrary subsets of the quaternions. Instead their natural domains of definition consist of axially symmetric subsets of $\hh$ so that they associate subspaces of $V_{R}$ not to sets of spectral values but to sets of spectral spheres. This is also consistent with the fact that the $S$-spectrum of an operator is axially symmetric.
\begin{definition}
We denote the sigma-algebra of axially symmetric Borel sets on $\hh$  by $\sBorel(\hh)$. Furthermore, we denote the set of all real-valued $\sBorel(\hh)$-$\Borel(\rr)$-measurable functions defined on $\hh$ by $\sMeas(\hh,\rr)$  and the set of all such functions that are bounded by $\bsMeas(\hh,\rr)$.
\end{definition}
\begin{remark}
The restrictions of functions in $\bsMeas(\hh,\rr)$ to a complex halfplane $\cc_{\I}^{\geq}$ are exactly the functions that were used to construct the spectral measure of a quaternionic normal operator in~\cite{Alpay:2016}, cf. \Cref{RemSpecThm}. 
\end{remark}

\begin{definition}\label{SpecMeas}
A quaternionic spectral measure on a quaternionic right Banach space $V_{R}$ is a function $E:\sBorel(\hh)\to\boundOP(V_{R})$ that satisfies
\begin{enumerate}[(i)]
\item \label{SM1} $E(\Delta)$ is a continuous projection  and $\|E(\Delta)\| \leq K$  for any $\Delta\in\sBorel(\hh)$ with a constant $K>0$ independent of $\Delta$,
\item \label{SM2} $E(\emptyset) = 0$ and $E(\hh) = \id$,
\item \label{SM3}  $E(\Delta_1\cap\Delta_2) = E(\Delta_1)E(\Delta_2)$  for any $\Delta_1,\Delta_2\in\sBorel(\hh)$  and
\item  \label{SM4}for any sequence $(\Delta_n)_{n\in\nn}$ of pairwise disjoint sets in $\sBorel(\hh)$ we have 
\(E\left  (\bigcup_{n\in\nn} \Delta_n\right)\bv = \sum_{n=1}^{+\infty}E(\Delta_n)\bv\) for all $\bv\in V_{R}$.
\end{enumerate}
\end{definition}
\begin{corollary}
Let $E$ be a spectral measure on $V_R$ and let $V_{R}^{*}$ be its dual space, the left Banach space consisting of all continuous right linear mappings from $V_{R}$ to $\hh$. For any $\bv\in V_R$ and any $\bv^{\mathbf{*}}\in V_R^*$, the mapping $\Delta \mapsto \langle \bv^{\mathbf{*}}, E(\Delta)\bv\rangle$ is a quaternion-valued measure on $\sBorel(\hh)$.
\end{corollary}
\begin{remark}\label{SpecMeasOthers}
Following \cite{Viswanath:1971}, most authors considered quaternionic spectral measures that are defined on the Borel sets $\Borel\big(\cc_{\I}^{\geq}\big)$ of a closed complex halfplane $\cc_{\I}^{\geq}:= \{s_0 + \I s_1: s_0\in\rr, s_1\geq 0\}$. This is equivalent to $E$ being defined on $\sBorel(\hh)$. Indeed, if $\tilde{E}$ is defined on $\Borel\big(\cc_{\I}^{\geq}\big)$, then setting
\[ E(\Delta) : = \tilde{E}\big(\Delta\cap\cc_{\I}^{\geq}\big)\qquad\forall\Delta\in\sBorel(\hh)\]
yields a spectral measure in the sense of \Cref{SpecMeas} that is defined on $\sBorel(\hh)$. If on the other hand we start with a spectral measure $E$ defined on $\sBorel(\hh)$, then setting
\[
\tilde{E}(\Delta) := E([\Delta])\qquad\forall\Delta\in\Borel\big(\cc_{\I}^{\geq}\big)
\]
yields the respective measure on $\Borel\big(\cc_{\I}^{\geq}\big)$. Although both definitions are equivalent, we prefer $\sBorel(\hh)$ as the domain of $E$ because it does not suggest a dependence on the imaginary unit $\I$. 
\end{remark}
For a function $f\in\bsMeas(\hh,\rr)$, we can now define  the spectral integral with respect to a spectral measure $E$ as in the classical case \cite{Dunford:1958,Dunford:1958}. If $f$ is a simple function, i.e. $f(s) = \sum_{k=1}^{n}\alpha_{k}\chi_{\Delta_k}(s)$ with pairwise disjoint sets $\Delta_{k}\in\sBorel(\hh)$, where $\chi_{\Delta_k}$ denotes the characteristic function of $\Delta_{k}$, then we set
\begin{equation}\label{SimpleSpecInt}
\int_{\hh} f(s)\,dE(s) : = \sum_{k=1}^{n}\alpha_{k}E(\Delta_k).
\end{equation}
There exists a constant $C_{E}>0$ that depends only on $E$ such that
\begin{equation}\label{NormEst1}
\left\| \int_{\hh} f(s)\,dE(s) \right\| \leq C_{E}\|f\|_{\infty},
\end{equation}
where $\| .\|_{\infty}$ denotes the supremum-norm. If $f\in\bsMeas(\hh,\rr)$ is arbitrary, then we can find a sequence of simple functions $(f_n)_{n\in\nn}$ such that $\| f - f_n\|_{\infty} \to 0$ as $n\to+ \infty$. In this case we can set
\begin{equation}\label{RealSpecInt}
 \int_{\hh} f(s) \, dE(s) : = \lim_{n\to+\infty} \int_{\hh} f_n(s) \, dE(s),
 \end{equation}
where this sequence converges in the operator norm because of \eqref{NormEst1}.

\begin{lemma}\label{HomLem}
Let $E$ be a quaternionic spectral measure on $V_{R}$. The mapping $f\mapsto \int_{\hh} f(s)\, dE(s)$ is a continuous homomorphism from $\bsMeas(\hh,\rr)$ to $\boundOP(V_{R})$. Moreover, if $T$ commutes with $E$, i.e. it satisfies $TE(\Delta) = E(\Delta) T$ for all $\Delta\in\sBorel(\hh)$, then $T$ commutes with $\int_{\hh} f(s)\, dE(s)$ for any $f\in\bsMeas(\hh,\rr)$.
\end{lemma}
\begin{corollary}\label{DualRep}
Let $E$ be a quaternionic spectral measure on $V_{R}$ and let $f\in\bsMeas(\hh,\rr)$. For any $\bv\in V_R$ and any $\bv^{\mathbf{*}} \in V_R^*$, we have
\[ \left\langle \bv^{\mathbf{*}}, \left[\int_{\hh} f\, dE\right] \bv\right\rangle = \int_{\hh} f(s) \, d\langle \bv^{\mathbf{*}}, E(s)\bv \rangle .\]
\end{corollary}
\begin{proof}
Let $f_n = \sum_{k=1}^{N_{n}} \alpha_{n,k}\chi_{\Delta_{n,k}}\in\bsMeas(\hh,\rr)$ be such that $\|f - f_{n}\| \to 0$ as $n\to+ \infty$. Since all coefficients $\alpha_{n,k}$ are real, we have
\begin{gather*}
 \left\langle \bv^{\mathbf{*}}, \left[\int_{\hh} f\, dE\right] \bv\right\rangle = \lim_{n\to\infty} \left\langle \bv^{\mathbf{*}}, \left[\sum_{k=1}^{N_n}\alpha_{n,k}E(\Delta_{n,k})\right] \bv\right\rangle \\
= \lim_{n\to\infty} \sum_{k=1}^{N_n}\alpha_{n,k}  \left\langle \bv^{\mathbf{*}},E(\Delta_{n,k}) \bv\right\rangle =\int_{\hh} f(s) \, d\langle \bv^{\mathbf{*}}, E(s)\bv \rangle .
 \end{gather*}
\end{proof}
\begin{remark}
The fact that the above definitions are well-posed and the properties given in \Cref{HomLem} can be shown as in the classical case, so we omit their proofs. One can also deduce them directly from the classical theory: if we consider $V_R$ as a real Banach space and $E$ as a spectral measure with values in the space $\boundOP_{\rr}(V_{R})$ of bounded $\rr$-linear operators on $V_{R}$, which obviously contains $\boundOP(V_R)$, then $\int_{\hh} f(s) \, dE(s)$ defined in \eqref{SimpleSpecInt} resp. \eqref{RealSpecInt} is nothing but the spectral integral of $f$ with respect to $E$ in the classical sense. Since any $\alpha_k$ in \eqref{SimpleSpecInt} is real and since each $E(\Delta)$ is a quaternionic right linear projection, the integral of any simple function $f$ with respect to $E$ is  a quaternionic right linear operator and hence belongs to $\boundOP(V_R)$. The space $\boundOP(V_R)$ is closed in $\boundOP_{\rr}(V_R)$ and hence the property of being quaternionic linear survives the approximation by simple functions such that $\int_{\hh} f(s)\,dE(s)\in\boundOP(V_R)$ for any $f\in\bsMeas(\hh,\rr)$ even if we consider it as the integral with respect to a (real) spectral measure on the real Banach space $V_R$.
\end{remark}

The techniques introduced so far allow us to integrate real-valued functions with respect to a spectral measure. This is obviously insufficient, even for formulating the statement corresponding to \eqref{TIntC} in the quaternionic setting unless $\sigma_{S}(T)$ is real. In order to define spectral integrals for functions that are not real-valued, we need additional information.

This fits another observation: in contrast to the complex case, even for the simple case of a normal operator on a finite-dimensional quaternionic Hilbert space, a decomposition of the space $V_{R}$ into the eigenspaces of $T$ is not sufficient to recover the entire operator $T$. Let $\I,\J\in\SS$ with $\I\neq \J$ and consider for instance the operators $T_1$, $T_2$ and $T_3$ on $\hh^2$, which are given by their matrix representations 
\begin{equation}
T_1 = \begin{pmatrix} \I & 0\\0&\I\end{pmatrix}\qquad T_2 = \begin{pmatrix} \I & 0\\0&\J\end{pmatrix}\qquad T_3 = \begin{pmatrix} \J& 0\\0&\J\end{pmatrix}.
\end{equation}
For each of these operators, we have $\sigma_S(T_{\ell}) = \SS$ and that its only eigenspace is the entire space $\hh^2$. The spectral measure $E$ that is associated with $T_{\ell}$ is hence given by $E(\Delta) = 0$ if $\SS\not\subset\Delta$ and $E(\Delta) = \id$ if $\SS\subset\Delta$. Hence, the spectral measures associated with these operators agree, although these operators do obviously not the coincide. 

Since the eigenspace of an operator $T$ that is associated with some eigensphere $[s]$ contains eigenvectors associated with different eigenvalues, we need some additional information to understand \lq how to multiply the eigensphere onto the associated eigenspace\rq, i.e. to understand which vector in the eigenspace must be multiplied with which eigenvalue in the corresponding eigensphere  $[s]$. This information will be provided by a suitable imaginary operator. Such operators generalize the properties of the anti-selfadjoint partially unitary operator $\mathsf{J}_0$ in the decomposition \eqref{SplitTeasy} of a normal operator on a Hilbert space to the Banach space setting.

\begin{definition}\label{ImOp}
An operator $J\in\boundOP(V_{R})$ is called imaginary if  $-J^2$ is the projection of $V_R$ onto $\ran J$ along $\ker J$. We call $J$ fully imaginary if $-J^2 = \id$, i.e. if in addition $\ker J = \{0\}$.
\end{definition}

\begin{corollary} \label{ImOpCor1}
An operator $J\in\boundOP(V_{R})$ is an imaginary operator if and only if
\begin{enumerate}[(i)]
\item \label{ImOp1} $-J^2$ is a projection and 
\item \label{ImOp2} $\ker J = \ker J^2$. 
\end{enumerate}
\end{corollary}
\begin{proof}
If $J$ is an imaginary operator, then obviously \cref{ImOp1} and \cref{ImOp2} hold true. Assume on the other hand that \cref{ImOp1} and \cref{ImOp2} hold. Obviously $\ran(-J^2) \subset \ran J$. For any $\bu\in V_{R}$, we have $(-J^2)\bu - \bu \in \ker (-J^2)  = \ker J$ because
\[
(-J^2)\left((-J^2)\bu - \bu\right) = (-J^2)^2\bu - (-J)^2\bu = (-J^2)\bu - (-J)^2\bu = \bO
\]
as  $(-J^2)$ is a projection. Therefore 
 \[
 \bO = J(-J^2\bu - \bu) = (-J^2)J\bu - J\bu
 \]
 and hence $\bv = (-J^2)\bv$ for any $\bv = J\bu \in\ran J$. Consequently, $\ran(-J^2) \supset \ran J$ and in turn $\ran J = \ran (-J^2)$. Since $\ker J = \ker (-J^2)$, we find that $-J^2$ is the projection of $V_{R}$ onto $\ran J$ along $\ker J$, i.e. that $J$ is an imaginary operator.

\end{proof}

\begin{remark}
The above implies that any anti-selfadjoint partially unitary operator $\mathsf{J}_0$  on a quaternionic Hilbert space $\hil$ is an imaginary operator. Indeed, for any $\bv\in\ker \mathsf{J}_{0}$, we obviously have $-\mathsf{J}_{0}^2\bv = \bO$. Since the restriction of $\mathsf{J}_{0}$ to $\hil_{0}:=\ran \mathsf{J}_{0} = \ker \mathsf{J}_{0}^{\perp}$ is unitary and $\mathsf{J}_{0}$ is anti-selfadjoint, we furthermore have for $\bv \in \hil_{0}$ that $- \mathsf{J}_{0}^2 \bv= \mathsf{J}_{0}^\ast \mathsf{J}_{0} \bv = \left(\mathsf{J}_{0}|_{\hil_{0}}\right)^{*}\left(\mathsf{J}_{0}|_{\hil_{0}}\right)\bv = \left(\mathsf{J}_{0}|_{\hil_{0}}\right)^{-1}\left(\mathsf{J}_{0}|_{\hil_{0}}\right)\bv = \bv$. Hence $-\mathsf{J}_{0}^2$ is the orthogonal projection onto $\hil_{0} = \ran \mathsf{J}_0$ and so $\mathsf{J}_{0}$ is an imagainary operator. In particular any unitary anti-selfadjoint operator is fully imaginary.
\end{remark}

\begin{lemma}\label{JSpec}
If $J\in\boundOP(V_R)$ is an imaginary operator, then $\sigma_S(T)\subset \{0\}\cup\{\SS\}$.
\end{lemma}
\begin{proof}
Since the operator $-J^2$ is a projection, its $S$-spectrum $\sigma_S(-J^2)$ is a subset of $\{0,1\}$. Indeed, for any projection $P\in\boundOP(V)$, a simple calculation shows that the pseudo-resolvent of $P$ at any $s\in\hh\setminus\{0,1\}$ is given by 
\[\Q_{s}(P)^{-1} = -\frac{1}{|s|^2}\left(\frac{1-2\Re(s)}{1-2\Re(s) + |s|^2}P - \id\right)\]
such that $s\in\rho_S(P)$. As a consequence of the spectral mapping theorem, we find that 
\[ -\sigma_S(J)^2 = \{-s^2:  s\in\sigma_S(J)\} = \sigma_S(-J^2) \subset \{0,1\}.\]
But if $-s^2 \in\{0,1\}$, then $s\in \{0\}\cup\SS$ and hence $\sigma_S(J)\subset \{0\}\cup\SS$.

\end{proof}
\begin{remark}
If $J = 0$, then $J$ is an imaginary operator with $\sigma_{S}(T) = \{0\}$. If on the other hand $\ker J = \{\bO\}$ (i.e. if $J$ is fully imaginary), then $\sigma_{S}(T) = \SS$. In any other case we obviously have $\sigma_{S}(T) = \{0\}\cup\SS$.
\end{remark}

The following theorem gives a complete characterization of imaginary operators on $V_{R}$. 
\begin{theorem}\label{JVSplitThm}
Let $J\in\boundOP(V_R)$ be an imaginary operator. For any $\I\in\SS$, the Banach space $V_R$ admits a direct sum decomposition as 
\begin{equation}\label{J0VSplit} 
V_R = V_{J,0} \oplus V_{J,\I}^{+} \oplus  V_{J,\I}^{-}
\end{equation}
with
\begin{equation}\label{VSplitSpaces}
\begin{split}
 V_{J,0} =& \ker (J), \\
  V_{J,\I}^{+} =& \{ \bv\in V: J\bv = \bv\I\},\\
 V_{J,\I}^{-} =& \{\bv\in V: J\bv = \bv(-\I)\}.
 \end{split}
 \end{equation}
The spaces $V_{J,\I}^{+}$ and $V_{J,\I}^{-}$ are complex Banach spaces over $\cc_{\I}$ with the natural structure inherited from $V_R$ and for each $\J\in\SS$ with $\I\perp \J$ the map $\bv\mapsto \bv \J$ is a $\cc_{\I}$-antilinear and isometric bijection between $V_{J,\I}^+$ and $V_{J,\I}^{-}$.

Conversely, let $\I,\J\in\SS$ with $\I\perp\J$ and assume that $V_{R}$ is the direct sum $ V_R = V_{0} \oplus V_{+} \oplus V_{-}$ of a closed ($\hh$-linear) subspace $V_{0}$ and two closed $\cc_{\I}$-linear subspaces  $V_{+}$ and $V_{-}$ of $V_{R}$  such that  $\Psi:\bv \mapsto \bv \J$ is a bijection between $V_{+}$ and $V_{-}$. Let $E_+$ and $E_{-}$ be the $\cc_{\I}$-linear projections onto $V_{+}$ and $V_{-}$ along $V_{0}\oplus V_{-}$ resp. $V_{0}\oplus V_{+}$. The operator $J\bv := E_{+}\bv\I + E_{-}\bv (-\I)$ for $\bv\in V_R$ is an imaginary operator on $V_{R}$.
\end{theorem}
\begin{proof}
We first assume that $J$ is an imaginary operator and show the existence of the corresponding decomposition of $V_{R}$. Let $\I\in\SS$ and $V_{R,\I}$ denote   the space $V_{R}$ considered as a complex Banach over $\cc_{\I}$. Furthermore, let us assume that $J\neq 0$ as the statement is obviously true in this case. Then $J$ is a bounded $\cc_{\I}$-linear operator on $V_{R,\I}$ and by \Cref{ABC} and \Cref{JSpec} the spectrum of $J$ as an element of $\boundOP(V_{R,\I})$ is $\sigma_{\cc_{\I}}(J) = \sigma_S(J)\cap\cc_{\I} \subset \{0,\I,-\I\}$. We define now for $\tau\in\{0,\I,-\I\}$ the projection $E_{\tau}$ as the spectral projection associated with $\{\tau\}$ obtained from the Riesz-Dunford functional calculus. If we choose $0<\varepsilon < \frac12$, then the relation $R_{z}(J) = (\overline{z}\id_{V_{R,\I}}-J)\Q_{z}(J)^{-1}$ in \Cref{ABC} implies
\[ E_{\tau} \bv =  \int_{\partial U_{\varepsilon}(\tau;\cc_{\I})}  R_{z}(J)\bv \, dz \frac{1}{2\pi i } =  \int_{\partial U_{\varepsilon}(\tau;\cc_{\I})}  \Q_{z}(J)^{-1}(\bv \overline{z}  - J\bv) \, dz \frac{1}{2\pi i } ,\]
where $ U_{\varepsilon}(\tau;\cc_{\I})$ denotes the ball of radius $\varepsilon$ in $\cc_{\I}$  that is centered at $\tau$. (Since we assumed $\ker J \neq V$, the projections $E_{\I}$ and $E_{-\I}$ are not trivial. It might however happen that $E_{0} = 0$, but this is not a problem in the following argumentation.)

We set
\[ V_{J,0} = E_{0} V_{R,\I},\quad V_{J,\I}^{+} = E_{\I}V_{R,\I}\quad \text{and} \quad V_{J,\I}^{-} = E_{-\I} V_{R,\I}.\]
Obviously these are closed $\cc_{\I}$-linear subspaces of $V_{R,\I}$ resp. $V_{R}$ and \eqref{J0VSplit} holds true.

Let us now show that the relation \eqref{VSplitSpaces} holds true. We first consider the subspace $V_{J,\I}^+$. Since it is the range of the Riesz-projector $E_{\I}$ associated with the spectral set $\{\I\}$, this is an $\cc_{\I}$-linear subspace of $V_{R,\I}$ that is  invariant under $J$ and the restriction $J_{+}:= J|_{V_{J,\I}^+}$ has spectrum $\sigma(J_{+}) = \{\I\}$.  Now observe that $-J_+^2 =-J^2|_{V_{J,\I}^+} $ is the restriction of a projection onto an invariant subspace and hence a projection itself. Since $0\notin\sigma(-J_{+}^2) =- \sigma(J_+)^2 = \{1\}$, we find $\ker -J_{+}^2 = \{\bO\}$ and in turn $\id_{+}:= \id_{V_{J,\I}^+} = -J^2_{+}$. For $\bv\in V_{J,\I}^+$ we therefore have
\begin{gather*}
-\bv = J_+^2 \bv = (J_+ - \I\id_{+} + \I\id_{+})^2\bv = (J_+ - \I\id_{+} + \I\id_{+})((J_+ - \I\id_{+})\bv + \bv\I)\\
= (J_+ - \I\id_+)^2\bv + (J_+ - \I\id_{+})\bv\I + (J_+ - \I\id_{+})\bv\I + \bv \I^2.
\end{gather*}
As $\I^2 = -1$ this is equivalent to
\[
 (J_+ - \I\id_+)^2 \bv = (J_+-\I\id_{+})\bv (-2\I).
\]
Hence $(J_{+}- \I\id_{+})\bv$ is either $\bO$ or an eigenvector of $J_{+} - \I\id_{+}$ associated with the eigenvalue $-2\I$. By the spectral mapping theorem $\sigma(J_+ - \I\id_{+}) = \sigma(J_+) - \I = \{0\}$. Hence, $J_+ - \I\id_+$ cannot have an eigenvector with respect to the eigenvalue $-2\I$ and so $(J_{+}-\I\id_{+})\bv = \bO$. Therefore $J_+ = \id_+ i$ and $J\bv = J_+\bv = \bv \I$ for all $\bv\in V_{J,\I}^+$. 

With similar arguments, one shows that $J\bv = \bv(-\I)$ for any $\bv \in V_{J,\I}^-$. Finally, $\sigma(-J_0^2) = - \sigma(J_0)^2 = \{0\}$  for $J_{0}:= J|_{V_{J,0}}$. Since $- J_0^2 = -J^2|_{V_{J,0}}$ is the restriction of a projection to an invariant subspace and thus a projection itself, we find that $-J_0^2$ is the zero operator and hence $V_{J,0}= \ker(-J_0)^2 \subset \ker(J^2) = \ker J$. On the other hand $\ker J \subset V_{J,0}$ as $V_{J,0}$ is the invariant subspace associated with the spectral value $0$ of $J$. Thus $V_{J,0} = \ker J$ and so \eqref{VSplitSpaces} is true.

Finally, if $\J\in\SS$ with $\I\perp \J$ and $\bv\in V_+$ then $(J \bv \J) = J(\bv) \J = \bv \I\J= (\bv \J)(-\I)$. Hence $\Psi:\bv \to \bv \J$ maps $V_{J,\I}^+$ to $V_{J,\I}^-$. It is obviously $\cc_{\I}$-antilinear, isometric and a bijection as  $\bv = - (\bv \J)\J$ so that the proof of the first statement is finished.

Now let $\I,\J\in\SS$ with $\I\perp\J$ and assume that $V_{R} = V_0 \oplus V_{+} \oplus V_{-}$  with subspaces $V_0$, $V_+$ and $V_{-}$ as in the assumptions. We define $J\bv := E_{+} \bv \I + E_{-}\bv(-\I)$. Obviously, $J$ is a continuous $\cc_{\I}$-linear operator on $V_{R,\I}$. The mapping $\Psi:\bv\mapsto \bv \J$ maps $V_{+}$ bijectively to $V_{-}$, but since $\Psi^{-1} = - \Psi$ it also maps $V_{-}$ bijectively to $V_{+}$. Moreover, as an $\hh$-linear subspace, $V_0$ is invariant under $\Psi$. For $\bv = \bv_{0} + \bv_{+} + \bv_{-} \in V_{0} \oplus V_{+} \oplus V_{-} = V_{R}$, we therefore find
\begin{align*}
 J(\bv \J) =& E_{+}(\bv\J)\I + E_{-}(\bv\J) (-\I)  = \bv_{-}\J\I + \bv_{+}\J(-\I)  \\
 =& \bv_{-}(-\I)\J + \bv_{+} \I \J = \left(E_{-}\bv(-\I)\right)\J + \left(E_{+}\bv\I\right)\J = (J\bv )\J.
\end{align*}
If now $a\in\hh$, then we can write $a = a_{1} +  a_2\J$ with $a_{1},a_{2}\in\cc_{\I}$ and find due to the $\cc_{\I}$-linearity of $J$ that
\begin{align*}
J(\bv a) = J (\bv a_1) + J(\bv a_2 \J) = J(\bv) a_1 + J(\bv) a_2\J = J(\bv) (a_1+  a _2 \J) = J(\bv) a.
\end{align*}
Hence, $J$ is quaternionic linear and therefore belongs to $\boundOP(V_{R})$. 

As $E_{+}E_{-} = E_{-}E_{+} = 0$, we furthermore observe that
\begin{gather*}
-J^2\bv =-J (E_{+} \bv \I + E_{-}\bv(-\I)) \\
= - \left( E_{+}^2 \bv \I^2  + E_+E_{-}\bv(-\I^2) + E_{-}E_{+} \bv(-\I^2) + E_{-}^2\bv(-\I)^2\right) = (E_{+} + E_{-})\bv.
\end{gather*}
Hence, $-J^2$ is the projection onto $V_{+}\oplus V_{-} = \ran(J)$ along $\ker J = V_0$ such that $J$ is actually an imaginary operator.

\end{proof}

As pointed out already several times, invariant subspaces of an operator are in the quaternionic setting not associated with spectral values but with entire spectral spheres. Hence quaternionic spectral measures associate subspaces of $V_{R}$ with sets of entire spectral spheres and not with arbitrary sets of spectral values. If we want to integrate a function $f$ that takes non-real values with respect to a spectral measure $E$, then we need some additional information. We need to know how to multiply the different values that $f$ takes on a spectral sphere onto the vectors associated with the different spectral values in this sphere. This information is given by a suitable imaginary operator. Similar to \cite{Viswanath:1971}, we hence introduce now the notion of a spectral system.
\begin{definition}\label{SpecSys}
A spectral  system on $V_R$ is a couple $(E,J)$ consisting of a spectral measure and an imaginary operator $J$ such that
\begin{enumerate}[(i)]
\item \label{SSys1} $E$ and $J$ commute, i.e. $E(\Delta)J = JE(\Delta)$ for all $\Delta\in\sBorel(\hh)$ and
\item \label{SSys2}  $E(\hh\setminus\rr) = -J^2$, that is $E(\rr)$ is the projection onto $\ker J$ along $\ran J$ and $E(\hh\setminus\rr)$ is the projection onto $\ran J$ along $\ker J$.
\end{enumerate}
\end{definition}

\begin{definition}
We denote by $\bsIntrin(\hh)$ the set of all intrinsic slice functions on $\hh$ that are measurable with respect to the usual Borel sets $\Borel(\hh)$ on $\hh$.
\end{definition}
\begin{lemma}\label{IntStruct}
A function $f:\hh\to\hh$ belongs to $\bsIntrin(\hh)$ if and only if it is of the form $f(s) = \alpha(s) + \I_{s}\beta(s)$ with $\alpha,\beta\in\bsMeas(\hh,\rr)$ and $\beta(s) = 0$ for $s\in\rr$. 
\end{lemma}
\begin{proof}
If $f(s)=\alpha(s) + \I_s\beta(s)$ with $\alpha,\beta\in\bsMeas(\hh,\rr)$ and $\beta(s) = 0$ for $s\in\rr$, then we can set $\alpha(s_0,s_1):= \alpha(s_0 + \I s_1)$ and $\beta(s_0,s_1) = \beta(s+\I s_1)$ and $\beta(s_0,-s_1) := - \beta(s_0 + \I s_1)$ with $\I \in\SS$ arbitrary. Since $\alpha(s)$ and $\beta(s)$ are $\sBorel(\hh)$-measurable, they are constant on each sphere $[s]$ and so this definition is independent of the chosen imaginary unit $\I$. Since $\beta(s) = 0$ for real $s$, $\beta(s_0,s_1)$  is moreover well defined for $s_1 = 0$.  We find that $f(s) = \alpha(s) + \I_{s}\beta(s) = \alpha(s_0,s_1) + \I_{s}\beta(s_0,s_1)$ with $\alpha(s_0,s_1)$ and $\beta(s_0,s_1)$ taking real values and satisfying \eqref{CCond} such that $f$ is actually an intrinsic slice function. Moreover, the functions $\alpha(s)$ and $\beta(s)$ and the function $\varphi(s):= \I_{s}$ if $s\notin\rr$ and $\varphi(s):= 0$ if $s\in\rr$ are $\Borel(\hh)$-$\Borel(\hh)$-measurable. As $\beta(s) = 0$  if $s\in\rr$, we have $f(s)  = \alpha(s) + \I_s\beta(s) = \alpha(s) + \varphi(s)\beta(s)$ and hence the function $f$ is $\Borel(\hh)$-$\Borel(\hh)$-measurable too.

If on the other hand $f\in\bsIntrin(\hh)$ with $f(s) = \alpha(s_0,s_1) + \I_s\beta(s_0,s_1)$, then also $\alpha(s):= \frac{1}{2}\left(f(s) + f\left(\overline{s}\right)\right)=\alpha(s_0,s_1)$ and $\beta(s) := \frac{1}{2}\varphi(s) \left(f\left(\overline{s}\right) - f(s)\right)=\beta(s_0,s_1)$ with $\varphi(s)$ as above are $\Borel(\hh)$-$\Borel(\hh)$-measurable. Moreover $\beta(s) = 0$ if $s_1 = 0$. Since $f$ is intrinsic, these functions take values in $\rr$ and hence they are $\Borel(\hh)$-$\Borel(\rr)$-measurable. They are moreover constant on each sphere $[s]$ such that the preimages $\alpha^{-1}(A)$ and $\beta^{-1}(A)$ of each set  $A\in\Borel(\rr)$ are axially symmetric Borel sets in $\hh$. Consequently, $\alpha$ and $\beta$ are $\sBorel(\hh)$-$\Borel(\rr)$-measurable. Finally, $|f|^2 = |\alpha|^2 + |\beta|^2$ such that $f$ is bounded if and only if $\alpha$ and $\beta$ are bounded.
 
\end{proof}
\begin{corollary}\label{SymMeas}
Any function $f\in\bsIntrin(\hh)$ is $\sBorel(\hh)$-$\sBorel(\hh)$-measurable.
\end{corollary}
\begin{proof}
Let $\Delta\in\sBorel(\hh)$. Its inverse image $f^{-1}(\Delta)$ is a Borel set in $\hh$ because $f$ is $\Borel(\hh)$-$\Borel(\hh)$-measurable. If $s\in f^{-1}(\Delta)$, then $f(s) = \alpha(s_0,s_1)+\I_{s}\beta(s_0,s_1)\in\Delta$. The axially symmetry of $\Delta$ implies then that for any $s_{\I}=s_0 + \I s_1\in[s]$ with $\I\in\SS$ also $f(s_{\I}) = \alpha(s_0,s_1) + \I_{s} \beta(s_0,s_1) \in \Delta$ and hence $s_{\I} \in f^{-1}(\Delta)$. Thus $s\in f^{-1}(\Delta)$ implies $[s]\subset f^{-1}(\Delta)$ and so $f^{-1}(\Delta)\in\sBorel(\hh)$.

\end{proof}
We observe that the \Cref{IntStruct} implies that the spectral integrals of the component functions $\alpha$ and $\beta$ of any  $f = \alpha + \I_{s}\beta\in\bsIntrin(\hh)$ are defined by  \Cref{SpecMeas}.
\begin{definition}
Let $(E,J)$ be a spectral system on $V_R$. For $f \in \bsIntrin(\hh)$ with $f(s) = \alpha(s) + \I_s\beta(s)$ we define the spectral integral of $f$ with respect to $(E,J)$ as
\begin{equation}\label{EJInt}
 \int_{\hh} f(s) \,dE_{J}(s) :=  \int_{\hh} \alpha(s)\,dE(s) + J \int_{\hh} \beta(s)\,dE(s).
 \end{equation}
\end{definition}
The estimate \eqref{NormEst1} generalizes to
\begin{equation}\label{NormEst2}
\left\| \int_{\hh} f(s)\,dE(s) \right\|\leq C_{E}\|\alpha\|_{\infty} + C_{E}\|J\| \|\beta\|_{\infty} \leq C_{E,J}\|f\|_{\infty}
\end{equation}
with  
\[
C_{E,J} := C_{E}(1+\|J|).
\]
As a consequence of \Cref{HomLem} and the fact that $J$ and $E$ commute, we immediately obtain the following result.
\begin{lemma}\label{GenIntProp}
Let $(E,J)$ be a spectral system on $V_{R}$. The mapping $f\mapsto \int_{\hh} f(s) \, dE_J(s)$ is a continuous homomorphism from $(\bsIntrin(\hh),\|.\|_{\infty})$ to $\boundOP(V_{R})$. Moreover, if $T\in\boundOP(V_{R})$ commutes with $E$ and $J$, then it commutes with $\int_{\hh} f(s)\,dE_J(s)$ for any $f\in\bsIntrin(\hh)$.
\end{lemma}
From \Cref{DualRep} we furthermore immediately obtain the following lemma, which is analogue to \Cref{FabLem}.
\begin{corollary}\label{DualRepAll}
Let $(E,J)$ be a spectral system on $V_R$ and let $f = \alpha + \I\beta\in\bsIntrin(\hh)$. For an $\bv\in V_R$ and any $\bv^{\mathbf{\ast}}\in V_{R}^{*}$, we have
\[
\left \langle \bv^{\mathbf{\ast}}, \left[ \int_{\hh} f(s)\,dE_J(s)\right]\bv\right\rangle  = \int_{\hh} \alpha(s)\,d \left \langle \bv^{\mathbf{\ast}},E(s)\bv\right\rangle  + \int_{\hh} \beta(s)\,d \left \langle \bv^{\mathbf{\ast}}, E(s) J\bv\right\rangle.
\]
\end{corollary}

Similar to the what happens for the $S$-functional calculus, there exists a deep relation between quaternionic and complex spectral integrals on $V_{R}$.
\begin{lemma}\label{BackCompat}
Let $(E,J)$ be a spectral system on $V_{R}$, let $\I\in\SS$, let $E_+$  be the projection of $V_R$ onto $V_{J,\I}^+$ along $V_{J,0}\oplus V_{J,\I}^{-}$ 
and let $E_{-}$ be the projection of $V_R$ onto $V_{J,\I}^-$ along $V_{J,0}\oplus V_{J,\I}^{+}$, cf. \Cref{JVSplitThm}. For $\Delta\in\Borel(\cc_{\I})$, we set
\begin{equation}\label{CSpecMeas}
 E_{\I}(\Delta) : = \begin{cases} E_{+}E([\Delta]) & \text{if $\Delta \subset\cc_{\I}^+$ }\\
E(\Delta) & \text{if $\Delta \subset\rr$ }\\
E_{-}E(\Delta) & \text{if $\Delta \subset\cc_{\I}^-$ }\\
E_{\I}(\Delta\cap \cc_{\I}^+) + E_{\I}(\Delta\cap\rr) + E_{\I}(\Delta\cap\cc_{\I}^{-}) & \text{otherwise,}
\end{cases}
\end{equation}
where $\cc_{\I}^{+}$ and $\cc_{\I}^{-}$ are the open upper and lower halfplane in $\cc_{\I}$.
Then $E_{\I}$ is a spectral measure on $V_{R,\I}$. For any $f\in\bsIntrin(\hh)$, we have with $f_{\I} := f|_{\cc_{\I}}$ that
\begin{equation}\label{H-Ci-SpecInt}
\int_{\hh} f(s)\,dE_J(s) = \int_{\cc_{\I}} f_{\I}(z)\,dE_{\I}(s).
\end{equation}
\end{lemma}
\begin{proof}
Recall that $E$ and $J$ commute. For $\bv_+\in V_{J,\I}^{+}$, we thus have $JE(\Delta) \bv_{+} = E(\Delta) J \bv_{+} = E(\Delta) \bv_{+} \I$ such that $E(\Delta)\bv_{+}\in V_{J,\I}^{+}$ and in turn $E_+ E(\Delta)\bv_{+} = E(\Delta)\bv_{+}$. Similarly, we see that $E(\Delta) \bv_{\sim} \in V_{J,0}\oplus V_{J,\i}^{-}$ for $\bv_{\sim}\in V_{J,0}\oplus V_{J,\I}^{-}$ such that $E_+ E(\Delta) \bv_{\sim} = \bO$. Hence, if we decompose $\bv\in V_{R}$ as $\bv = \bv_{+} + \bv_{\sim}$ with $\bv_{+}\in V_{J,\I}^{+}$ and $\bv_{\sim}\in  V_{J,0}\oplus V_{J,\I}^{-}$ according to \Cref{JVSplitThm}, then  $E_+ E(\Delta) \bv = E_+ E(\Delta)\bv_{+} + E_+E(\Delta) \bv_{\sim} = E(\Delta)\bv_{+}$ and $E(\Delta)E_{+}\bv = E(\Delta)\bv_+$ such that altogether $E(\Delta)E_+ \bv = E_+ E(\Delta) \bv$. Analogous arguments show that $E_{-}E(\Delta) = E(\Delta) E_{-}$ and hence $E_+$, $E_{-}$, and $E(\Delta)$, $\Delta\in\sBorel(\hh)$, commute mutually.

Let us now show that $E_{\I}$ is actually a $\cc_{\I}$-complex linear spectral measure on $V_{R,\I}$. For each $\Delta\in\Borel(\cc_{\I})$ set $\Delta_+:= \Delta\cap\cc_{\I}^+$, $\Delta_{-}:=\Delta\cap\cc_{\I}^{-}$ and $\Delta_{\rr}:=\Delta\cap\rr$ for neatness and recall that $[\,\cdot\,]$ denotes the axially symmetric hull of  a set. For any $\Delta,\sigma\in\sBorel(\hh)$, we have then
\begin{align}\label{ZZzz}
E([\Delta_{+}])E(\sigma_{\rr}) &= E(\Delta_{\rr}) E([\sigma_{+}]) = 0 & E([\Delta_{-}])E(\sigma_{\rr}) &= E(\Delta_{\rr}) E([\sigma_{-}])  = 0
\end{align}
because of \cref{SM3} in \Cref{SpecMeas}. Moreover, $E_{+}$ and $E_{-}$ as well as $E([\Delta+])$, $E([\Delta_{-}])$ and $E(\Delta_{\rr})$ are projections that commute mutually, as we just showed. Since in addition $E_{+}E_{-} =  E_{-}E_{+} = 0$, we have
\begin{equation}\label{ZZzz1}
\begin{split}
 E_{\I}(\Delta)^2 = &\left(E_{+}E([\Delta_+]) + E(\Delta_{\rr}) + E_{-}E([\Delta_{-}])\right)^2\\
 = &E_{+}^2E([\Delta_+])^2+ E_{+}E([\Delta_+])E(\Delta_{\rr}) + E_{+}E_{-}E([\Delta_+])E([\Delta_{-}])\\
&+ E_{+} E(\Delta_{\rr}) E([\Delta_+]) + E(\Delta_{\rr})^2 + E_{-}E(\Delta_{\rr})E([\Delta_{-}])\\
&+ E_{-} E_{+}E([\Delta_{-}])E([\Delta_+]) + E_{-}E([\Delta_{-}])E(\Delta_{\rr}) + E_{-}^2E([\Delta_{-}])^2\\
=& E_{+}E([\Delta_+]) + E(\Delta_{\rr}) + E_{-}E([\Delta_{-}]) = E_{\I}(\Delta).
 \end{split}
 \end{equation}
Hence, $E_{\I}(\Delta)$ is a projection that is moreover continuous as $ \| E_{\I}(\Delta)\| \leq K(1+ \|E_{+}\| + \|E_{-}\|)$, where $K>0$ is the constant in \Cref{SpecMeas}. Althogether, we find that $E$ has takes values that are uniformly bounded projections in $\boundOP(V_{R,\I})$.

We obviously have $E_{\I}(\emptyset) = 0$. Since $ E_{+} +E_{-} = E(\hh\setminus\rr)$ because of \cref{SSys2} in \Cref{SpecSys} also
\[ 
E_{\I}(\cc_{\I}) = E_{+}E([\cc_{\I}^+]) + E(\rr) + E_{-}E([\cc_{\I}^{-}]) = (E_{+} + E_{-})E(\hh\setminus \rr) + E(\rr) = E(\hh) = \id.
\]
Using the same properties of $E_{+}$, $E_{-}$ and $E(\Delta)$ as in \eqref{ZZzz1}, we find that for $\Delta,\sigma\in\Borel(\cc_{\I})$
\begin{align*}
 & E_{\I}(\Delta)E(\sigma) = \\
 = &\big(E_{+}E([\Delta_+]) + E(\Delta_{\rr}) + E_{-}E([\Delta_{-}])\big)\big(E_{+}E([\sigma_+]) + E(\sigma_{\rr}) + E_{-}E([\sigma_{-}])\big)\\
 = &E_{+}^2E([\Delta_+])E([\sigma_+])+ E_{+}E([\Delta_+])E(\sigma_{\rr}) + E_{+}E_{-}E([\Delta_+])E([\sigma_{-}])\\
&+ E_{+} E(\Delta_{\rr}) E([\sigma_+]) + E(\Delta_{\rr})E(\sigma_{\rr}) + E_{-}E(\Delta_{\rr})E([\sigma_{-}])\\
&+ E_{-} E_{+}E([\Delta_{-}])E([\sigma_+]) + E_{-}E([\Delta_{-}])E(\sigma_{\rr}) + E_{-}^2E([\Delta_{-}])E([\sigma_{-}])\\
=& E_{+}E([\Delta_+]\cap[\sigma_+]) + E(\Delta_{\rr}\cap\sigma_{\rr}) + E_{-}E([\Delta_{-}]\cap[\sigma_{-}]).
 \end{align*}
In general it not true that $[A] \cap [B] = [A\cap B]$ for $A,B\subset \cc_{\I}$. (Just think for instance about $A = \{\I\}$ and $B=\{-\I\}$ with $[A]\cap[B] = \SS\cap\SS = \SS$ and $[A\cap B] = [\emptyset] = \emptyset$.) For any axially symmetric set $C$ we have however
\[
C = \left[C \cap \cc_{\J}^{\geq}\right]\qquad \forall \J\in\SS
\]
If $A$ and $B$ belong to the same complex halfplane $\cc_{\J}^{\geq}$, then 
\begin{equation}\label{JJjj}
[A]\cap[B] = \left[ ([A]\cap [B])\cap\cc_{\J}^{\geq}\right] = \left[ \left([A]\cap\cc_{\J}^{\geq}\right)\cap \left([B]\cap \cc_{\J}^{\geq}\right)\right] = [A\cap B].
\end{equation}
Hence $ [\Delta_{+} ]\cap [\sigma_{+}] =  [(\Delta\cap\sigma)_+]$ and $[\Delta_{-}]\cap[\sigma_{-}] = [(\Delta\cap\sigma)_{-}]$ such that altogether
\[ E_{\I}(\Delta)E_{\I}(\sigma) = E_{+}E([(\Delta\cap\sigma)_{+}]) + E(\Delta_{\rr}\cap\sigma_{\rr}) + E_{-}E([(\Delta\cap\sigma)_{-}]) = E_{\I}(\Delta\cap\sigma).\]
Finally, we find  for $\bv \in V_{R,\I} = V_{R}$ and any countable family $(\Delta_{n})_{n\in\nn}$ of pairwise disjoint sets that
\begin{align*}
&E_{\I}\left(\bigcup_{n\in\nn} \Delta_n\right)\bv =\\
=& E_{+}E\left(\left[ \bigcup_{n\in\nn} \Delta_{n,+}\right]\right)\bv + E\left(\bigcup_{n\in\nn} \Delta_{n,\rr}\right)\bv +E_{-} E\left(\left[\bigcup_{n\in\nn} \Delta_{n,-}\right]\right)\bv\\
= &E_{+}E\left( \bigcup_{n\in\nn} \left[ \Delta_{n,+}\right]\right)\bv + E\left(\bigcup_{n\in\nn} \Delta_{n,\rr}\right)\bv +E_{-} E\left(\bigcup_{n\in\nn} \left[\Delta_{n,-}\right]\right)\bv.
\end{align*}
Since the sets $\Delta_{n,+}$, $n\in\nn$ are pairwise disjoint sets in the upper halfplane $\cc_{\I}^+$, also their axially symmetric hulls are because of \eqref{JJjj}.
Similarly, also the axially symmetric hulls of the sets $\Delta_{n,-}$, $n\in\nn$ are pairwise disjoint such that
\begin{align*}
&E_{\I}\left(\bigcup_{n\in\nn} \Delta_n\right)\bv =\\
 =& \sum_{n\in\nn} E_{+}E_{\I}\left( \left[ \Delta_{n,+}\right]\right)\bv + \sum_{n\in\nn}E\left( \Delta_{n,\rr}\right)\bv + \sum_{n\in\nn}E_{-} E\left( \left[\Delta_{n,-}\right]\right)\bv \\
 = & \sum_{n\in\nn} E_{\I}(\Delta_n)\bv.
\end{align*}
Altogether, we see that $E_{\I}$ is actually a $\cc_{\I}$-linear spectral measure on $V_{\rr,\I}$.

Now let us consider spectral integrals. We start with the simplest real-valued function possible: $f = \alpha \chi_{\Delta}$ with $a\in\rr$ and $\Delta\in\sBorel(\hh)$. As $f_{\I} = \alpha \chi_{\Delta\cap\cc_{\I}}$ and $E(\Delta) = E_{\I}(\Delta_{\I}\cap\cc_{\I})$, we have for such function 
\[
\int_{\hh} f(s) \, dE(s) =  \alpha E(\Delta) = \alpha E_{\I}(\Delta\cap\cc_{\I}) = \int_{\cc_{\I}} f_{\I}(z) \,dE(z).
\]
By linearity we find that \eqref{H-Ci-SpecInt} holds true for any simple function $f(s) = \sum_{\ell=1}^{n} \alpha_{k}\chi_{\Delta(s)}$ in $\bsMeas(\hh,\rr)$. Since these functions are dense in $\bsMeas(\hh,\rr)$, it even holds true for any function in $\bsMeas(\hh,\rr)$. Now consider the function $\varphi(s) = \I_{s}$ if $s\in\hh\setminus\rr$ and $\varphi(s) = 0$ if $s\in\rr$. Since $\varphi_{\I}(z) = \I \chi_{\cc_{\I}^+} + (-\I)\chi_{\cc_{\I}^{-}}$ and $E_{\I}(\cc_{\I}^{+}) = E_{+}$ and $E_{\I}^{-} = E_{-}$, the integral of $\varphi_{\I}$  with respect to $E_{\I}$ is 
\begin{gather*}
 \int_{\cc_{\I}}\varphi(z) \,dE_{\I}(z) \bv=  \left(\I E_{\I}(\cc_{\I}^{})\right)\bv  + \left((-\I)E_{\I}(\cc_{\I}^{-})\right)\bv\\
  = E_{+}\bv\I + E_{-}\bv (-\I)  = J\bv
 \end{gather*}
for all  $\bv\in V_{R,\I} = V_{R}$. If $f$ is now an arbitrary function in $\bsIntrin(\hh)$, then $f(s) = \alpha(s) + \varphi(s)\beta(s)$  with $\alpha,\beta \in\bsMeas(\hh,\rr)$ and $\beta(s) = 0$ if $s\in\rr$ by \Cref{IntStruct}. By what we have shown so far and the homomorphism properties of both quaternionic and the complex spectral integrals, we thus find
\begin{align*}
 &\int_{\hh} f(s)\,dE_J(s) =\\
 =& \int_{\hh}\alpha(s)\,dE(s) + J \int_{\hh}\beta(s)\,dE(s) \\
 =& \int_{\cc_{\I}} \alpha_{\I}(z) \,dE_{\I}(z) + \left(\int_{\cc_{\I}} \varphi_{\I}(z)\, dE_{\I}(z)\right)\left(\int_{\cc_{\I}} \beta_{\I}(z)\, dE_{\I}(z)\right)\\
 =& \int_{\cc_{\I}} \alpha_{\I}(z) + \varphi_{\I}(z)\beta_{\I}(z) \,dE_{\I}(z) = \int_{\cc_{\I}} f_{\I}(z) \,dE_{\I}(z).
\end{align*}

\end{proof}
Working on a quaternionic Hilbert space, one might consider only spectral measures whose values are orthogonal projections. If $J$ is an anti-selfadjoint partially unitary operator as it happens for instance in the spectral theorem for normal operators in \cite{Alpay:2016}, then $E_{\I}$ has values that are orthogonal projections.
\begin{corollary}
Let $\hil$ be a quaternionic  Hilbert space, let $(E,J)$ be a spectral system on $\hil$, let $\I\in\SS$ and let $E_{\I}$ be the spectral measure defined in \eqref{CSpecMeas}. If $E(\Delta)$ is for any $\Delta\in\sBorel(\hh)$ an orthogonal projection on $\hil$ and $J$ is an anti-selfadjoint partially unitary operator, then $E_{\I}(\Delta_{\I})$ is  for any $\Delta_{\I}\in\Borel(\cc_{\I})$ an orthogonal projection on $(\hil_{\I}, \langle\cdot,\cdot,\rangle_{\I}$ with $\langle\bu,\bv\rangle_{\I} = \{\langle\bu,\bv\rangle\}_{\I}$ as in \Cref{realHil}.
\end{corollary}
\begin{proof}
If $\bu,\bv\in \hil_{J,\I}^{+}$, then 
\[
\langle  \bu,\bv\rangle = \langle \bu, -J^2\bv \rangle = \langle J\bu,J\bv\rangle = \langle \bu \I, \bv\I\rangle = (-\I)\langle \bu,\bv\rangle\I
\]
such that $\I\langle \bu,\bv\rangle = \langle\bu,\bv \rangle \I$. Since a quaternion commutes with $\I\in\SS$ if and only if it belongs to $\cc_{\I}$, we have $\langle\bu,\bv\rangle \in\cc_{\I}$. Hence, if we choose $\J\in\SS$ with $\I\perp\J$, then $\langle\bu,\bv\J\rangle = \langle \bu,\bv\rangle {\J} \in\cc_{\I}\J$ such that in turn $\langle\bu ,\bv\J\rangle_{\I} = \{\langle\bu,\bv\rangle\}_{\I} =0$ for $\bu,\bv\in \hil_{J,\I}^{+}$. Since $\hil_{J,\I}^{-} = \{\bv\J: \bv\in \hil_{J,\I}^{+}\}$ by \Cref{JVSplitThm}, we find $\hil_{J,\I}^{-}\perp_{\I}\hil_{J,\I}^{+}$, where $\perp_{\I}$ denotes orthogonality in $\hil_{\I}$. Furthermore, we have for $\bu\in \hil_{0}=\ker J$ and $\bv\in \hil_{J,\I}^{+}$ that
\[\langle\bu,\bv\rangle = \langle \bu,J\bv\rangle (-\I) = \langle J\bu,\bv\rangle\I = \langle \bO,\bv\rangle \I = 0\]
and so $\langle\bu,\bv\rangle_{\I} = \{\langle\bu,\bv\rangle\}_{\I} = 0$ and in turn $\hil_{J,\I}^{+}\perp \hil_{0}$. Similarly, we see that also $\hil_{J,\I}^{-}\perp _{\I}\hil_{0}$. Hence, the direct sum decomposition $\hil_{\I}= \hil_{J,0}\oplus\hil_{J,\I}^{+}\oplus\hil_{J,\I}^{-}$ in \eqref{J0VSplit} is actually a decomposition into orthogonal subspaces of $\hil_{\I}$. The projection $E_+$ of $\hil$ onto $\hil_{J,\I}^{+}$ along $\hil_{J,0}\oplus\hil_{J,\I}^{-}$ and the projection $E_{-}$ of $\hil$ onto $\hil_{J,\I}^{-}$ along $\hil_{J,0}\oplus\hil_{J,\I}^{+}$ are hence orthogonal projections on $\hil_{\I}$.

Since the operator $E(\Delta)$ is for $\Delta\in\sBorel(\hh)$ an orthogonal projection on $\hil$, it is an orthongal  projection on $\hil_{\I}$. A projection is orthogonal if and only if it is self-adjoint. Since $E_{+}$, $E_{-}$ and $E$ commute mutually, we find for any $\Delta\in\Borel(\cc_{\I})$ and any $\bu,\bv\in\hil_{\I} = \hil$ that
\begin{align*}
&\langle \bu,E_{\I}(\Delta)\bv \rangle_{\I} \\
=& \langle \bu, E_{+}E([\Delta\cap\cc_{\I}^+])\bv\rangle_{\I} + \langle \bu, E(\Delta\cap\rr)\bv\rangle_{\I} + \langle \bu, E_{-}E([\Delta\cap\cc_{\I}^-])\bv\rangle_{\I}\\
=& \langle  E_{+}E([\Delta\cap\cc_{\I}^+])\bu,\bv\rangle_{\I} + \langle E(\Delta\cap\rr)\bu,\bv\rangle_{\I} + \langle E_{-}E([\Delta\cap\cc_{\I}^-])\bu,\bv\rangle_{\I} \\
=& \langle E_{\I}(\Delta)\bu,\bv\rangle_{\I}.
\end{align*}
Hence, $E_{\I}(\Delta)$ is an orthogonal projection on $\hil_{\I}$.

\end{proof}
We present two easy examples of spectral systems that illustrate the intuition behind the concept of a spectral system. 
\begin{example}\label{ExCompact}
We consider a compact normal operator  $T$ on a quaternionic Hilbert space $\hil$. The spectral theorem for compact normal operators in \cite{Ghiloni:2014a} implies that the $S$-spectrum consists of a (possibly finite) sequence $[s_n] = s_{n,0} + \SS s_{n,1}, n\in \Upsilon\subset\nn$ of spectral spheres that are (apart from possibly the sphere $[0]$) isolated in $\hh$. Moreover it implies the existence of an orthonormal basis of eigenvectors $(\bb_{\ell})_{\ell\in\varLambda}$ associated with eigenvalues $ s_{\ell} = s_{\ell,0} + \I_{s_{\ell}}s_{\ell,1}$ with $\I_{s_{\ell}} = 0$ if $s_{\ell}\in\rr$ such that
\begin{equation}\label{COP}
T\bv = \sum_{\ell\in\varLambda} \bb_{\ell}  s_{\ell} \langle\bb_{\ell},\bv\rangle.
\end{equation}
Each eigenvalue $ s_{\ell}$ obviously belongs to one spectral sphere, namely to $[ s_{n(\ell)}]$ with $s_{n(\ell),0} = s_{\ell,0}$ and $s_{n(\ell),1} = s_{\ell,1}$, and for $[s_{n}]\neq \{ 0 \}$ only finitely many eigenvalues belong to the spectral sphere $[s_n]$. We can hence rewrite \eqref{COP} as
\begin{equation*}
T\bv = \sum_{[s_n]\in\sigma_{S}(T)} \sum_{s_{\ell}\in[s_{n}]}  \bb_{\ell}  s_{\ell} \langle\bb_{\ell},\bv\rangle = \sum_{n\in\Upsilon} \sum_{n(\ell) =n}  \bb_{\ell}  s_{\ell} \langle\bb_{\ell},\bv\rangle
\end{equation*}
The spectral measure $E$ of $T$ is then given by 
\[
E (\Delta)\bv = \sum_{\substack{n\in\Upsilon\\ [s_n]\subset \Delta}} \sum_{n(\ell) = n} \bb_{\ell}\langle\bb_{\ell},\bv\rangle \qquad\forall\Delta \in \sBorel(\hh)
\]
If $f\in\bsMeas(\hh,\rr)$, then obviously
\begin{equation}\label{SSulk}
\int_{\hh}f(s)\,dE(s)\bv = \sum_{n\in\Upsilon} E([s_n])\bv f(s_n)  = \sum_{n\in\Upsilon} \sum_{ n(\ell) = n} \bb_{\ell}\langle\bb_{\ell},\bv\rangle f(s_n).
\end{equation}
In particular
\[
\int_{\hh}s_0\,dE(s)\bv  = \sum_{n\in\Upsilon}\sum_{n(\ell) = n}  \bb_{\ell}   \langle\bb_{\ell},\bv\rangle s_{\ell,0}
\]
and
\[
\int_{\hh}s_1\,dE(s)\bv =  \sum_{n\in\Upsilon}\sum_{n(\ell) = n}   \bb_{\ell}   \langle\bb_{\ell},\bv\rangle s_{\ell,1}.
\]

 If we  define
\[
\mathsf{J}_{0}\bv := \sum_{n\in\Upsilon} \sum_{n(\ell) =n}\bb_{\ell}\I_{s_{\ell}}\langle \bb_{\ell},\bv \rangle,
\]
then $\mathsf{J}_{0}$ is an anti-selfadjoint partially unitary operator and $(E,\mathsf{J}_{0})$ is a spectral system. One can check easily that $E$ and $\mathsf{J}_{0}$ commute and, as $\I_{s_{\ell}} = 0$ for $s_{\ell} \in\rr$ and $\I_{s_{\ell}}\in\SS$ with $\I_{s_{\ell}}^2 = -1$, otherwise one has
\[
- \mathsf{J}_{0}^2\bv = - \sum_{n\in\Upsilon} \sum_{n(\ell) =n}\bb_{\ell}\I_{s_{\ell}}^2\langle \bb_{\ell},\bv \rangle = \sum_{n\in\Upsilon: [s_n]\subset\hh\setminus\rr} \sum_{n(\ell) =n}\bb_{\ell}\langle \bb_{\ell},\bv \rangle = E(\hh\setminus\rr) \bv.
\]
In particular $\ker \mathsf{J}_{0} = \clos{\linspan{\hh}\{ \bb_{\ell}: s_{\ell}\in\rr\}} = E(\rr)$. Note moreover that $\mathsf{J}_{0}$ is completely determined by $T$. 

For any function $f = \alpha + \I\beta \in\bsIntrin(\hh)$, we have because of \eqref{SSulk} and as $\langle \bb_{\ell},\bb_{\kappa}\rangle = \delta_{\ell,\kappa}$ that
\begin{align*}
&\int_{\hh}f(s)\,dE_{\mathsf{J}_{0}}(s) \bv =  \int_{\hh}\alpha(s)\,dE(s) \bv + \mathsf{J}_{0} \int_{\hh}f(s)\,dE(s) \bv\\
=  & \sum_{n\in\Upsilon} \sum_{ n(\ell) = n} \bb_{\ell}\langle\bb_{\ell},\bv\rangle \alpha (s_{n,0},s_{n,1}) \\
&+  \sum_{m,n\in\Upsilon} \sum_{ \substack{n(\ell) =n \\  n(\kappa) = m}}\bb_{\ell}\I_{s_{\ell}}\langle \bb_{\ell}, \bb_{\kappa}\rangle\langle\bb_{\kappa},\bv\rangle \beta (s_{m,0},s_{m,1}) \\
=  & \sum_{n\in\Upsilon} \sum_{ n(\ell) = n} \bb_{\ell} \alpha (s_{\ell,0},s_{\ell,1}) \langle\bb_{\ell},\bv\rangle \\
&+  \sum_{n\in\Upsilon} \sum_{ n(\ell) =n}\bb_{\ell}\I_{s_{\ell}}\beta (s_{\ell,0},s_{\ell,1}) \langle\bb_{\ell},\bv\rangle \\
=  & \sum_{n\in\Upsilon} \sum_{ n(\ell) = n} \bb_{\ell} ( \alpha (s_{\ell,0},s_{\ell,1}) + \I_{s_{\ell}}\beta (s_{\ell,0},s_{\ell,1})) \langle\bb_{\ell},\bv\rangle 
\end{align*}
and so
\begin{align}\label{AseL}
\int_{\hh}f(s)\,dE_{\mathsf{J}_{0}}(s) \bv =  \sum_{n\in\Upsilon} \sum_{ n(\ell) = n} \bb_{\ell} f( s_{\ell}) \langle\bb_{\ell},\bv\rangle.
\end{align}
In particular
\begin{align*}
\int_{\hh}s\,dE_{\mathsf{J}_{0}}(s) = \sum_{n\in\Upsilon} \sum_{ n(\ell) = n} \bb_{\ell} s_{\ell} \langle\bb_{\ell},\bv\rangle = T\bv.
\end{align*}
We in particular have $T = A + \mathsf{J}_{0} B$ with $A= \int_{\hh} s_0\,dE(s)$ self-adjoint, $B = \int_{\hh} s_1\,dE(s)$ positive and $\mathsf{J}_{0}$ anti-selfadjoint and partially unitary as in \eqref{SplitTeasy}. Moreover, $E$ corresponds via \Cref{SpecMeasOthers} to the spectral measure obtained from \Cref{SpecThm1}.

We choose now $\I,\J\in\SS$ with $\I\perp\J$ and for each $\ell\in\varLambda$ with $s_{\ell}\notin\rr$ we choose $h_{\ell}\in\hh$ with $|h_{\ell}| = 1$ such that $h_{\ell}^{-1}\I_{s_{\ell}} h_{\ell} = \I$ and in turn 
\[
h_{\ell}^{-1}s_{\ell}h_{\ell} = s_{\ell,0} + h_{\ell}^{-1}\I_{s_{\ell}}h_{\ell}s_{1} = s_{\ell,0} + \I s_{\ell,1} =: s_{\ell,\I}.
\]
 In order to simplify the notation we also set $h_{\ell} = 1$ and $\I_{s_{\ell}} = 0$ if $s_{\ell}\in\rr$. Then $\tilde{\bb}_{\ell} := \bb_{\ell} h_{\ell}, \ell\in\varLambda$ is another orthonormal basis consisting of eigenvector of $T$ and as $h_{\ell}^{-1} = \overline{h_{\ell}}/|h_{\ell}|^2 = \overline{h_{\ell}}$ we have
 \begin{equation}\label{RuTzi}
\begin{split}
T\bv =& \sum_{n\in\Upsilon} \sum_{ n(\ell) = n}\bb_{\ell}(h_{\ell} h_{\ell}^{-1}) s_{\ell}(h_{\ell} h_{\ell}^{-1})\langle \bb_{\ell},\bv\rangle \\
= & \sum_{n\in\Upsilon} \sum_{ n(\ell) = n}(\bb_{\ell}h_{\ell}) (h_{\ell}^{-1} s_{\ell}h_{\ell}) \langle \bb_{\ell}h_{\ell},\bv\rangle = \sum_{n\in\Upsilon} \sum_{ n(\ell) = n}\tilde{\bb}_{\ell}s_{\ell,\I} \langle \tilde{\bb}_{\ell},\bv\rangle
\end{split}
\end{equation}
and similarly
\begin{align*}
J_{0}\bv = & \sum_{n\in\Upsilon} \sum_{ n(\ell) = n}\bb_{\ell}(h_{\ell} h_{\ell}^{-1}) \I_{\ell}(h_{\ell} h_{\ell}^{-1})\langle \bb_{\ell},\bv\rangle \\
=& \sum_{n\in\Upsilon} \sum_{ n(\ell) = n}(\bb_{\ell}h_{\ell}) (h_{\ell}^{-1} \I_{\ell}h_{\ell}) \langle \bb_{\ell}h_{\ell},\bv\rangle =  \sum_{n\in\Upsilon} \sum_{ n(\ell) = n}\tilde{\bb}_{\ell}\I \langle \tilde{\bb}_{\ell},\bv\rangle.
\end{align*}
Recall that $\I \lambda = \lambda\I$ for any $\lambda\in\cc_{\I}$ and $\I\J = - \J\I$. The splitting of $\hil$ obtained from \Cref{JVSplitThm} is therefore given by
\[
\hil_{\mathsf{J}_{0},0} = \ker J_{0} = \clos{\linspan{\hh}\{\tilde{\bb}_{\ell}: s_{\ell} \in\rr\}},\qquad \hil_{\mathsf{J}_{0},\I}^{+} := \clos{\linspan{\cc_{\I}}\{\tilde{\bb}_{\ell}:s_{\ell}\notin\rr\}}
\]
and
\[
 \hil_{\mathsf{J}_{0},\I}^{-}= \clos{\linspan{\cc_{\I}}\{\tilde{\bb}_{\ell}\J:s_{\ell}\notin\rr\}} = \hil_{\mathsf{J}_{0},\I}^{+}\J.
\]
If $\langle \bb_{\ell}, \bv\rangle  = a_{\ell} = a_{\ell,1} + a_{\ell,2}\J$ with $a_{\ell,1},a_{\ell,2}\in\cc_{\I}$, then \eqref{RuTzi} implies
\begin{equation}\label{TrafL}
\begin{split}
 &T\bv =  \sum_{n\in\Upsilon} \sum_{ n(\ell) = n}\tilde{\bb}_{\ell}s_{\ell,\I} a_{\ell} \\
 =& \sum_{\substack{n\in\Upsilon \\ [s_n]\subset\rr}} \sum_{ n(\ell) = n}\tilde{\bb}_{\ell} a_{\ell} s_{\ell} + \sum_{\substack{n\in\Upsilon \\ [s_n]\subset\hh\setminus \rr}} \sum_{ n(\ell) = n}\tilde{\bb}_{\ell} a_{\ell,1} s_{\ell,\I} + \sum_{\substack{n\in\Upsilon\\ [s_n]\subset\hh\setminus \rr}} \sum_{ n(\ell) = n}\tilde{\bb}_{\ell} a_{\ell,2}\J \overline{s_{\ell,\I}}.
 \end{split}
 \end{equation}
If $f\in\bsIntrin(\hh)$, then the representation \eqref{AseL} of $\int_{\hh} f(s)\,dE_{\mathsf{J}_{0}}(s)$ in the basis $\tilde{\bb}_{\ell},\ell \in\varLambda$ implies
\begin{align}
\notag &\int f(s)\,dE_{\mathsf{J}_{0}}(s)\bv =  \sum_{n\in\Upsilon} \sum_{ n(\ell) = n}\tilde{\bb}_{\ell}f(s_{\ell,\I}) a_{\ell} \\
\notag =& \sum_{\substack{n\in\Upsilon \\ [s_n]\subset\rr}} \sum_{ n(\ell) = n}\tilde{\bb}_{\ell} a_{\ell} f(s_{\ell}) \\
\notag &+ \sum_{\substack{n\in\Upsilon \\ [s_n]\subset\hh\setminus \rr}} \sum_{ n(\ell) = n}\tilde{\bb}_{\ell} a_{\ell,1} f(s_{\ell,\I}) + \sum_{\substack{n\in\Upsilon\\ [s_n]\subset\hh\setminus \rr}} \sum_{ n(\ell) = n}\tilde{\bb}_{\ell} a_{\ell,2}\J \overline{f(s_{\ell,\I})}\\
\notag=& \sum_{\substack{n\in\Upsilon \\ [s_n]\subset\rr}} \sum_{ n(\ell) = n}\tilde{\bb}_{\ell} a_{\ell} f(s_{\ell})\\
\label{CuKAMA} & + \sum_{\substack{n\in\Upsilon \\ [s_n]\subset\hh\setminus \rr}} \sum_{ n(\ell) = n}\tilde{\bb}_{\ell} a_{\ell,1} f(s_{\ell,\I}) + \sum_{\substack{n\in\Upsilon\\ [s_n]\subset\hh\setminus \rr}} \sum_{ n(\ell) = n}\tilde{\bb}_{\ell} a_{\ell,2}\J f(\overline{s_{\ell,\I}}).
  \end{align}
 as $f(s_{\ell})\in\rr$ for $s_{\ell}\in\rr$ and $\overline{f(s_{\ell,\I})} = f(\overline{s_{\ell,\I}})$ because $f$ is intrinsic. Note that the representation \eqref{TrafL} and \eqref{CuKAMA} show clearly that $f(T)$ is actually defined by letting $f$ act on the right eigenvalues of $T$. 
 \end{example}
\begin{example}\label{ExMulti}
Let us consider the space $L^2(\rr,\hh)$ of all quaternion-valued functions on $\rr$ that are square-integrable with respect to the Lebesgue measure $\lambda$. Endowed with the pointwise multiplication $(fa)(t) = f(t)a$ for $f\in L^2(\rr,\hh)$ and $a\in\hh$ and with the scalar product 
\begin{equation}\label{L2}
\langle g, f\rangle = \int_{\rr}\overline{g(t)}f(t)\,d\lambda(t)\qquad\forall f,g\in L^2(\rr,\hh),
\end{equation}
this space is a quaternionic Hilbert space. Let us now consider a bounded measurable function $\varphi: \rr \to \hh$ and the multiplication operator $(M_{\varphi} f)(s) := \varphi(s) f(s)$. This operator is normal with $(M_{\varphi})^{\ast} = M_{\overline{\varphi}}$ and its $S$-spectrum is the set $\overline{\varphi(\rr)}$. Indeed, writing $\varphi(t) = \varphi_{0}(t) + \I_{\varphi(t)} \varphi_{1}(t)$ with $\varphi_{0}(t)\in\rr$, $\varphi_{1}(t)>0$ and $\I_{\varphi(t)} \in\SS$ for $\varphi(t)\in\hh\setminus\rr$ and $\I_{\varphi(t)} = 0$ for $\varphi(t)\in\rr$, we find that
\begin{align*}
\Q_{s}(M_{\varphi})f(t) = &M_{\varphi}^2f(t) -2s_0 M_{\varphi}f(t) + |s|^2 f(t) \\
= &(\varphi^2(t) - 2s_0 \varphi(t) + |s|^2)f(t)\\
= & (\varphi(t) - s_{\I_{\varphi(t)}}) (\varphi(t) - \overline{s_{\I_{\varphi(t)}}}) f(t)
\end{align*}
with $s_{\I_{\varphi(t)}} = s_0 + \I_{\varphi(t)} s_1$ and hence 
\[
\Q_{s}(M_{\varphi})^{-1}f(t) =  (\varphi(t) - s_{\I_{\varphi(t)}})^{-1} (\varphi(t) - \overline{s_{\I_{\varphi(t)}}})^{-1} f(t)
\]
is a bounded operator if $s\notin \overline{\varphi(\rr)}$. If we define $E(\Delta)  = M_{\chi_{\varphi^{-1}(\Delta)}}$ for all $\Delta \in\sBorel(\hh)$ then we obtain a spectral measure on $\sBorel(\hh)$, namely
\[
E(\Delta)f(t) = \chi_{\varphi^{-1}(\Delta)}(t) f(t).
\]
If we set 
\[
\mathsf{J}_{0}:= M_{\I_{\varphi}}\quad \text{i.e. }\quad (\mathsf{J}_{0} f)(t) = \I_{\varphi(t)} f(t),
\]
then we find that $(E,\mathsf{J}_0)$ is a spectral system. Obviously $\mathsf{J}_{0}$ is anti-selfadjoint and partially unitary and hence an imaginary operator that commutes with $E$. Since $\I_{\varphi(t)} = 0$ if $\varphi(t)\in\rr$  and $\I_{\varphi(t)} \in\SS$ otherwise, we have moreover
\[
(- \mathsf{J}_{0}^2 f) (t) = -\I_{\varphi(t)}^{2} f(t) = \chi_{\varphi^{-1}(\hh\setminus\rr)}f(t) = (E(\hh\setminus\rr) f)(t).
\]

If $g\in\bsMeas(\hh,\rr)$, then let $g_{n}(s) = \sum_{\ell =1}^{N_{n}}a_{n,\ell} \chi_{\Delta_{n,\ell}}(s)\in\bsMeas(\hh,\rr)$ be a sequence of simple functions that converges uniformly to $g$. Then 
\begin{gather*}
\int_{\hh} g(s)\,dE(s) f(t) = \lim_{n\to\infty}\sum_{\ell =1}^{N_{n}}a_{n,\ell} E(\Delta_{n,\ell}) f(t) =  \lim_{n\to\infty}\sum_{\ell =1}^{N_{n}}a_{n,\ell} \chi_{\varphi^{-1}(\Delta)}(t) f(t) \\
=  \lim_{n\to\infty}\sum_{\ell =1}^{N_{n}}a_{n,\ell} \chi_{\Delta}(\varphi(t)) f(t)  =  \lim_{n\to\infty}(g_{n}\circ\varphi)(t)f(t)  = (g\circ \varphi)(t) f(t).
\end{gather*} 
Hence, if $g(s) = \alpha(s) + \I_{s}\beta(s)\in\bsIntrin(\hh)$, then
\begin{align*}
&\int_{\hh}g(s)\,dE_{\mathsf{J}_0}(s)f(t) = \int_{\hh}g(s)\,dE(s)f(t) \\
=& \int_{\hh}\alpha(s)\,dE(s)f(t) + \mathsf{J}_{0}\int_{\hh}\beta(s)\,dE(s)f(t)\\
=& \alpha(\varphi(t))f(t) + \I_{\varphi(t)}\beta(\varphi(t)) f(t)\\
=&  (\alpha(\varphi(t)) + \I_{\varphi(t)}\beta(\varphi(t)) f(t) = (g\circ \varphi)(t)f(t)
\end{align*}
and so
\[
\int_{\hh}g(s)\,dE_{\mathsf{J}_0}(s) = M_{g\circ\varphi}. 
\]
Choosing $g(s) = s$, we in particular find $T = A + \mathsf{J}_{0} B$ with $A= \int_{\hh} s_0\,dE(s)$ self-adjoint, $B = \int_{\hh} s_1\,dE(s)$ positive and $\mathsf{J}_{0}$ anti-selfadjoint and partially unitary as in \eqref{SplitTeasy}. $E$ corresponds via \Cref{SpecMeasOthers} to the spectral measure obtained from \Cref{SpecThm1}.
\end{example}

\subsection{On the different approaches to spectral integration}\label{DiffSpInt}
Our approach to spectral integration specified some ideas in \cite{Viswanath:1971}. We conclude this section with a comparison of this approach with the approaches in \cite{Alpay:2016} and \cite{Ghiloni:2017}, which were explained quickly in \Cref{SpectThmSect}. All three approaches are consistent, if things are interpreted correctly. Let us first consider a spectral measure $E$ over $\cc_{\I}^{\geq}$ for some $\I\in\SS$  in the sense of \Cref{SpecMeasHil}, the values of which are orthogonal projections on a quaternionic Hilbert space $\hil$. Let furthermore $\mathsf{J}$ be a unitary anti-selfadjoint operator on $\hil$ that commutes with $E$ and let us interpret the application of $\mathsf{J}$ as a multiplication with $\I$ from the left as in \cite{Alpay:2016}. By \Cref{SpecMeasOthers}, we obtain a quaternionic spectral measure on $\sBorel(\hh)$ if we set $\tilde{E}(\Delta) := E\big(\Delta\cap\cc_{\I}^{\geq}\big)$ for $\Delta\in\sBorel(\hh)$ and obviously we have 
\[
 \int_{\hh} f(s) \,d\tilde{E}(s) = \int_{\cc_{\I}}f_{\I}(z)\,dE(z) \qquad\forall f\in\bsMeas(\hh,\rr),
 \]
where $f_{\I} = f|_{\cc_{\I}^{\geq}}$. If we set $J := \mathsf{J}\tilde{E}(\hh\setminus\rr) = \mathsf{J}E(\cc_{\I}^{+})$, then $J$ is an imaginary operator and we find that $(\tilde{E},J)$ is a spectral system on $\hil$. Now let $f(s) = \alpha(s) +\I\beta(s)\in\bsIntrin(\hh)$ and let again $f_{\I} = f|_{\cc_{\I}^{\geq}}$, $\alpha_{\I} = \alpha|_{\cc_{\I}^{\geq}}$ and $\beta_{\I} = \beta|_{\cc_{\I}^{\geq}}$. Since $\beta(s) = 0$ if $s\in\rr$, we have $\beta(s) = \chi_{\hh\setminus\rr}(s)\beta(s)$ and in turn
\begin{equation}\label{Kuakka}
\begin{split}
\int_{\cc_{\I}^{\geq}} f_{\I}(z)\, dE(z) =& \int_{\cc_{\I}^{\geq}} \alpha_{\I}(z)\, dE(z) +\mathsf{J} \int_{\cc_{\I}^{\geq}} \beta_{\I}(z)\, dE(z) \\
=&\int_{\hh} \alpha(s)\, d\tilde{E}(s) +\mathsf{J} \int_{\hh}\chi_{\hh\setminus\rr}(s) \beta(s)\, d\tilde{E}(s)\\
=&\int_{\hh} \alpha(s)\, d\tilde{E}(s) +\mathsf{J} E(\hh\setminus\rr)\int_{\hh} \beta(s)\, d\tilde{E}(s)\\
=&\int_{\hh} \alpha(s)\, d\tilde{E}(s) + J \int_{\hh} \beta(s)\, d\tilde{E}(s) = \int_{\hh}f(s)\,d\tilde{E}_{J}(s).
\end{split}
\end{equation}
Hence, for any measurable intrinsic slice function $f$, the spectral integral of $f$ with respect to the spectral system $(\tilde{E},J)$ coincides with the spectral integral of $f|_{\cc_{\I}^{\geq}}$ with respect to $E$, where we interpret the application of $\mathsf{J}$ as a multiplication with $\I$ from the left. Since the mapping $f\mapsto f|_{\cc_{\I}^{\geq}}$ is a bijection between the set of all measurable intrinsic slice functions on $\hh$ and the set of all measurable $\cc_{\I}$-valued functions on $\cc_{\I}^{\geq}$ that map the real line into itself, both approaches are equivalent for real slice functions if we identify $\tilde{E}$ with $E$ and $f$ with $f_{\I}$. The same identifications show that the approach in \cite{Ghiloni:2017} is equivalent to our approach, as long as we only consider intrinsic slice functions. Indeed, if $\mathcal{E} = (E,\mathcal{L})$ is an iqPVM over $\cc_{\I}^{\geq}$ on $\hil$, then $\mathsf{J}\bv:= L_{\I}\bv = \I\bv$ is a unitary and anti-selfadjoint operator on $\hil$. As before, we can set $\tilde{E}(\Delta) = E(\Delta\cap\cc_{\I}^{\geq})$ and  $J := \mathsf{J}\tilde{E}(\hh\setminus\rr) = L_{\I}E(\cc_{\I}^+)$. We then find as in \eqref{Kuakka} that
\begin{equation}
\int_{\cc_{\I}^{\geq}} f_{\I}(z)\, d\mathcal{E}(z)  = \int_{\hh}f(s)\,d\tilde{E}_{J}(s)\qquad\forall f\in\bsIntrin(\hh).
\end{equation}
For intrinsic slice functions, all three approaches are hence consistent.

Let us continue our discussion of how our approach to spectral integration fits into the existing theory. We recall that any normal operator $T$ on $\hil$ can be decomposed as
\[
 T = A + \mathsf{J}_{0}B ,
\]
with the selfadjoint operator $A= \frac{1}{2}(T+T^*)$, the positive operator $B = \frac{1}{2}|T-T^*|$ and the anti-selfadjoint partially unitary operator $\mathsf{J}_{0}$ with $\ker J_{0} = \ker(T-T^*) $ and $\ran J_{0} = \ker (T-T^*)^{\perp}$. Let $\mathcal{E} = (E,\mathcal{L})$ be the iqPVM of $T$ obtained from \Cref{SpecThm2}. From \cite[Theorem~3.13]{Ghiloni:2017}, we know that $\left(\int_{\cc_{\I}^{\geq}}\varphi(z)\,d\mathcal{E}(z)\right)^{\ast} = \int_{\cc_{\I}^{\geq}}\overline{\varphi(z)}\,d\mathcal{E}(z)$ and $\ker \int_{\cc_{\I}^{\geq}}\varphi(z)\,d\mathcal{E}(z) = \ran E(\varphi^{-1}(0))$. Hence
\[
T - T^* = \int_{\cc_{\I}^{\geq}} z\,d\mathcal{E}(z) -  \int_{\cc_{\I}^{\geq}} \overline{z}\,d\mathcal{E}(z) = \int_{\cc_{\I}^{\geq}} 2\I z_1\, d\mathcal{E}(z) .
\]
As $z_1=0$ if and only if $z\in\rr$, we find that $\ker \mathsf{J}_0 = \ker (T - T^*) =  \ran E(\rr)$ and in turn $\ran \mathsf{J}_{0} = \ker (T - T^*) ^{\perp} = \ran E\big(\cc_{\I}^{\geq}\setminus\rr\big) = \ran E(\cc_{\I}^{+})$. 

If we identify $E$ with the spectral measure $\tilde{E}$ on $\sBorel(\hh)$ that is given by $\tilde{E}(\Delta) = E\big(\Delta\cap\cc_{\I}^{\geq}\big)$, then $J = L_{\I}E(\cc_{\I}^{+})$ is an imaginary operator such that $(\tilde{E},J)$ is a spectral system as we shoed above. The spectral integral of any measurable intrinsic slice function $f$ with respect to $(\tilde{E},J)$ coincides with the spectral integral of $f|_{\cc_{\I}^{\geq}}$ with respect to $\mathcal{E}$. Since $\ran E(\cc_{\I}^{+}) = \ker(T-T^*)^{\perp} = \ran \mathsf{J}_{0}$ and $L_{\I}\bv = \mathsf{J}_{0}$ for all $\bv\in\ker(T-T^*)^{\perp}$, we moreover find $J = \mathsf{J}_0$. Therefore $(\tilde{E},\mathsf{J}_{0})$ is the spectral system that is for integration of intrinsic slice functions equivalent to $\mathcal{E}$. We can hence rewrite \Cref{SpecThm1} and \Cref{SpecThm2} in the terminology of spectral systems as follows.
\begin{theorem}
Let $T = A + \mathsf{J}_0B\in\boundOP(\hil)$ be a normal operator. There exists a unique quaternionic spectral measure $E$ on $\sBorel(\hh)$ with $E(\hh\setminus\sigma_{S}(T)) = 0$, the values of which are orthogonal projections on $\hil$, such that $(E,\mathsf{J}_{0})$ is a spectral system and such that
\[
T = \int_{\hh}s\,dE_{\mathsf{J}_{0}}(s).
\]
\end{theorem}
We want to point out that the spectral system $(E,\mathsf{J}_{0})$ is completely determined by $T$---unlike the unitary anti-selfadjoint operator $\mathsf{J}$ that extends $\mathsf{J}_{0}$ used in \cite{Alpay:2016} and unlike the iqPVM used in \cite{Ghiloni:2017}. We also want to stress that the proof of the spectral theorem in \cite{Alpay:2016} translates directly into the language of spectral systems:  the spectral measure for $T$ is constructed using only real-valued functions, hence the extension  $\mathsf{J}$ of $\mathsf{J}_{0}$ is irrelevant in this proof. Indeed, in the article \cite{Alpay:2016}, only functions that are restrictions of intrinsic slice functions are integrated so that one can pass to the language of spectral systems by the identification described above without any problems.  
\begin{example}\label{ExCompact2}
In order to discuss the relations described above let us return to \Cref{ExCompact}, in which we considered normal compact operator on a quaternionic Hilbert space given by
\[
T\bv =  \sum_{n\in\Upsilon} \sum_{n(\ell) =n}  \bb_{\ell}  s_{\ell} \langle\bb_{\ell},\bv\rangle,
\]
the spectral system $(E,\mathsf{J}_{0})$ of which was\[
E (\Delta)\bv = \sum_{\substack{n\in\Upsilon \\ [s_n]\subset \Delta}} \sum_{n(\ell) = n} \bb_{\ell}\langle\bb_{\ell},\bv\rangle\qquad \text{and} \qquad \mathsf{J}_{0}\bv = \sum_{n\in\Upsilon} \sum_{n(\ell) =n}\bb_{\ell}\I_{s_{\ell}}\langle \bb_{\ell},\bv \rangle.
\]
The integral of a function $f\in\bsIntrin(\hh)$ with respect to $(E, \mathsf{J}_{0})$ is then given by \eqref{AseL} as
\begin{equation}\label{RASL}
\int_{\hh}f(s)\,dE_{\mathsf{J}_{0}}(s) \bv =  \sum_{n\in\Upsilon} \sum_{ n(\ell) = n} \bb_{\ell} f( s_{\ell}) \langle\bb_{\ell},\bv\rangle.
\end{equation}

Let $\I\in\SS$. If we set $\tilde{E}(\Delta) = E([\Delta])$ for any $\Delta\in\Borel\big(\cc_{\I}^{\geq}\big)$, then $\tilde{E}$ is a quaternionic spectral measure over $\cc_{\I}$. In \cite{Alpay:2015} the authors extend $\mathsf{J}_{0}$ to an anti-selfadjoint and unitary operator $\mathsf{J}$ that commutes with $T$ and interpret applying this operator as a multiplication with $\I$ from the left in order to define spectral integrals. One possibility to do this is to define $\iota(\ell) = \I_{s_{\ell}}$ if $s_{\ell}\not\in\rr$ and $\iota(\ell) \in\SS$ arbitrary if $s_{\ell}\in\rr$ and to set
\[
 \mathsf{J}\bv = \sum_{n\in\Upsilon} \sum_{n(\ell) =n}\bb_{\ell}\iota(\ell)\langle \bb_{\ell},\bv \rangle
\]
and $\mathcal{\I} \bv = \mathsf{J}\bv$

In \cite{Ghiloni:2017} the authors go even one step further and extend this multiplication with scalars from the left to a full left multiplication $\mathcal{L} = (L_{a})_{a\in\hh}$ that commutes with $E$ in order obtain an iqPVM $\mathcal{E} = (E,\mathcal{L})$. We can do this by choosing for each $\ell\in\varLambda$ an imaginary unit $\j(\ell)\in \SS$ with $\j(\ell)\perp\iota(\ell)$ and by defining
\[
 \mathsf{K}\bv = \sum_{n\in\Upsilon} \sum_{n(\ell) =n}\bb_{\ell}\j(\ell)\langle \bb_{\ell},\bv \rangle.
\]
If we choose $\J\in\SS$ and define for $a = a_0 + a_1 \I + a_2 \J + a_3\I\J \in \hh$
\begin{align*}
a\bv = &L_{a}\bv  := \bv a_0 +  \mathsf{\J}\bv a_1 +  \mathsf{K}\bv a_2 + \mathsf{JK}\bv a_3 \\
=& \sum_{n\in\Upsilon} \sum_{n(\ell) =n}\bb_{\ell}(a_0 + a_1\iota(\ell) + a_2\j(\ell) + a_3\iota(\ell)\j(\ell)\langle \bb_{\ell},\bv \rangle,
\end{align*}
then $\mathcal{L} = (L_{a})_{a\in\hh}$ is obviously a left multiplication that commutes with $E$ and hence $\mathcal{E} = (\tilde{E},\mathcal{L})$ is an iqPVM over $\cc_{\I}^{\geq}$.

Set $s_{n,\I} =[s_{n}]\cap\cc_{\I}$. For $f_{\I}:\cc_{\I}^{\geq}\to\hh$, the integral of $f_{\I}$ with respect to $\mathcal{E}$ is 
\begin{equation}\label{RumZL}
\begin{split}
&\int_{\cc_{\I}^{\geq}} f_{\I}(z) \,d\mathcal{E}(z) \\
=&  \sum_{n\in\Upsilon}  f_{\I}(s_{n,\I}) \tilde{E}(\{s_{n,\I}\})\bv =\sum_{n\in\Upsilon} f_{\I}(s_{n,\I})E([s_n]) \bv\\
=& \sum_{n\in\Upsilon}\big(f_0( s_{n,\I}) + f_{1}(s_{n,\I})\mathsf{J} + f_{2}(s_{n,\I})\mathsf{K} + f_{3}(s_{n,\I})\mathsf{J}\mathsf{K}\big) \sum_{ n(\ell) = n} \bb_{\ell}  \langle\bb_{\ell},\bv\rangle\\
=& \sum_{n\in\Upsilon} \sum_{ n(\ell) = n} \bb_{\ell} \big(f_0( s_{n,\I}) + f_{1}(s_{n,\I})\iota(\ell) + f_{2}(s_{n,\I})\j(\ell) + f_{3}(s_{n,\I}) \iota(\ell)\j(\ell)\big) \langle\bb_{\ell},\bv\rangle,
\end{split}
\end{equation}
where $f_{0},\ldots, f_{3}$ are the real-valued component functions such that $f_{\I}(z) = f_{0}(z) + f_{1}(z)\I + f_{2}(z)\J + f_{3}(z)\I\J$. If now $f_{\I}$ is the restriction of an intrinsic slice function $f(s) = \alpha(s) + \I_{s} \beta(s)$, then $f_{0} (s_{n(\ell),\I}) = \alpha(s_{\ell,\I}) = \alpha(s_{\ell})$ and $f_{1}(s_{n(\ell),\I}) = \beta(s_{\ell,\I}) = \beta(s_{\ell})$ and $f_{2}(z) = f_{3}(z) = 0$. As moreover $f_1(s_{n(\ell),\I}) =  \beta(s_{\ell}) = 0$ if $s_{\ell}\in\rr$ and $\iota(\ell) = \I_{s_{\ell}}$ if $s_{\ell}\notin\rr$, we find that actually \eqref{RumZL} equals \eqref{RASL} in this case. Note however that for any other function $f_{\I}$, the integral \eqref{RumZL} depends on the random choice of the functions $\iota(\ell)$ and $\j(\ell)$, which are not fully determined by $T$.

Let us now investigate the relation of \eqref{RumZL} with the right linear structure of $T$. Let us therefore change to the eigenbasis $\tilde{\bb}_{\ell}, \ell\in\Lambda$ with $T \tilde{\bb}_{\ell} = \tilde{\bb}_{\ell}s_{\ell,\I}$ defined in \Cref{ExCompact}. For convenience let us furthermore choose $\iota(\ell)$ and $\j(\ell)$ such that 
\[
 \mathsf{J}\bv = \sum_{n\in\Upsilon} \sum_{n(\ell) =n}\tilde{\bb}_{\ell}\I\langle \tilde{\bb}_{\ell},\bv \rangle\qquad\text{and}\qquad \mathsf{K}\bv = \sum_{n\in\Upsilon} \sum_{n(\ell) =n}\tilde{\bb}_{\ell}\J\langle \tilde{\bb}_{\ell},\bv \rangle.
\]
The left-multiplication $\mathcal{L}$ is hence exactly the left-multiplication induced by the basis $\tilde{\bb}_{\ell},\ell\in\varLambda$ and multiplication of $\bv$ with $a\in\hh$ from the left exactly corresponds to multiplying the coordinates $\langle\tilde{\bb}_{\ell},\bv\rangle$ with $a$ from the left, i.e. $a\bv = \sum_{n\in\Upsilon} \sum_{n(\ell) =n}\tilde{\bb}_{\ell}a\langle \tilde{\bb}_{\ell},\bv \rangle$. (Note however that unlike multiplication with scalars from the right, the multiplication with scalars from the left does only in this basis correspond to a multiplication of the coordinates. This relation is lost if we change the basis.)

Let us denote $\langle \tilde{\bb}_{\ell}, \bv \rangle = a_{\ell}$ with $a_{\ell} = a_{\ell,1} + a_{\ell,2}\J$ with $a_{\ell,1},a_{\ell,2}\in\cc_{\I}$ and let $f_{\I}:\cc_{\I}^{\geq} \to  \hh $. If we write   $f_{\I}(z) = f_{1}(z) + f_{2}(z)\J$, this time with $\cc_{\I}$-valued components $f_{1},f_{2}:\cc_{\I}^{\geq}\to\cc_{\I}$, then \eqref{RumZL} yields
\begin{equation}
\begin{split}
\int_{\cc_{\I}^{\geq}} f_{\I}(z) \,d\mathcal{E}(z) =& \sum_{n\in\Upsilon} \sum_{ n(\ell) = n}\tilde{\bb}_{\ell} \big(f_{1} (s_{n,\I}) +  f_{2}(s_{n,\I})\J \big) (a_1 + a_2\J)\\
=& \sum_{n\in\Upsilon} \sum_{ n(\ell) = n} \tilde{\bb}_{\ell} \big(a_1 f_{1} (s_{n,\I}) + \overline{a_1}f_{2}(s_{n,\I})\J \big) \\
&+ \sum_{n\in\Upsilon} \sum_{ n(\ell) = n} \tilde{\bb}_{\ell} \big(a_2\J \overline{f_{2} (s_{n,\I})} -  \overline{a_2}f_{2}(s_{n,\I}) \big).
\end{split}
\end{equation}
If we compare this with \eqref{TrafL}, then we find that $\int_{\cc_{\I}^{\geq}} f_{\I}(z) \,d\mathcal{E}(z)$ does only correspond to an application of $f_{\I}$ to the right eigenvalues of $T$ if $f_2 \equiv 0$ and $f_{1}$ can be extended to a function on all of $\cc_{\I}$ such that $f_{1}(\overline{s_{\ell,\I}}) = \overline{f_{1}(s_{\ell,\I})}$. This is however only the case if and only if $f_{\I} = f_{1}$ is the restriction of an intrinsic slice function to $\cc_{\I}^{\geq}$. 
\end{example}
As pointed out above, spectral integrals of intrinsic slice functions defined in the sense of \cite{Alpay:2016} or \cite{Ghiloni:2017} can be considered as spectral integrals with respect to a suitably chosen spectral system. The other two approaches---in particular the approach using iqPVMs in \cite{Ghiloni:2017}---allow however the integration of a larger class of functions. Nevertheless we have several arguments in favour of the approach using spectral systems:
\begin{enumerate}[(i)]
\item The approach using spectral systems generalizes to quaternionic right Banach spaces. Moreover, defining a measurable functional calculus with this approach does not require choosing any arbitrary additional structure such as a left multiplication. 

If we consider a normal operator $T = A+ \mathsf{J}_0 B$ on a quaternionic Hilbert space, then only its spectral system $\mathsf{J_0}$ is uniquely defined. The extension of $\mathsf{J}_{0}$ to a unitary anti-selfadjoint operator $\mathsf{J}$ that can be interpreted as a multiplication $L_{\I} = J$ with some $\I\in\SS$ from the left  and the multiplication $L_{\J}$ with some $\J\in\SS$ with $\J\perp\I$ that determine the left multiplication $\mathcal{L}$ in a iqPVM $\mathcal{E} = (E,\mathcal{L})$ associated with $T$ are not determined by $T$. Their construction in \cite{Ghiloni:2013} and \cite{Ghiloni:2017} is  based on the spectral theorems for quaternionic  selfadjoint operators and for complex linear normal operators.

 As we shall see in the next section, the spectral orientation $J$ of a spectral operator $T$ on a right Banach space can be constructed once the spectral measure $E$ associated with $T$ is known. Since the spectral theorem is not available on Banach spaces, it is however not clear how to extend $J$ to a fully imaginary operator or even further to something that generalizes an iqPVM and whether this is possible at all.

\item\label{Phen2} The natural domain of a right linear operator is a right Banach space. If a left multiplication is defined on the Banach space, then the operator's spectral properties should be independent of this left multiplication. Integration with respect to a spectral system $(E,J$) has a clear and intuitive interpretation with respect to the right linear structure of the space: the spectral measure $E$ associates (right) linear subspaces to spectral spheres and the spectral orientation determines how to multiply the spectral values in the corresponding spectral spheres (from the right) onto the vectors in these subspaces.

The role of the left multiplication in an iqPVM in terms of the right linear structure is less clear. Indeed, we doubt that there exists a similarly clear and intuitive interpretation in view of the fact that no relation between left and right eigenvalues has been discovered up to now.

\item\label{Phen3} Extending the class of admissible functions for spectral integration beyond the class of measurable intrinsic slice functions seems to add little value to the theory. As pointed out above, the proof of the spectral theorem in \cite{Alpay:2016} translates directly into the language of spectral systems and hence spectral systems offer a framework that is sufficient in order to prove the most powerful result of spectral theory. 

Even more, as we shall see in \Cref{Conclusions} (and already observed for a special case in \Cref{ExCompact2}), spectral integrals of functions that are not intrinsic slice functions cannot follow the basic intuition of spectral integration. In particular if we define a measurable functional calculus via spectral integration, then this functional calculus does only follow the fundamental intuition of a functional calculus, namely that $f(T)$ should be defined by action of $f$ on the spectral values of $T$,  if the underlying class of functions consists of intrinsic slice functions.

\end{enumerate}

\section{Bounded Quaternionic Spectral Operators}\label{SpecOpSect}
Let us now turn our attention to the study of quaternionic linear spectral operators. We will follow the complex linear theory in \cite{Dunford:1958}. A complex spectral operator is a bounded operator $A$ that has a spectral resolution, i.e. there exists a spectral measure $E$ defined on the Borel sets $\Borel(\cc)$ on $\cc$ such that $\sigma_{S}(A|_{\Delta})\subset\clos{\Delta}$ with $A_{\Delta} = A|_{\ran E(\Delta)}$ for all $\Delta\in\Borel(\cc)$. \Cref{SpecIntSect} showed  that spectral systems take over  the role of spectral measures in the quaternionic setting. If $E$ is a spectral measure that reduces an operator $T\in\boundOP(V_{R})$, then there will in general exist infinitely many imaginary operators $J$  such that $(E,J)$ is a spectral system. We thus have to find a criterion for identifying the one among them that fits the operator $T$ and can hence serve as its spectral orientation. A first and quite obvious requirement is that $T$ and $J$  commute. This is however not sufficient. Indeed, if $J$ and $T$ commute, then also $-J$ and $T$ commute. More general, any operator of the form $\tilde{J} := -E(\Delta)J + E(\hh\setminus \Delta) J$ with $\Delta\in\sBorel(\hh)$ is an imaginary operator such that $(E,\tilde{J})$ is a spectral system that commutes with $T$. 

We develop a second criterion by analogy with the finite-dimensional case. Let $T\in\boundOP(\hh^{n})$ be the operator on $\hh^n$ that is given by the diagonal matrix $T {= \diag(\lambda_1,\ldots,\lambda_{n})}$ and let us assume $\lambda_{\ell}\notin\rr$ for $\ell=1,\ldots,n$. We intuitively identify the operator $J =  \diag(\I_{\lambda_{1}},\ldots,\I_{\lambda_{n}})$  as the spectral orientation for $T$, cf. also \Cref{ExCompact}. Obviously $J$ commutes with $T$. Moreover, if $s_0\in\rr$ and $s_1>0$ are arbitrary, then 
the operator $ (s_0 \id - s_1 J)-T$ is invertible. Indeed, one has 
\[
(s_0 \id - s_1 J) -T = \diag( \overline{s_{\I_{\lambda_{1}}}} - \lambda_{1} ,\ldots, \overline{s_{\I_{\lambda_n}}}- \lambda_{n}),
\]
where $s_{\I_{\lambda_\ell}} = s_0 + \I_{\lambda_{\ell}} s_{1}$. Since $\overline{s_{\I_{\lambda_{\ell}}}} - \lambda_{\ell} = (s_{0}-\lambda_{\ell,0}) + \I_{\lambda_{\ell}}(-s_{1} - \lambda_{\ell,1})$ and  both $s_{1}> 0$ and $\lambda_{\ell,1}>0$ for all $\ell = 1,\ldots,n $, each of the diagonal elements has an inverse and so
\[
( (s_0 \id - s_1 J) - T)^{-1} = \diag\left( (\overline{s_{\I_{\lambda_{1}}}} - \lambda_{1})^{-1} ,\ldots, ( \overline{s_{\I_{\lambda_n}}} - \lambda_{n} )^{-1}\right).
\]
This invertibility is the criterion that uniquely identifies $J$.
\begin{definition}\label{SpecOP}
An operator $T\in\boundOP(V_R)$ is called a spectral operator if there exists a spectral decomposition for~$T$, i.e. a spectral system $(E,J)$ on $V_{R}$ such that the following three conditions hold:
\begin{enumerate}[(i)]
\item \label{SOi} The spectral system $(E,J)$ commutes with  $T$ , i.e. $E(\Delta)T = TE(\Delta)$ for all $\Delta\in\sBorel(\hh)$ and $TJ = J T$,
\item \label{SOii} If we set $T_{\Delta}:= T|_{V_{\Delta}}$ with $V_{\Delta} = E(\Delta)V_{R}$ for $\Delta\in\sBorel(\hh)$, then
\[ \sigma_S(T_{\Delta}) \subset \clos{\Delta}\quad \text{for all $\Delta\in\sBorel(\hh)$.}\]
\item \label{SOiii} For any $s_0,s_1\in\rr$ with $s_1> 0$, the operator $((s_0\id - s_1 J )-T)|_{V_1} $ has a bounded inverse on $V_1 := E(\hh\setminus\rr)V_R = \ran J$.
\end{enumerate}
The spectral measure $E$ is called a spectral resolution for $T$ and the imaginary operator $J$ is called a spectral orientation of $T$.
\end{definition}
A first easy result, which we shall use frequently, is that the restriction of a spectral operator to an invariant subspace $E(\Delta)V_{R}$ is again a spectral operator. 
\begin{corollary}\label{SubSpCor}
Let $T\in\boundOP(V_{R})$ be a spectral operator on $V_R$ and let $(E,J)$ be a spectral decomposition for $T$. For any $\Delta\in\sBorel(\hh)$, the operator $T_{\Delta} = T|_{V_{\Delta}}$ with $V_{\Delta} = \ran E(\Delta)$ is a spectral operator on $V_{\Delta}$. A spectral decomposition for $T_{\Delta}$ is $(E_{\Delta}, J_{\Delta})$ with $E_{\Delta}(\sigma) = E(\sigma)|_{V_{\Delta}}$ and $J_{\Delta} = J|_{V_{\Delta}}$. 
\end{corollary}
\begin{proof}
Since $E(\Delta)$ commutes with $E(\sigma)$ for $\sigma\in\sBorel(\hh)$ and $J$, the restrictions $E_{\Delta}(\sigma) = E(\sigma)|_{V_\Delta}$ and $J_{\Delta} = J|_{V_{\Delta}}$ are right linear operators on $V_{\Delta}$. It is immediate that $E_{\Delta}$ is a spectral measure on $V_{\Delta}$. Moreover
\[
\ker J_{\Delta} = \ker J \cap V_{\Delta} = \ran E(\rr)\cap V_{\Delta} = \ran E_{\Delta}(\rr)
\]
and 
\[
\ran J_{\Delta} = \ran J \cap V_{\Delta} = \ran E(\hh\setminus\rr) \cap V_{\Delta} =  \ran E_{\Delta}(\hh\setminus\rr).
\]
Since 
\[
-J_{\Delta}^2 = -J^2|_{V_{\Delta}} = E(\hh\setminus\rr)|_{V_{\Delta}} = E_{\Delta}(\hh\setminus\rr),
\]
the operator $-J_{\Delta}^2$ is the projection of $V_{\Delta}$ onto $\ran J_{\Delta}$ along $\ker J_{\Delta}$. Hence $J_{\Delta}$ is an imaginary operator on $V_{\Delta}$. Moreover $(E_{\Delta},J_{\Delta})$ is a spectral system. As
\[
E_{\Delta}(\sigma)T_{\Delta}E(\Delta) = E(\sigma)T E(\Delta) = TE(\sigma)E(\Delta) = T_{\Delta}E_{\Delta}(\sigma)E(\Delta)
\]
and similar
\[
J_{\Delta}T_{\Delta}E(\Delta) = JT E(\Delta) = TJE(\Delta) = T_{\Delta}J_{\Delta}E(\Delta),
\]
this spectral system commutes with $T_{\Delta}$.

If $\sigma\in \sBorel(\hh)$ and we set $V_{\Delta,\sigma} = \ran E_{\Delta}(\sigma)$, then
\[
V_{\Delta,\sigma}  = \ran E(\sigma)|_{V_{\Delta}} = \ran E(\sigma) E(\Delta) = \ran E(\sigma\cap\Delta) = V_{\Delta\cap\sigma}.
\]
Thus $T_{\Delta}|_{V_{\Delta,\sigma}} = T |_{V_{\sigma\cap\Delta}}$ and so $\sigma_{S}(T_{\Delta,\sigma}) = \sigma_{S}(T_{\Delta\cap\sigma}) \subset \Delta\cup\sigma \subset \sigma$. Hence $E_{\Delta}$ is a spectral resolution for $T_{\Delta}$. Finally, for $s_0,s_1\in\rr$ with $s_1>0$, the operator $s_0\id - s_1J - T$ leaves the subspace $V_{\Delta,1} := \ran E_{\Delta}(\hh\setminus\rr) =\ran E(\Delta\cap(\hh\setminus\rr))$ invariant because it commutes with $E$. Hence, the restriction of $(s_0\id - s_1J - T)|_{V_1}^{-1}$ to $V_{\Delta,1}\subset V_{1} = \ran E(\hh\setminus\rr)$ is abounded linear operator on $V_{\Delta,1}$. It obviously is the inverse of $(s_0\id - s_1J_{\Delta}-T_{\Delta})|_{V_{\Delta,1}}$. Therefore $(E_{\Delta},J_{\Delta})$ is actually a spectral decomposition for $T_{\Delta}$, which hence is in turn spectral operator.

\end{proof}

The remainder of this section considers the questions of uniqueness and existence of the spectral decomposition $(E,J)$ of $T$. We recall the $V_{R}$-valued right slice-hyperholomorphic function $\RRes_{s}(T;\bv):=\Q_{s}(T)^{-1}\bv\overline{s} - T\Q_{s}(T)^{-1}\bv$ on $\rho_S(T)$ for $T\in\closOP(V_R)$ and $\bv\in V_R$, which was defined in \Cref{CRes} in order to give a representation of the $S$-functional calculus for intrinsic functions in terms of the right multiplication on $V_{R}$. If $T$ is bounded, then $\Q_{s}(T)^{-1}$ and $T$ commute and we have
\[\RRes_{s}(T;\bv) := \Q_{s}(T)^{-1}(\bv\overline{s}-T\bv).\]

\begin{definition}\label{SHExt}
Let $T\in\boundOP(V_R)$ and let $\bv\in V_R$. A $V_R$-valued right slice-hyperholomorphic function $\bof$ defined on an axially symmetric open set $\dom(\bof)\subset \hh$ with $\rho_S(T)\subset \dom(\bof)$  is called a slice-hyperholomorphic extension of $\RRes_{s}(T;\bv)$ if 
\begin{equation}
\label{SRExt}(T^2 -2s_0T + |s|^2\id)\bof(s) = \bv \overline{s}-T\bv \qquad \forall s\in\dom(\bof).
\end{equation}
\end{definition}
Obviously such an extension satisfies
\[ \bof(s) = \RRes_s(T;\bv)\qquad \text{for }s\in\rho_S(T).\]

\begin{definition}\label{SpecVec}
Let $T\in\boundOP(V_R)$ and let $\bv\in V_R$. The function $\RRes_{s}(T;\bv)$ is said to have the single valued extension property if any two slice hyperholomorphic extensions $\bof$ and $\bg$ of $\RRes_{s}(T;\bv)$ satisfy $\bof(s) = \bg(s)$ for $s\in \dom(\bof)\cap\dom(\bg)$. In this case 
\[\rho_S(\bv) := \bigcup\{\dom(\bof): \bof\text{ is a slice hyperholomorphic extension of $\RRes_{s}(T;\bv)$}\}\]
is called the $S$-resolvent set of $\bv$ and $\sigma_S(\bv) = \hh\setminus\rho_S(\bv)$ is called the $S$-spectrum of $\bv$.
\end{definition}
Since it is the union of axially symmetric sets, $\rho_S(\bv)$ is axially symmetric. Moreover, there exists a unique maximal extension of $\RRes_{s}(T;\bv)$ to $\rho_{S}(\bv)$. We shall denote this extension by $\bv(s)$.

We shall soon see that the single valued extension property holds for $\RRes_s(T;\bv)$ for any $\bv\in V_{R}$ if $T$ is a spectral operator. This is however not true for an arbitrary operator $T\in\boundOP(V_R)$. A counterexample can be constructed analogue to \cite[p. 1932]{Dunford:1958}.

\begin{lemma}\label{Lemma1}
Let $T\in\boundOP(V_R)$ be a spectral operator and let $E$ be a spectral resolution for $T$. Let $s\in\hh$ and let $\Delta\subset\hh$ be a closed axially symmetric set such that $s\notin \Delta$.  If $\bv\in V_{R}$ satisfies ${ (T^2 - 2s_0 T + |s|^2\id)\bv = \bO}$ then
\[E(\Delta)\bv = \bO \qquad\text{ and }\qquad E([s])\bv = \bv.\]
\end{lemma}
\begin{proof}
Assume that $\bv\in V_{R}$ satisfies ${ (T^2 - 2s_0 T + |s|^2\id)\bv = \bO}$ and let $T_{\Delta}$ be the restriction of $T$ to the subspace $V_\Delta = E(\Delta)V$. As $s\notin\Delta$, we have $s\in\rho_S(T_\Delta)$ and so $\Q_{s}(T_{\Delta})$ is invertible. Since $\Q_{s}(T_{\Delta})^{-1} = \Q_{s}(T)^{-1}|_{V_{\Delta}}$, we have
\[\Q_s(T_\Delta)^{-1} (T^2 - 2s_0 T + |s|^2\id)E(\Delta) = E(\Delta), \]
from which we deduce
\begin{align*} 
E(\Delta)\bv &= \Q_s(T_\Delta)^{-1}(T^2 - 2s_0 T + |s|^2\id)E(\Delta)\bv \\
&= \Q_s(T_\Delta)^{-1}E(\Delta)(T^2 - 2s_0 T + |s|^2\id)\bv = \bO.
\end{align*}
Now define for $n\in\nn$ the closed axially symmetric set
\[\Delta_n = \left\{ p\in\hh : \dist (p, [s]) \geq \frac1n\right\}.\]
By the above $E(\Delta_n)\bv = \bO$ and in turn
\[ (\id-E([s])) \bv = \lim_{n\to\infty}E(\Delta_n)\bv = \bO\]
so that $\bv = E([s])\bv$.

\end{proof}

\begin{lemma}
If $T\in\boundOP(V_{R})$ is a spectral operator, then for any $\bv\in V_{R}$ the function $\RRes_{s}(T;\bv)$ has the single valued extension property.
\end{lemma}
\begin{proof}
Let $\bv\in V_{R}$ and let $\bof$ and $\bg$ be two slice hyperholomorphic extensions of $\RRes_{s}(T;\bv)$. We set $\bh(s) = \bof(s) - \bg(s)$ for $s\in \dom(\bh) = {\dom(\bof)\cap\dom(\bg)}$. 

If $s\in \dom(\bh)$ then there exists an axially symmetric neighborhood $U\subset\dom(\bh)$ of $s$ and for any $p\in U$ we have
\begin{align*} 
&(T^2 - 2p_0T + |p|^2\id) \bh(p)  \\
=&(T^2-2p_0T+|p|^2\id)\bof(p) - (T^2-2p_0T+|p|^2)\bg(p)\\
=& (\bv\overline{p}-T\bv) - (\bv\overline{p}-T\bv) = \bO.
\end{align*}
If $E$ is a spectral resolution of $T$, then we can conclude from the above and \Cref{Lemma1}  that $E([p])\bh(p) = \bh(p)$ for $p\in U$. We consider now a sequence $s_n\in U$ with $s_n\neq s$ for $n\in  U$ such that $s_n\to s$ as $n\to\infty$ and find 
\begin{align*}
\bO &= E([s])E([s_n])\bh(s_n) = E([s])\bh(s_n) \to E([s])\bh(s) = \bh(s).
\end{align*}
Hence, $\bof(s) = \bg(s)$ and $\RRes_{s}(T,\bv)$ has the single valued extension property.

\end{proof}
\begin{corollary}
If $T\in\boundOP(V_{R})$ is a spectral operator, then for any $\bv\in V_{R}$ the function $\RRes_{s}(T;\bv)$ has a unique maximal slice hyperholomorphic extension to $\rho_S(\bv)$. We denote this maximal slice hyperholomorphic extension of $\RRes_{s}(T;\bv)$ by $\bv(\cdot)$.
\end{corollary}

\begin{corollary}\label{EmptySpec}
Let $T\in\boundOP(V_R)$ be a spectral operator and let $\bv\in V_R$. Then $\sigma_S(\bv)=\emptyset$ if and only if~$\bv = \bO$.
\end{corollary}
\begin{proof}
If $\bv = \bO$ then $\bv(s) =\bO$ is the maximal slice hyperholomorphic extension of $\RRes_{s}(T;\bv)$. It is defined on all of $\hh$ and hence $\sigma_S(\bv) =\emptyset$.

Now assume that $\sigma_S(\bv) = \emptyset$ for some $\bv\in V_R$ such that the maximal slice hyperholomorphic extension  $\bv(\cdot)$ of $\RRes_{s}(T;\bv)$ is defined on all of $\hh$. For any $\bw^*\in V_R^{*}$, the function $s\to \langle \bw^* ,\bv(s)\rangle$ is an entire right slice-hyperholomorphic function. From the fact that $\RRes_{s}(T;\bv)$ equals the resolvent of $T$ as a bounded operator on $V_{R,\I_{s}}$, we deduce $\lim_{s\to\infty}\RRes_{s}(T;\bv) = \bO$ and then 
\[\lim_{s\to\infty} \langle \bw^*, \bv(s)\rangle = \lim_{s\to\infty}\langle \bw^*, \RRes_{s}(T;\bv)\rangle =  0.\]
Liouville's Theorem for slice hyperholomorphic functions in \cite{Colombo:2016b} hence  implies $\langle \bw^*, \bv(s)\rangle= 0$ for all $s\in\hh$. Since $\bw^*$ was arbitrary, we obtain $\bv(s) =\bO$ for all $s\in\hh$. 

Finally, we can choose $s\in\rho_S(T)$ such that the operator $\Q_{s}(T) = T^2 - 2s_0T + |s|^2\id$ is invertible and we find because of \eqref{SRExt} that
\begin{align*}
\bO =&  \bv(s)s - T\bv(s) =   \Q_{s}(T)^{-1}\Q_{s}(T)\bv(s)s - T\Q_{s}(T)^{-1}\Q_{s}(T)\bv(s) \\
=& \Q_s(T)^{-1}\left( \Q_{s}(T)\bv(s)s - T \Q_s(T)\bv(s)\right) = \\
=&\Q_{s}(T)^{-1}( (\bv\overline{s} - T\bv)s - T(\bv\overline{s} - T\bv)) \\
= &\Q_{s}(T)^{-1}\big(T^2\bv - T\bv 2s_0+ \bv|s|^2\big) = \Q_{s}(T)^{-1}\Q_{s}(T) \bv = \bv.  
\end{align*}

\end{proof}
\begin{theorem}\label{SubSpec}
Let $T\in\boundOP(V_{R})$ be a spectral operator and let $E$ be a spectral resolution for $T$. If $\Delta\in\sBorel(\hh)$ is closed, then 
\[E(\Delta)V_{R} = \{\bv\in V_R: \sigma_S(\bv)\subset \Delta\}.\]
\end{theorem}
\begin{proof}
Let $V_\Delta = E(\Delta)V_{R}$  and let $T_\Delta$ be the restriction of $T$ to $V_{\Delta}$. Since $\Delta$ is closed, \Cref{SpecOP} implies $\sigma_S(T_\Delta)\subset  \Delta$. Moreover $\Q_{s}(T_{\Delta}) = \Q_{s}(T)|_{V_{\Delta}}$ for $s\in\hh$. If  $ \bv \in V_{\Delta}$, then 
\[
\Q_{s}(T) \RRes_{s}(T;\bv) = \Q_{s}(T_{\Delta}) \Q_{s}(T_{\Delta})^{-1}(\bv\overline{s}- T_{\Delta}\bv) = \bv\overline{s}-T\bv
\]
for $ s\in\rho_S(T_{\Delta})$ and hence $\RRes_{s}(T_{\Delta};\bv)$ is a slice hyperholomorphic extension of $\RRes_{s}(T;\bv)$ to $\rho_S(T_{\Delta}) \supset \hh\setminus\Delta$. Thus $\sigma_S(\bv)\subset \Delta$. Since $\bv\in V_{R}$ was arbitrary, we find $E(\Delta)V_R\subset \{\bv\in V_R: \sigma_S(\bv)\subset \Delta\}$.

In order to show the converse relation, we assume that $\sigma_S(\bv)\subset\Delta$. We consider a closed subset $\sigma\in\sBorel(\hh)$ of the complement of $\Delta$ and set $T_{\sigma} = T|_{V_{\sigma}}$ with $V_{\sigma} = E(\sigma)V_R$. As above $\RRes_s(T_{\sigma}; E(\sigma)\bv)$ is  then a slice hyperholomorphic extension of $\RRes_s(T;E(\sigma)\bv)$ to $\hh\setminus\sigma$. If on the other hand $\bv(s)$ is the unique maximal slice hyperholomorphic extension of $\RRes_{s}(T;\bv)$, then
\begin{gather*}
 \Q_{s}(T) E(\sigma) \bv(s) = E(\sigma) \Q_{s}(T)\bv(s) \\
 = E(\sigma)(\bv\overline{s} - T\bv) = (E(\sigma)\bv)\overline{s} - T(E(\sigma)\bv )
\end{gather*}
for $s\in\hh\setminus\Delta$ and hence $E(\sigma)\bv(s)$ is a slice hyperholomorphic extension of $\RRes_{s}(T;E(\sigma)\bv)$ to $\hh\setminus\Delta$. Combining these two extensions, we find that $\RRes_s(T;E(\sigma)\bv)$ has a slice hyperholomorphic extension to all of $\hh$. Hence, $\sigma_S(E(\sigma)\bv) = \emptyset$ and in turn $E(\Delta)\bv = \bO$ by \Cref{EmptySpec}.

Let us now choose an increasing sequence of closed subsets $\sigma_n\in\sBorel(\hh)$ of $\hh\setminus\Delta$ such that $\bigcup_{n\in\nn}\sigma_n = \hh\setminus\Delta$. By the above arguments $E(\sigma_n)\bv = \bO$ for any $n\in\nn$. Hence
\[
 E(\hh\setminus\Delta)\bv = \lim_{n\to\infty} E(\Delta_n)\bv = \bO,
 \]
so that in turn $E(\Delta)\bv = \bv$. We thus obtain $E(\Delta)V_R\supset \{\bv\in V_R: \sigma_S(\bv)\subset \Delta\}$.

\end{proof}

\begin{corollary}\label{Support}
Let $T\in\boundOP(V_R)$ be a spectral operator and let $E$ be a spectral resolution of $T$. Then $E(\sigma_S(T)) = \id$.
\end{corollary}

\begin{corollary}
Let $T\in\boundOP(V_R)$ be a spectral operator and let $\Delta\in\sBorel(\hh)$ be closed. The set of all $\bv\in V_R$ with $\sigma_S(\bv)\subset \Delta$ is a closed right subspace of $V_{R}$. 
\end{corollary}

\begin{lemma}\label{Commut}
Let $T\in\boundOP(V_R)$ be a spectral operator. If $A\in\boundOP(V_{R})$ commutes with $T$, then $A$ commutes with every spectral resolution $E$ for $T$. Moreover, $\sigma_S(A\bv)\subset\sigma_S(\bv)$ for all $\bv\in V_R$.
\end{lemma}
\begin{proof}
For $\bv\in  V_R$ we have
\begin{align*}
& (T^2 - 2s_0T + |s|^2\id)A\bv(s) \\
=& A(T^2 - 2s_0T + |s|^2\id)\bv(s) \\
=& A(\bv\overline{s}-T\bv) = (A\bv)\overline{s}- T(A\bv).
\end{align*}
The function $A\bv(s)$ is therefore a slice hyperholomorphic extension of $\RRes_{s}(T;A\bv)$ to $\rho_S(\bv)$ and so $\sigma_S(A\bv)\subset \sigma_S(\bv)$. From \Cref{SubSpec} we deduce that
\[ A E(\Delta) V \subset E(\Delta) V\]
for any closed axially symmetric subset $\Delta$ of $\hh$. 

If $\sigma$ and $\Delta$ are two disjoint closed axially symmetric sets we therefore have
\[ E(\Delta) A E(\Delta) = A E(\Delta) \qquad \text{and}\qquad E(\Delta) A E (\sigma) = E(\Delta) E (\sigma) A E(\sigma) = 0.\]
If we choose again an increasing sequence of closed sets $\Delta_n\in\sBorel(\hh)$ with $\hh\setminus\Delta = \bigcup_{n\in\nn} \Delta_n$, we therefore have
\[ E(\Delta) A E(\hh\setminus \Delta)\bv = \lim_{n\to\infty}E(\Delta) A E(\Delta_n)\bv = \bO\qquad\forall \bv\in V_R \]
and hence
\begin{equation}\label{EA=AE}
E(\Delta) A = E(\Delta)A[E(\Delta) + E(\hh\setminus\Delta)] = E(\Delta)AE(\Delta) = AE(\Delta).
\end{equation}
Since $\Delta$ was an arbitrary closed set in $\sBorel(\hh)$ and since the sigma-algebra $\sBorel(\hh)$ is generated by sets of this type, we finally conclude that \eqref{EA=AE} holds true for any set $\sigma\in\sBorel(\hh)$.

\end{proof}

\begin{lemma}\label{EUnique}
The spectral resolution $E$ of a spectral operator $T\in\boundOP(V_R)$ is uniquely determined.
\end{lemma}
\begin{proof}
Let $E$ and $\tilde{E}$ be two spectral resolutions of  $T$. For any closed set $\Delta\in\sBorel(\hh)$, \Cref{SubSpec} implies
\[ \tilde{E}(\Delta)E(\Delta) = E(\Delta)\qquad\text{and}\qquad E(\Delta) \tilde{E}(\Delta)= \tilde{E}(\Delta)\]
and we deduce from \Cref{Commut} that $E(\Delta) = \tilde{E}(\Delta)$. Since the sigma algebra $\sBorel(\hh)$ is generated by the closed sets in $\sBorel(\hh)$, we obtain $E=\tilde{E}$ and hence the spectral resolution of $T$ is uniquely determined.

\end{proof}

Before we consider the uniqueness of the spectral orientation, we observe that for certain operators, the existence of a spectral resolution already implies the existence of a spectral orientation and is hence sufficient for them to be a spectral operator.

\begin{proposition}\label{SSOr}
Let $T\in\boundOP(V_R)$ and assume that there exists a spectral resolution $E$ for $T$. If $\sigma_{S}(T) \cap \rr = \emptyset$, then there exists an imaginary operator $J\in\boundOP(V_R)$ that is a spectral orientation for $T$ such that $T$ is a spectral operator with spectral resolution $(E,J)$. Moreover, this spectral orientation is unique.
\end{proposition}
\begin{proof}
Since $\sigma_S(T)$ is closed, the assumption $\sigma_{S}(T)\cap\rr = \emptyset$ implies that $\dist(\sigma_S(T),\rr)>0$. We choose $\I\in\SS$ and consider $T$ as a complex linear operator on $V_{R,\I}$. Because of \Cref{ABC}, the spectrum of $T$ as a $\cc_{\I}$-linear operator on $V_{R,\I}$ is $\sigma_{\cc_{\I}}(T) = \sigma_{S}(T)\cap\cc_{\I}$. As $\dist(\sigma_{S}(T),\rr)>0$, the sets 
\[
\sigma_{+} = \sigma_{\cc_{\I}}(T)\cap\cc_{\I}^{+}\qquad \text{and}\qquad \sigma_{-} = \sigma_{\cc_{\I}}(T)\cap\cc_{\I}^{-}
\]
are open and closed subsets of $\sigma_{\cc_{\I}}(T)$ such that $\sigma_{+} \cup \sigma_{-} = \sigma_{\cc_{\I}}(T)$.  Via the Riesz-Dunford functional calculus we can hence associate spectral projections $E_{+}$ and $E_{-}$ onto closed invariant $\cc_{\I}$-linear subspaces of $V_{R,\I}$ to $\sigma_{+}$ and $\sigma_{-}$. The resolvent of $T$ as a $\cc_{\I}$-linear operator on $V_{R,i}$ at $z\in\rho_{\cc_{\I}}(T)$ is $R_{z}(T)\bv := \Q_{z}(T)^{-1}(\bv\overline{z}-T\bv)$, and hence these projections are given by 
\begin{equation}\label{KuJaKu}
 \begin{split}
 &E_{+}\bv:= \int_{\Gamma_{+}} \Q_{z}(T)^{-1}(\bv\overline{z} - T \bv) \, dz\frac{1}{2\pi\I}\\
 & E_{-}\bv:= \int_{\Gamma_{-}} \Q_{z}(T)^{-1}(\bv\overline{z} - T \bv) \, dz\frac{1}{2\pi\I},
 \end{split}
 \end{equation}
where $\Gamma_+$ is a positively oriented Jordan curve that surrounds $\sigma_{+}$ in $\cc_{\I}^{+}$ and $\Gamma_{-}$ is a positively oriented Jordan curve that surrounds $\sigma_{-}$ in $\cc_{\I}^{-}$.  We set 
\[
 J\bv : = E_{-} \bv (-\I) + E_{+}\bv \I.
\]
From \Cref{JVSplitThm}  we deduce that $J$ is an imaginary operator on $V_{R}$ if $\Psi:\bv \mapsto \bv\J$  is a bijection between $V_{+}:= E_{+}V_R$ and $V_{-} := E_{-}V_{R}$ for $\J\in\SS$ with $\J\perp\I$. This is indeed the case: due to the symmetry of $\sigma_{\cc_{\I}}(T) = \sigma_{S}(T)\cap\cc_{\I}$ with respect to the real axis, we find $\sigma_{+}= \overline{\sigma_{-}}$  so that we can choose $\Gamma_{-}(t) = \overline{\Gamma_{+}(1- t)}$ for $t\in[0,1]$ in \eqref{KuJaKu}. Because of the relation \eqref{RLambdaCon} established in \Cref{ABC}, the resolvent $R_{z}(T)$ of $T$ as an operator on $V_{R,\I}$ satisfies  $R_{\overline{z}}(T)\bv = - \left[R_{z}(T)(\bv \J)\right]\J$ and so
\begin{gather*}
E_{-}\bv = \int_{\Gamma_{-}} R_{z}(T)\bv \, dz\frac{1}{2\pi\I} = -\int_{\Gamma_{+}} R_{\overline{z}}(T)\bv \,d\overline{z}\frac{1}{2\pi\I} \\
= \int_{\Gamma_{+}}\left[ R_{z}(T) (\bv\J)\right]\J \,d\overline{z}\frac{1}{2\pi\I} =  \int_{\Gamma_{+}}\left[ R_{z}(T) (\bv\J)\right] \,dz\frac{1}{2\pi\I} (-\J) = \left[E_{+}(\bv\J)\right](-\J).
\end{gather*}
Hence, we have
\begin{equation}\label{CUmuLA}
 \left(E_{-}\bv\right)\J =  E_{+}(\bv\J)\qquad\forall \bv\in V_{R}.
\end{equation}
If $\bv\in V_{-}$, then $\bv\J = (E_{-}\bv)\J = E_{+}(\bv \J)$ and so $\bv\J \in V_{+}$. Replacing $\bv$ by $\bv\J$ in \eqref{CUmuLA}, we find that also $ \left(E_{-}\bv\J\right)\J =  -E_{+}(\bv)$ and in turn $E_{-}(\bv \J) = E_{+}(\bv )\J$. For $\bv \in V_{+}$ we thus find $\bv\J = E_{+}(\bv)\J = E_{-}(\bv \J)$ and so $\bv\J\in V_{-}$. Hence, $\Psi$ maps $V_{+}$ to $V_{-}$ and $V_{-}$ to $V_{+}$ and as $\Psi^{-1} = - \Psi$ it is even bijective. We conclude that $J$ is actually an imaginary operator.

Let us now show that \ref{SOi} in \Cref{SpecOP} holds true.
For any $\Delta\in\sBorel(\hh)$ the operator $\Q_{z}(T)^{-1}$ commutes with  $E(\Delta)$. Hence
\begin{equation}
\begin{split}
 E(\Delta) E_{+}\bv =& \int_{\Gamma_{+}}E(\Delta) \Q_{z}(T)^{-1}(\bv\overline{z} - T \bv) \, dz\frac{1}{2\pi\I}  \\
 =& \int_{\Gamma_{+}} \Q_{z}(T)^{-1}(E(\Delta)\bv\overline{z} - TE(\Delta) \bv) \, dz\frac{1}{2\pi\I} = E_+ E(\Delta)\bv
\end{split}
\end{equation}
for any $\bv \in V_{R,\I} = V_{R}$ and so $E_{+}E(\Delta)=E(\Delta)E_{+}$.  Similarly, one can show that also $E(\Delta)E_{-} = E_{-}E(\Delta)$. By construction, the operator $J$ hence commutes with $T$ and with $E(\Delta)$ for any $\Delta\in\sBorel(\hh)$ as
\[
TJ\bv =  TE_{-} \bv (-\I) + TE_{+}\bv \I = E_{-}T\bv(-\I) + E_{+}T\bv\I = JT\bv
\]
and 
\begin{gather*}
E(\Delta)J\bv =  E(\Delta)E_{-} \bv (-\I) + E(\Delta)E_{+}\bv \I\\
 = E_{-}E(\Delta)\bv(-\I) + E_{+}E(\Delta)\bv\I = JE(\Delta)\bv.
\end{gather*}
Moreover, as $\sigma_{S}(T) \cap\rr = \emptyset$, \Cref{Support} implies $\ran E(\rr) = \{\bO\} = \ker J$ and  $\ran E(\hh\setminus\rr) = V_{R} = \ran J$. Hence, $(E,J)$ is actually a spectral system that moreover commutes with $T$.

Let us now show condition \ref{SOiii} of \Cref{SpecOP}. If $s_0,s_1\in\rr$ with $s_1 > 0$, then set $s_{\I} := s_{0} + \I s_1$. As $E_+ + E_- = \id$, we then have
\begin{align*}
&((s_0\id - s_1 J) - T )\bv = \\
=&  (E_{+} + E_{-})\bv s_0 - (E_{+}\bv)\I s_1 - (E_{-}\bv)(-\I)s_1 - T(E_{+} + E_{-})\bv \\
=& (E_+\bv) (s_0  -  s_1\I)  - T (E_{+} \bv)  + (E_{-}\bv)) (s_0 + s_1\I) - T (E_{-} \bv) \\
=& (E_{+}\bv) \overline{s_{\I}} - T(E_{+}\bv) +  (E_{-}\bv)s_{\I}- T(E_{-}\bv) \\
=& (\overline{s_{\I}}\id_{V_{R,\I} } - T) E_{+} \bv + ( s_{\I} \id_{V_{R,\I}}   - T)E_{-}\bv.
\end{align*}
Since $E_{+}$ and $E_{-}$ are the Riesz-projectors associated to $\sigma_{+}$ and $\sigma_{-}$, the spectrum $\sigma(T_{+})$ of  $T_+ := T|_{V_{+}}$ is $\sigma_{+} \subset \cc_{\I}^{+}$  and the spectrum $\sigma(T_{-})$ of $T_- := T|_{ V_{-}}$ is $\sigma_{-} \subset \cc_{\I}^{-}$. As $s_{\I}$ has positive imaginary part, we find $\overline{s_{\I}}\in \cc_{\I}^{-}\subset\rho(T_+)$ and $s_{\I}\in\cc_{\I}^{+}\subset\rho(T_{-})$ such that  $R_{\overline{s}_{\I}}(T_+):= \left(\overline{s_{\I}}\id_{V_{+} } - T_{+}\right)^{-1}\in\boundOP(V_{+})$ and $R_{s_{\I}}(T)^{-1} := \left(s_{\I}\id_{V_{-} } - T_{-}\right)^{-1}\in\boundOP(V_{-})$ exist.   As $E_{+}|_{V_{+}} = \id_{V_{+}}$ and $E_{-}|_{V_{+}} = 0$, they satisfy the relations
\begin{equation}\label{KKkk}
E_{+} R_{\overline{s_{\I}}}(T_{+})E_{+} = R_{\overline{s_{\I}}}(T_{+})E_{+} \qquad \text{and}\qquad E_{-} R_{\overline{s_{\I}}}(T_{+})E_{+} = 0
\end{equation}
and similarly also
 \begin{equation}\label{LLll}
 E_{-} R_{s_{\I}}(T_{-})E_{-} = R_{s_{\I}}(T_{-})E_{-} \qquad\text{and}\qquad E_{+} R_{s_{\I}}(T_{-})E_{-} = 0.
 \end{equation}
 Setting $R(s_0,s_1) : = R_{\overline{s_{\I}}}(T_{+})E_{+} + R_{s_{\I}}(T_{-})E_{-}$, we obtain a bounded $\cc_{\I}$-linear operator that is defined on the entire space $V_{R,\I} = V_{R}$. Because $E_{+}$ and $E_{-}$ commute with $T$ and satisfy $E_{+}E_{-} = E_{-}E_{+} =0$ and because \eqref{KKkk} and \eqref{LLll} hold true, we find for any $\bv\in V_{R}$
\begin{align*}
&R(s_0,s_1)((s_0\id - s_1 J) - T )\bv   \\
= &\left[R_{\overline{s_{\I}}}(T_{+})E_{+} + R_{s_{\I}}(T_{-})E_{-}\right]\left[(\overline{s_{\I}}\id_{V_{R,\I} } - T) E_{+} \bv + ( s_{\I} \id_{V_{R,\I}}   - T)E_{-}\bv \right]\\
=&R_{\overline{s_{\I}}}(T_{+})(\overline{s_{\I}}\id_{V_{R,\I} } - T_{+}) E_{+} \bv + R_{s_{\I}}(T_{-})E_{-}( s_{\I} \id_{V_{R,\I}}   - T_{-})E_{-}\bv \\
=& E_{+}\bv + E_{-}\bv = \bv
\end{align*}
and
\begin{align*}
&((s_0\id - s_1 J) - T )R(s_0,s_1)\bv    \\
=& \left[(\overline{s_{\I}}\id_{V_{R,\I} } - T) E_{+}  + ( s_{\I} \id_{V_{R,\I}}   - T)E_{-} \right] \left[R_{\overline{s_{\I}}}(T_{+})E_{+} + R_{s_{\I}}(T_{-})E_{-}\right]\bv\\
= &(\overline{s_{\I}}\id_{V_{+} } - T_{+}) R_{\overline{s_{\I}}}(T_{+}) E_{+} \bv + ( s_{\I} \id_{V_{-}}   - T_{-}) R_{s_{\I}}(T_{-}) E_{-}\bv \\
= &E_{+}\bv + E_{-}\bv = \bv.
\end{align*}
Hence, $R(s_0,s_1)\in\boundOP(V_{R,\I})$ is the $\cc_{\I}$-linear bounded inverse of $(s_0\id - s_1 J) - T$. Since $(s_0\id - s_1 J) - T$ is quaternionic right linear, its inverse is quaternionic right linear too so that even $((s_0\id - s_1 J) - T)^{-1}\in\boundOP(V_R)$.  $J$ is therefore actually a spectral orientation for $T$ and $T$ is in turn a spectral operator with spectral decomposition $(E,J)$.

Finally, we show the uniqueness of the spectral orientation $J$. Assume that $\widetilde{J}$ is an arbitrary spectral orientation for $T$. We show that $\widetilde{V_{+}} := V_{\widetilde{J},\I}^{+}$ equals $V_{+} = V_{J,\I}^{+}$. \Cref{JVSplitThm} imlies then $J = \widetilde{J}$ because $\ker J = \ker \widetilde{J} = \ran E(\rr) = \{\bO\}$ and $V_{J, \I}^{-} = V_{+}\J = \widetilde{V_{+}}\J = V_{\widetilde{J},\I}^{-}$. 

Since $\widetilde{J}$ commutes with $T$, we have $\widetilde{J}E_{+} = E_{+}\widetilde{J}$ as
\begin{equation}\label{RompAL}
\begin{split}
\widetilde{J} E_{+}\bv =& \int_{\Gamma_{+}}\widetilde{J} \Q_{z}(T)^{-1}(\bv\overline{z} - T \bv) \, dz\frac{1}{2\pi\I}  \\
 =& \int_{\Gamma_{+}} \Q_{z}(T)^{-1}(\widetilde{J}\bv\overline{z} - T\widetilde{J} \bv) \, dz\frac{1}{2\pi\I} = E_+ \widetilde{J}\bv.
\end{split}
\end{equation}
The projection $E_{+}$ therefore leaves $\widetilde{V_{+}}$ invariant because $\widetilde{J} (E_{+}\bv) = E_{+}\big(\widetilde{J}\bv\big) =  (E_{+}\bv)\I \in \widetilde{V_{+}}$ for any $\bv \in \widetilde{V_{+}}$. Hence, $E_{+}|_{\widetilde{V_{+}}}$ is a projection on $\widetilde{V_{+}}$. We show now that $\ker E_{+}|_{\widetilde{V_{+}}} = \{\bO\}$ such that  $E_{+}|_{\widetilde{V}_{+}} = \id_{V_{+}}$ and hence  $\widetilde{V_{+}}\subset\ran E_{+} = V_{+}$. We do this by constructing a slice hyperholomorphic extension of $\RRes_{s}(T;\bv)$ that is defined on all of $\hh$ and applying \Cref{EmptySpec} for any $\bv\in \ker E_{+}|_{\widetilde{V_{+}}}$.

 Let  $\bv \in \ker E_{+}|_{\widetilde{V_{+}}}$. As $\ker E_{+}|_{\widetilde{V_{+}}}\subset \ker E_{+} = \ran E_{-} = V_{-}$, we find  $\bv\in V_{-}$. For $z = z_0 + z_1\I \in \cc_{\I}$, we define the function
\[
\bof_{\I}(z;\bv) :=\begin{cases} R_{z}(T_{-})\bv, & z_{1} \geq 0 \\ \left(z_0 \id + z_1 \widetilde{J} - T)\right)^{-1}\bv, & z_1 < 0.
\end{cases}
\]
This function is (right) holomorphic on $\cc_{\I}$. On $\cc_{\I}^{\geq}$ this is obvious because the resolvent of $T_{-}$ is a holomorphic function. For $z_1 <0$, we have
\begin{align*}
&\frac{1}{2}\left(\frac{\partial}{\partial z_0} \bof_{\I}(z;\bv) + \frac{\partial}{\partial z_1} \bof_{\I}(z;\bv)\I\right) \\
=& \frac{1}{2}\left(-\left(z_0 \id + z_1 \widetilde{J} - T)\right)^{-2}\bv - \left(z_0 \id + z_1 \widetilde{J} - T)\right)^{-2}\widetilde{J}\bv\I\right)\\
=&  \frac{1}{2}\left(-\left(z_0 \id + z_1 \widetilde{J} - T)\right)^{-2}\bv - \left(z_0 \id + z_1 \widetilde{J} - T)\right)^{-2}\bv\I^2\right) = \bO,
\end{align*}
as $\widetilde{J}\bv = \bv\I$ because $\bv\in \widetilde{V_{+}} = V_{\widetilde{J},\I}^{+}$. The slice extension $\bof(s;\bv)$ of $\bof_{\I}(s;\bv)$ obtained from \Cref{extLem} is a slice hyperholomorphic extension of $\RRes_{s}(T;\bv)$ to all of $\hh$ in the sense of \Cref{SHExt}. Indeed, as $\Q_{z}(T)|_{V_{-}} = \Q_{z}(T_{-}) = (\id_{V_{-}}\overline{z} - T_{-})(\id_{V_{-}}z -T_{-})$, we find for $s\in\cc_{\I}^{\geq}$ that
\begin{align*}
&\Q_{s}(T)\bof(s;\bv) = \Q_{s}(T_{-}) \bof_{\I}(s;\bv) \\
=& (\overline{s} \id_{V_{-}} - T_{-})(s \id_{V_{-}} -T_{-})R_{s}(T_{-})\bv\\
 = &(\overline{s} \id_{V_{-}} - T_{-})\bv = \bv s - T_{-}\bv = \bv s - T\bv.
\end{align*}
On the other hand, the facts that $T$ and $\widetilde{J}$ commute and that $-\widetilde{J}^2 = \id$ because $\widetilde{J}$ is an imaginary operator with $\ran\widetilde{J} = V_{R}$ imply
\begin{align*}
&\left(s_0\id  + s_1 \widetilde{J} - T \right)\left(s_0\id  - s_1 \widetilde{J} - T \right) \\
=&  s_0^2\id - s_0s_{1} \widetilde{J} - s_0 T + s_0s_1\widetilde{J}-s_1^2\widetilde{J}^2 - s_{1}\widetilde{J}T  -s_0T + s_1 T\widetilde{J} + T^2 \\
=& |s|^2\id  -2 s_0 T   + T^2 = \Q_{s}(T).
\end{align*}
 For $s = s_1 + (-\I)s_1\in \cc_{\I}^{-}$, we find thus because of $\bv \in\widetilde{V_{+}} = V_{\widetilde{J},\I}^{+}$ that
\begin{align*}
\Q_{s}(T)\bof(s;\bv) =& \left(s_0\id  + s_1 \widetilde{J} - T \right)\left(s_0\id  - s_1 \widetilde{J} - T\right)\bof_{\I}(s;\bv) \\
=& \left(s_0\id  + s_1 \widetilde{J} - T \right)\left(s_0\id  - s_1 \widetilde{J} - T \right)\left(s_0\id  - s_1 \widetilde{J} - T \right)^{-1} \bv\\
=& \left(s_0\id  + s_1 \widetilde{J} - T \right) \bv = \bv s_{0} + \bv \I s_{1} - T = \bv\overline{s} - T\bv.
\end{align*}
Finally, for $s\notin\cc_{\I}$,  \Cref{RepFo} yields
\begin{align*}
 \Q_{s}(T)\bof(s;\bv) =&  \Q_{s}(T)\bof_{\I}(s_{\I};\bv) (1 - \I\I_{s}) \frac{1}{2} +  \Q_{s}(T) \bof_{\I}\left(\overline{s_{\I}};\bv\right) (1 + \I\I_{s})\frac{1}{2} \\
=& \left(\bv \overline{s_{\I}} - T\bv\right)  (1 - \I\I_{s}) \frac{1}{2}+  (\bv s_{\I} - T\bv)  (1 + \I\I_{s})\frac{1}{2}\\
=& \bv \left(\overline{s_{\I}}(1-\I\I_{s}) + s(1+\I\I_{s})\right)\frac{1}{2} - T\bv\left((1 - \I\I_{s}) + (1+\I\I_{s})\right) \frac{1}{2}\\
=& \bv(s_{\I} + \overline{s_{\I}} + (s_{\I} - \overline{s_{\I}})\I\I_{s})\frac{1}{2} - T\bv = \bv (s_0 - s_1\I_{s}) - T\bv = \bv \overline{s} - T\bv.
\end{align*}
From \Cref{EmptySpec}, we hence deduce that $\bv = \bO$ and so  $\ker E_{+}|_{\widetilde{V_{+}}} = \{\bO\}$. Since $E_{+}|_{\widetilde{V_{+}}}$ is a projection on $\widetilde{V_{+}}$, we have $\widetilde{V_{+}} = \ker E_{+}|_{\widetilde{V_{+}}} \oplus \ran E_{+}|_{\widetilde{V_{+}}} = \{\bO\} \oplus \ran E_{+}|_{\widetilde{V_{+}}}$. We conclude $\widetilde{V_{+}} = \ran E_{+}|_{\widetilde{V_{+}}} \subset \ran E_{+} = V_{+}$. We therefore have
\[
 V_{R} = \widetilde{V_{+}}\oplus \widetilde{V_{+}}\J \subset V_{+} \oplus V_{+}\J = V_{R}.
\]
This implies $V_{+}   = \widetilde{V_{+}}$ and in turn $J = \widetilde{J}$.

\end{proof}
\begin{corollary}\label{SSOr2}
Let $T\in\boundOP(V_R)$ and assume that there exists a spectral resolution for $T$ as in \Cref{SSOr}. If $\sigma_{S}(T) = \Delta_{1}\cup\Delta_{2}$ with closed sets $\Delta_{1},\Delta_{2}\in\sBorel(\hh)$ such that $\Delta_{1}\subset\rr$ and $\Delta_{2}\cap\rr = \emptyset$, then there exists a unique imaginary operator $J\in\boundOP(V_R)$ that is a spectral orientation for $T$ such that $T$ is a spectral operator with spectral decomposition $(E,J)$.
\end{corollary}
\begin{proof}
Let $T_{2} = T_{2}|_{V_2}$, where $V_{2} = \ran E(\hh\setminus \rr) = \ran E(\Delta_{2})$. Then $E_{2}(\Delta) :=  E(\Delta)|_{V_2}$ for $\Delta\in\sBorel(\hh)$ is by \Cref{SubSpCor} a spectral resolution for $T_2$. Since $\sigma_{S}(T_2) \subset \Delta_{2}$ and  $\Delta_{2}\cap\rr = \emptyset$,  \Cref{SSOr} implies the existence of a unique spectral orientation $J_2$ for $T_{2}$. 

The fact that $(E_2,J_2)$ is a spectral system implies $\ran J_2 = \ran E_{2}(\hh\setminus\rr)V_{2} = V_{2}$ because $E_{2}(\hh\setminus\rr) = E(\hh\setminus\rr)|_{V_2} = \id_{V_2}$. If we set $J = J_2 E(\hh\setminus\rr)$,  we find that $\ker J = \ran E(\rr)$ and $\ran J = V_{2} = \ran E(\hh\setminus\rr)$. We also have
\begin{align*}
E(\Delta) J =&  E(\Delta\cap\rr) J_2 E(\hh\setminus \rr)  + E(\Delta\setminus\rr)J_2 E(\hh\setminus\rr) \\
= &E_2(\Delta\setminus\rr) J_2 E(\hh\setminus\rr) = J_2 E_2(\Delta\setminus\rr) E(\hh\setminus\rr) \\
= &  J_2E(\Delta\setminus\rr)E(\hh\setminus\rr) =  J_2E(\hh\setminus\rr)E(\Delta\setminus\rr) = J E(\Delta),
\end{align*}
where the last identity used that $E(\hh\setminus\rr)E(\Delta\cap\rr) = 0$. Moreover, we have
\[ - J^2 = - J_2 E(\hh\setminus\rr)J_2E(\hh\setminus\rr) = - J_2^2 E(\hh\setminus\rr) = E(\hh\setminus\rr)\]
so that $-J^2$ is a projection onto $\ran J = \ran E(\hh\setminus\rr)$ along $\ker J = \ran E(\rr)$. Hence, $J$  is  an imaginary operator and $(E,J)$ is a spectral system on $V_{R}$. Finally, for any $s_0,s_1\in\rr$ with $s_1 >0$, we have
\[
\left((s_0\id - s_1 J - T)|_{V_2}\right)^{-1} = (s_0 \id_{V_{2}} - s_1 J_2 - T_2)^{-1} \in\boundOP(V_2)
\]
and hence $(E,J)$ is actually a spectral decomposition for $T$. 

In order to show the uniqueness of $J$ we consider an arbitrary spectral orientation $\widetilde{J}$  for $T$. Then 
\begin{equation}\label{TumPO}
\ker\widetilde{J} =  E(\rr)V_R = \ker J \quad\text{and}\quad \ran\widetilde{J} = E(\hh\setminus\rr)V_{R} = \ran J.
\end{equation}
By \Cref{SubSpCor}, the operator $\widetilde{J}|_{V_2}$ is a spectral orientation for $T_2$. The spectral orientation of $T_2$ is however unique by \Cref{SSOr} and hence $\widetilde{J}|_{V_2} = J_2 = J|_{V_2}$. We conclude $\widetilde{J} = J$.

\end{proof}

Finally we can now show the uniqueness of the spectral orientation of an arbitrary spectral operator.

\begin{theorem}\label{UniqueSpecDec}
The spectral decomposition $(E,J)$ of a spectral operator $T\in\boundOP(V_{R})$ is uniquely determined
\end{theorem}
\begin{proof}
The uniqueness of the spectral resolution $E$ has already been shown in \Cref{EUnique}. Let $J$ and $\widetilde{J}$ be two spectral orientations for $T$. Since \eqref{TumPO} holds true also in this case, we can reduce the problem to showing $J|_{V_1} = \widetilde{J}|_{V_1}$ with $V_{1}: = \ran E(\hh\setminus\rr)$. The operator $T_1:= T|_{V_{1}}$ is a spectral operator on $V_1$. By \Cref{SubSpCor}, $(E_1,J_1)$ and $(E_1, \widetilde{J}_1)$ with $E_{1}(\Delta) = E(\Delta)|_{V_1}$ and $J_1 = J|_{V_{1}}$ and $\widetilde{J}_1:= \widetilde{J}|_{V_1}$ are spectral decompositions of $T_1$. As $E_{0}(\rr) = 0$, it is hence sufficient to show the uniqueness of the spectral orientation of a spectral operator under the assumption $E(\rr) = 0$.

Let hence $T$ be a spectral operator with spectral decomposition $(E,J)$ such that $E(\rr) = 0$. If $\dist(\sigma_{S}(T),\rr) > 0$, then we already know that the statement holds true. We have shown this in \Cref{SSOr}. Otherwise we choose a sequence of  pairwise disjoint sets $\Delta_{n}\in\sBorel(\hh)$ with $\dist(\Delta_{n}, \rr)>0$ that cover $\sigma_{S}(T)\setminus \rr$. We can choose for instance
\[
 \Delta_{n} := \left\{ s \in \hh : - \| T \| \leq s_0 \leq \|T\|, \  \frac{\| T\|}{n+1}  < s_1 \leq \frac{\| T\|}{n}  \right\}.
\]
By \Cref{Support} and as $E(\rr) = 0$, we have 
\[
E(\sigma_{S}(T)\setminus\rr) = E(\sigma_{S}(T)\setminus\rr) + E(\sigma_{S}(T)\cap\rr) = E(\sigma_{S}(T)) = \id
\]
We therefore find $\sum_{n=0}^{\infty} E(\Delta_n)\bv = E\left(\bigcup_{n\in\nn}\Delta_n\right)\bv = \bv$ for any $\bv\in V_{R}$ as $\sigma_{S}(T)\setminus\rr \subset \bigcup_{n\in\nn}\Delta_n$. 

Since $E(\Delta_{n})$ and $J$ commute, the operator $J$ leaves $V_{\Delta_n} := \ran  E(\Delta_{n})$ invariant. Hence $J_{\Delta_{n}} = J|_{V_{\Delta_{n}}}$ is a bounded operator on $V_{\Delta_n}$ and we have
 \[
 J\bv = J\sum_{n=0}^{+\infty} E(\Delta_{n})\bv = \sum_{n=1}^{+\infty}JE(\Delta_n)\bv = \sum_{n=1}^{+\infty}J_{\Delta_n}E(\Delta_n)\bv.
 \]
Similarly, we see that also $\widetilde{J_{\Delta_n}}:= \widetilde{J}|_{V_{\Delta_n}}$ is a bounded operator on $V_{\Delta_n}$ and that  $\widetilde{J}\bv  = \sum_{n=1}^{+\infty}\widetilde{J}_{\Delta_n}E(\Delta_n)\bv$. 
 
 Now observe that $T_{\Delta_{n}}$ is a spectral operator. Its spectral resolution is given by $E_{n}(\Delta) : = E(\Delta)|_{V_{\Delta_n}}$ for $\Delta\in\sBorel(\hh)$ as one can check easily. Its spectral orientation is given by $J_{\Delta_n}$: for any $\Delta\in\sBorel(\hh)$, we have
 \[
 E_{n}(\Delta) J_{\Delta_{n}} E(\Delta_n)  = E(\Delta)J E(\Delta_n) = J E(\Delta) E(\Delta_{n}) = J_{\Delta_n}E_{n}(\Delta)E(\Delta_n)
 \]
and hence $E_{n}(\Delta) J_{\Delta_{n}} = J_{\Delta_{n}} E(\Delta_n)$ on $V_{\Delta_n}$. Since $\ker J_{\Delta_n} = \{\bO\} = E_n(\rr)$ and $\ran J_{\Delta_n} = V_{\Delta_n} = E_{n}(\hh\setminus\rr)$, the couple $(E, J_{\Delta_n})$ is actually a spectral system. Furthermore, the operators $T_{\Delta_n}$ and $J_{\Delta_n}$ commute as
\[
T_{\Delta_n} J_{\Delta_n}E(\Delta_n) = T J E(\Delta_n) = JT E(\Delta_n) = J_{\Delta_n}T_{\Delta_n}E(\Delta_n).
\]
Finally, for any $s_0,s_1\in\rr$ with $s_1>0$, we find 
 \[
 (s_0\id_{V_{\Delta_n}} - s_1 J_{\Delta_n} - T_{\Delta_n})^{-1} = (s_0\id - s_1J - T)^{-1}|_{V_{\Delta_n}}
 \]
  such that $(E_n, J_{\Delta_n})$ is actually a spectral decomposition for $T_{\Delta_n}$.  However, the same arguments show that also $\left(E_{n}, \widetilde{J}_{\Delta_n}\right)$ is a spectral decomposition for $T_{\Delta_n}$. Since $\sigma_{S}\left(T_{\Delta_n}\right)\subset \Delta_n$  and $\dist(\Delta_n,\rr)>0$, \Cref{SSOr} implies that the spectral orientation of $T_{\Delta_n}$ is unique such that $J_{\Delta_n} = \widetilde{J}_{\Delta_n}$. We thus find
\[
 J\bv  = \sum_{n=1}^{+\infty}J_{\Delta_n}E(\Delta_n)\bv = \sum_{n=1}^{+\infty}\widetilde{J}_{\Delta_n}E(\Delta_n)\bv = \widetilde{J}\bv.
\]

\end{proof}

\begin{remark}
In \Cref{SSOr} and \Cref{SSOr2} we showed that under certain assumptions the existence of a spectral resolution $E$ for $T$  already implies the existence of a spectral orientation and is hence a sufficient condition for $T$ to be a spectral operator. One may wonder whether this is true in general. An intuitive approach for showing this follows the idea of the proof of \Cref{UniqueSpecDec}. We can cover $\sigma_{S}(T)\setminus\rr$ by pairwise disjoint sets $\Delta_n\in\sBorel(\hh)$ with $\dist(\Delta_n,\rr)>0$ for each $n\in\nn$.  On each of the subspaces $V_{n} := \ran E(\Delta_n)$, the operator $T$ induces the operator $T_n := T|_{V_n}$ with $\sigma_S(T_n)\subset \clos{\Delta_{n}}$. Since $\dist(\Delta_n,\rr)>0$, we can then define $\Delta_{n,+} := \Delta_n\cap\cc_{\I}^{+}$ and $\Delta_{n,-}:= \Delta_n\cap\cc_{\I}^{-}$ for an arbitrary imaginary unit $\I\in\SS$ and consider the Riesz-projectors $E_{n,+} := \chi_{\Delta_{n,+}}(T_n)$ and $E_{n,-}:=\chi_{\Delta_{n,-}}(T_n)$ of $T_n$ on $V_{n,\I}$ associated with $\Delta_{n,+}$ and $\Delta_{n,-}$. Just as we did it in the proof of \Cref{SSOr}, we can then construct a spectral orientation for $T_n$ by setting $J_{n}\bv = E_{n,+}\bv \I + E_{n,-}\bv(-\I)$ for $\bv \in V_{n}$. The spectral orientation of $J$ must then be 
\begin{equation}\label{FailSpecOr}
J\bv = \sum_{n=1}^{+\infty}J_{n}E(\Delta_n) \bv = \sum_{n=1}^{+\infty}E_{n,+}E(\Delta_n)\bv\I + E_{n,-}E(\Delta_n)\bv(-\I).
\end{equation}

If $T$ is a spectral operator, then  $E_{n,+}= E_+ |_{V_n}$ and $E_{n,-} = E_{-}|_{V_n}$, where $E_+$ and $E_-$ are as usual the projections of $V_R$ onto $V_{J,\I}^{+}$ and $V_{J,\I}^{-}$ along $V_{0}\oplus V_{J,\I}^{-}$ resp.  $V_{0}\oplus V_{J,\I}^{+}$. Hence the Riesz-projectors $E_{n,+}$ and $E_{n,-}$ are uniformly bounded in $n\in\nn$ and the above series converges. The spectral orientation of $T$ can therefore be constructed as described above if $T$ is a spectral operator.

This procedure however fails if the Riesz-projectors $E_{n,+}$ and $E_{n,-}$ are not uniformly bounded because the convergence of the above series  is  in this case not guaranteed. The next example presents an operator for which the above series does actually not converge for this reason although the operator has a quaternionic spectral resolution. Hence, the existence of a spectral resolution does in general not imply the existence of a spectral orientation.
\end{remark}
\begin{example}\label{CEx}
Let $\ell^{2}(\hh)$ be the space of all square-summable sequences with quaternionic entries and choose $\I,\J\in\SS$ with $\I\perp\J$. We define an operator $T$ on $\ell^2(\hh)$ by the following rule: if $(b_n)_{n\in\nn} = T(a_n)_{n\in\nn}$, then
\begin{equation}\label{ZumPSTi}
\begin{pmatrix} b_{2m-1}\\ b_{2m}\end{pmatrix} = \frac{1}{m^2}\begin{pmatrix} \I & 2m\I \\ 0 & -\I \end{pmatrix}\begin{pmatrix}a_{2m-1} \\ a_{2m}\end{pmatrix}.
\end{equation}
For neatness, let us denote the matrix in the above equation by $J_m$ and let us set $T_{m} := \frac{1}{m^2} J_m$, that is
\[
J_{m} := \begin{pmatrix} \I & 2m\I \\ 0 & -\I \end{pmatrix}\qquad\text{and}\qquad T_m := \frac{1}{m^2} \begin{pmatrix} \I & 2m\I \\ 0 & -\I \end{pmatrix}.
\]
Since all matrix norms are equivalent, there exists a constant $C>0$ such that 
\begin{equation}\label{MatNormEquiv}
\| M \| \leq C \max_{\ell,\kappa \in \{1,2\} }\left|m_{\ell,\kappa}\right| \qquad \forall M = \begin{pmatrix} m_{1,1} & m_{1,2} \\ m_{2,1} & m_{2,2}\end{pmatrix} \in\hh^{2\times 2}
\end{equation}
 such that $\|J_{m}\| \leq 2Cm$. We thus find for \eqref{ZumPSTi} that 
\[
\| (b_{2m-1}, b_{2m})^T \|_{2} \leq \frac{2C}{m}\|(a_{2m-1},a_{2m})^T\|_{2} \leq 2C \|(a_{2m-1},a_{2m})^T\|_{2} .
\]
and in turn 
\begin{equation}\label{TzDM}
\begin{split}
\| T (a_n)_{n\in\nn}\|_{\ell^2(\hh)}^2 & = \sum_{m=1}^{+\infty} |b_{2m-1}|^2 + |b_{2m}|^2\\
&\leq \sum_{m=1}^{+\infty} 4C^2 \left(|a_{2m-1}|^2 + |a_{2m}|^2\right) = 4C^2 \left\|(a_{n})_{n\in\nn}\right\|_{\ell^2(\hh)}^2.
\end{split}
\end{equation}
Hence $T$ is a bounded right-linear operator on $\ell^2(\hh)$.

We show now that the $S$-spectrum of $T$ is the set $\Lambda = \{0\} \cup \cup_{n\in\nn} \frac{1}{n^2}\SS$. For $s\in\hh$, 
 the operator $\Q_{s}(T) = T^2 - 2s_0 T + |s|^2$ is given by the following relation: if  $(c_n)_{n\in\nn} = \Q_{s}(T)(a_n)_{n\in\nn}$ then
\begin{equation}\label{JuHKA}
\begin{pmatrix} c_{2m-1}\\ c_{2m}\end{pmatrix} = \begin{pmatrix} - \frac{1}{m^2} - 2\I \frac{s_0}{m^2} + |s|^2 & -4\I \frac{s_0}{m} \\ 0 &  - \frac{1}{m^2} - 2\I \frac{s_0}{m^2} + |s|^2  \end{pmatrix}\begin{pmatrix}a_{2m-1} \\ a_{2m}\end{pmatrix}.
\end{equation}
The inverse of the above matrix is
\[
\begin{split}
\Q_{s}(T_m)^{-1} = &
\begin{pmatrix}
\frac{m^4}{|s|^2 m^4-2 i s_0 m^2-1} & \frac{4 i m^7 s_0}{|s|^4 m^8+2 \left(s_0^2-s_1^2\right) m^4+1} \\
 0 & \frac{m^4}{|s|^2 m^4+2 i s_0 m^2-1} \\
\end{pmatrix}
\\
= &\begin{pmatrix}
 \frac{1}{\left(s_{\I} - \frac{\I}{m^2}\right)\left(\overline{s_{\I}} - \frac{\I}{m^2}\right) }& \frac{4 i  s_0}{m \left( s_{\I} + \frac{\I}{m^2}\right)\left( s_{\I} - \frac{\I}{m^2}\right)\left( \overline{s_{\I}} + \frac{\I}{m^2}\right)\left( \overline{s_{\I}} - \frac{\I}{m^2}\right)} \\
 0 & \frac{1}{\left( s_{\I} +\frac{\I}{m^2}\right)\left( \overline{s_{\I}} +\frac{\I}{m^2}\right)} 
\end{pmatrix}.
\end{split}
\]
with $s_{\I} = s_0 + \I s_{1}$. Hence, $\Q_{s}(T_m)^{-1}$ exists for $s_{\I}\neq \frac{1}{m^2}\I$. We have 
\[
\left| s_{\I}  - \frac{\I}{m^2}\right| \left|\overline{s_{\I}} - \frac{\I}{m^2}\right|  = \left| s_{\I}  + \frac{\I}{m^2}\right| \left|\overline{s_{\I}} + \frac{\I}{m^2}\right| \geq 2 \left|s_{\I} - \frac{\I}{m^2}\right|  = 2 \dist\left(s, \left[\frac{\I}{m}\right]\right)
\]
and so 
\begin{equation}\label{RSrs}
\left\| \Q_{s}(T_m)^{-1} \right\|  \leq  C \max\left\{
 \frac{1}{2\dist\left(s, \left[\frac{\I}{m^2}\right]\right) }  ,  \frac{|s_0|}{ m \left(\dist\left(s, \left[\frac{\I}{m^2}\right]\right) \right)^2}
 \right\},
\end{equation}
where $C$ is the constant in \eqref{MatNormEquiv}.  If $s\notin\Lambda$, then $0<\dist(s,\varLambda) \leq \dist\left(s, \left[\frac{\I}{m^2}\right]\right)$ and hence the matrices $\Q_{s}(T_m)^{-1}$ are for $m\in\nn$ uniformly bounded by
\[
\left\| \Q_{s}(T_m)^{-1} \right\|  \leq  C \max\left\{
 \frac{1}{2\dist\left(s, \Lambda \right) }  ,  \frac{|s_0|}{ \dist\left(s, \Lambda \right)^2}
 \right\}.
\]
The operator $\Q_{s}(T)^{-1}$ is then given by the relation
\begin{equation}
\begin{pmatrix} a_{2m-1} \\ a_{2m}\end{pmatrix} = \Q_{s}(T_m)^{-1} \begin{pmatrix}c_{2m-1}\\ c_{2m}\end{pmatrix},
\end{equation}
for $(a_n)_{n\in\nn} = Q_{s}(T)^{-1}(c_n)_{n\in\nn}$. A computation similar to the one in \eqref{TzDM} shows that this operator is bounded on $\ell^2(\hh)$. Thus $s\in\rho_{S}(T)$ if $s\notin\Lambda$ and in turn $\sigma_{S}(T)\subset \Lambda$.

For any $m\in\nn$, we set $s_m = \frac{1}{m^2}\I$. The sphere  $[s_{m}] = \frac{1}{m^2} \SS$ is an eigensphere of $T$ and the associated eigenspace $V_{m}$ is the right-linear span of $\mathbf{e}_{2m-1}$ and $\mathbf{e}_{2m}$, where $\mathbf{e}_{n} = (\delta_{n,\ell})_{\ell\in\nn}$, as one can see easily from \eqref{JuHKA}. A straightforward computation moreover shows that the vectors $\bv_{2m-1} := \mathbf{e}_{2m-1}$ and $\bv_{2m}:= -\mathbf{e}_{2m-1}\J + \frac{1}{m} \mathbf{e}_{2m}\J$ are eigenvectors of $T$ with respect to the eigenvalue $s_{m}$. Hence $[s_m] \subset \sigma_{S}(T)$. Since $\sigma_{S}(T)$ is closed, we finally find $\Lambda = \clos{\bigcup_{m\in\nn} [s_m] } \subset \sigma_{S}(T)$ and in turn $\sigma_{S}(T) = \Lambda$.

Let $E_{m}$ for $m\in\nn$ be the orthogonal projection of $\ell^2(\hh)$ onto  $V_{m}:=\linspan{\hh}\{\mathbf{e}_{2m-1},\mathbf{e}_{2m}\}$, that is $E_{m}(a_{n})_{n\in\nn} = \mathbf{e}_{2m-1}a_{2m-1} + \mathbf{e}_{2m}a_{2m}$. We set for $\Delta\in\sBorel(\hh)$
\[
E(\Delta) = \sum_{m\in I_{\Delta}} E_{m} \qquad \text{with }\qquad I_{\Delta}:=\left\{m\in\nn: \frac{1}{m^2}\SS \subset \Delta\right\}.
\]
It is immediate that $E$ is a spectral measure on $\ell^2(\hh)$, that $\|E(\Delta)\| \leq 1$ for any $\Delta\in\sBorel(\hh)$ and that $E(\Delta)$ commutes with $T$ for any $\Delta\in\sBorel(\hh)$. Moreover, if $s \notin\clos{\Delta}$, then the pseudo-resolvent $\Q_{s}(T_{\Delta})^{-1}$ of $T_{\Delta} = T|_{V_{\Delta}}$ with $V_{\Delta} = \ran E(\Delta)$ is given by
\[ \Q_{s}(T_{\Delta})^{-1} = \left.\left( \sum_{m\in I_{\Delta}} \Q_{s}(T_m)^{-1}E_{m}\right)\right|_{\ran E(\Delta)}.\]
Since $0 < \dist\left(s,\bigcup_{m\in I_{\Delta}} \left[\frac{\I}{m^2}\right]\right) = \inf_{m\in I_{\Delta}} \dist\left(s, \left[\frac{\I}{m^2}\right]\right)$, the operators $\Q_{s}(T_{m})^{-1}$ are uniformly bounded for $m\in I_{\Delta}$. Computations similar to \eqref{TzDM} show that $\Q_{s}(T_{\Delta})^{-1}$ is a bounded operator on $V_{\Delta}$. Hence, $s\in\rho_{S}(T_{\Delta})$ and in turn $\sigma_{S}(T_{\Delta})\subset \clos{\Delta}$. Altogether we obtain that $E$ is a spectral resolution for $T$.

In order to construct a spectral orientation for $T$, we first observe that $J_m$ is a spectral orientation for $T_m$. For $s_0,s_1\in\rr$ with $s_1>0$, we have
\[
s_0\id_{\hh^2} - s_1J_{m} -T_m = \begin{pmatrix} s_0  - \left(s_1 + \frac{1}{m^2}\right)\I & -\left(s_1 + \frac{1}{m^2}\right) 2m\I\\ 0 & s_0 + \left(s_1 + \frac{1}{m^2}\right)\I\end{pmatrix},
\]
the inverse of which is given by the matrix
\[
\left( s_0\id_{\hh^2} - s_1J_{m} -T_m\right)^{-1} = \begin{pmatrix} \frac{1}{s_0-\left( s_1 + \frac{1}{m^2}\right)\I} & \frac{2\I m\left(\frac{1}{m^2}+s_{1}\right)}{s_0^2 + \left(\frac{1}{m^2} + s_1\right)^{2}
}\\
 0 & \frac{1}{s_0+ \left(\frac{1}{m^2}+s_1\right)\I}
\end{pmatrix}.
\]
Since $s_1>0$, each entry has non-zero denominator and hence we have  that $\left( s_0\id_{\hh^2} - s_1J_{m} -T_m\right)^{-1}\in\boundOP(\hh^2)$. 

If $J\in\boundOP(\ell^2(\hh))$ is a spectral orientation for $T$, then the restriction $J|_{V_{m}}$ of $J$ to $V_{m} = \linspan{\hh}\{\mathbf{e}_{2m-1}, \mathbf{e}_{2m}\}$ is also a spectral orientation for $T_{m}$. The uniqueness of the spectral orientation implies $J|_{V_m} = J_m$ and hence $J =  \sum_{m=1}^{+\infty} J|_{V_{m}} E\left(\frac{1}{m^2}\SS\right) =  \sum_{m=1}^{+\infty} J_m E_m$.
This series does however not converge because the operators $J_{V_m}$ are not uniformly bounded. Hence, it does not define a bounded operator on $\ell^2(\hh)$. Indeed, the sequence $a_{2m-1} = 0$, $a_{2m} = m^{-\frac{3}{2}}$ for instance belongs to $\ell^2(\hh)$, but 
\begin{gather*}
\left\| \sum_{m=1}^{+\infty} J_{m} E_{m} (a_{n})_{n\in\nn}\right\|_{\ell^2(\hh)}^2 = \sum_{m=1}^{+\infty} \left\|  \begin{pmatrix} \I & 2m\I \\ 0 & -\I\end{pmatrix}\begin{pmatrix} 0 \\ m^{-\frac{3}{2}}\end{pmatrix} \right\|_{2}^{2} \\
=2  \sum_{m=1}^{+\infty} 4\frac{1}{m} + \frac{1}{m^3} = +\infty.
\end{gather*}
Hence there cannot exist a spectral orientation for $T$ and in turn $T$ is not a spectral operator on $\ell^{2}(\hh)$.

We conclude this example with a remark on its geometric intuition. Let us identify $\hh^2 \cong \cc_{\I}^4$, which is for any $\J\in\SS$ with $\J\perp\I$ spanned by the basis vectors 
\[
\mathbf{b}_1 = \begin{pmatrix}1 \\0 \end{pmatrix},\quad \mathbf{b}_2 = \begin{pmatrix}\J \\0 \end{pmatrix},\quad \mathbf{b}_{3} = \begin{pmatrix}0 \\ 1 \end{pmatrix}\quad \text{and} \quad \mathbf{b}_4 = \begin{pmatrix}0 \\\J \end{pmatrix}.
\]
The vectors $\bv_{m,1} = \bb_{1}$ and $\bv_{m,2} = -\bb_{2} + \frac{1}{m}\bb_{4}$ are eigenvectors of $J_{m}$ with respect to $\I$ and the vectors $\bv_{1}\J = \bb_{2}  $ and $\bv_{m,2} = \bb_{1} -\frac{1}{m}\bb_{3}$ are eigenvectors of $J_{m}$ with respect to $-\I$. We thus find $V_{J_{m},\I}^{+} = \linspan{\cc_{\I}}\{\bb_{1}, -\bb_{2} + \frac{1}{m}\bb_{4}\}$ and $V_{J_{m},\I}^{-} = V_{J_{m},\I}^{+}\J$. However, as $m$ tends to infinity, the vector $\bv_{2}$ tends to $\bv_{1}\J$  and $\bv_{2}\J$ tends to $\bv_{1}$. Hence, intuitively, in the limit $V_{J_{m},\I}^{-} = V_{J_{m},\I}^{+}\J = V_{{J_m},\I}^{+}$ and consequently the projections of $\hh^2 = \cc_{\I}^4$ onto $V_{J_m,\I}^{+}$ along $V_{J_m,\I}^{-}$ become unbounded.
\end{example}
Finally, the notion of quaternionic spectral operator is again backwards compatible with complex theory on $V_{R,\I}$. 
\begin{theorem}\label{SOpBackwards}
An operator $T\in\boundOP(V_{R})$ is a quaternionic spectral operator if and only if it is a spectral operator on $V_{R,\I}$ for some (and hence any) $\I\in\SS$. (See \cite{Dunford:1958} for the complex theory.) If furthermore $(E,J)$ is the quaternionic spectral decomposition of $T$ and $E_{\I}$ is the spectral resolution of $T$ as a complex $\cc_{\I}$-linear operator on $V_{R,\I}$, then 
\begin{equation}\label{QWerID} 
\begin{split}
&E(\Delta) = E_{\I}(\Delta\cap\cc_{\I})\quad\forall\Delta\in\sBorel(\hh)\\
 &J\bv = E_{\I}(\cc_{\I}^{+})\bv \I + E_{\I}(\cc_{\I}^{-})\bv(-\I)\quad\forall\bv\in V_{R}.
 \end{split}
\end{equation}
Conversely, $E_{\I}$ is the spectral measure on $V_{R}$ determined by $(E,J)$ that was constructed in \Cref{BackCompat}.
\end{theorem}
\begin{proof}
Let us first assume that $T\in\boundOP(V_{R})$ is a quaternionic spectral operator with spectral decomposition $(E,J)$ in the sense of \Cref{SpecOP} and let $\I\in\SS$. Let $E_{+}$ be the projection of $\ran J = V_{J,\I}^{+}\oplus V_{J,\I}^{-}$ onto $V_{J,\I}^{+}$ along $V_{J,\I}^{-}$ and let $E_{-}$ be the projection of $\ran J$ onto $V_{J,\I}^{-}$ along $V_{J,\I}^{+}$, cf. \Cref{JVSplitThm}. Since $T$ and $E(\Delta)$ for $\Delta\in\sBorel(\hh)$ commute with $J$, they leave the spaces $V_{J,\I}^{+}$ and $V_{J,\I}^{-}$ invariant and hence they commute with $E_{+}$ and $E_{-}$. By \Cref{BackCompat} the set function $E_{\I}$ on $\cc_{\I}$ defined in \eqref{CSpecMeas}, which is given by
\begin{equation}\label{RumSti}
  E_{\I}(\Delta) = E_{+}E\left(\left[\Delta\cap\cc_{\I}^{+}\right]\right) + E(\Delta\cap\rr) + E_{-}E\left(\left[\Delta\cap\cc_{\I}^{-}\right]\right),
\end{equation}
for $\Delta\in\Borel(\cc_{\I})$ is a spectral measure on $V_{R,\I}$. Since the spectral measure $E$ and the projections $E_{+}$ and $E_{-}$ commute with $T$,  the spectral measure $E_{\I}$ commutes with $T$ too. 

If $\Delta\in\Borel(\cc_{\I})$ is a subset of $\cc_{\I}^{+}$, then $ J\bv= \bv\I$ for $\bv\in V_{\I,\Delta} := \ran E_{\I}(\Delta)$ as $\ran E_{\I}(\Delta) = \ran (E_{+}E([\Delta])) \subset V_{J,\I}^{+}$. For $z = z_0 + \I z_1\in\cc_{\I}$ and $\bv\in V_{\I,\Delta}$, we thus have
\begin{gather*}
(z\id_{V_{\I,\Delta}} - T)\bv = \bv z_0  + \bv\I z_1 - T\bv \\
= \bv z_0 + J \bv z_1 - T\bv = (z_0\id_{V_{\I,\Delta}} + z_1J - T)\bv.
\end{gather*}
If $z\in\cc_{\I}^{-}$, then  the inverse of $(z_0\id_{V_{R,\I}} + z_{1}J - T)|_{\ran J}$ exists because $J$ is the spectral orientation of $T$. We thus have $R_{z}(T_{\Delta}) = (z_0\id_{V_{R,\I}} + z_{1}J - T)^{-1}|_{V_{\I,\Delta}}$  and so $\cc_{\I}^{-}\subset \rho(T_{\Delta})$. If on the other hand $z\in\cc_{\I}^{\geq}\setminus \clos{\Delta}$, then $z\in\rho_S(T_{[\Delta]})$ where $T_{[\Delta]} = T|_{ V_{[\Delta]}}$ with $V_{[\Delta]} = \ran E([\Delta])$. Hence, $\Q_{z}(T_{[\Delta]})$ has a bounded inverse on $V_{[\Delta]}$. By the construction of $E_{\I}$ we have $V_{\I,\Delta} = E_{+}V_{[\Delta]}$ and since $T_{[\Delta]}$ and $E_{+}$ commute $\Q_{z}(T_{[\Delta]})^{-1}$ leaves $V_{\I,\Delta}$ invariant so that $\left.\Q_{z}(T_{[\Delta]})^{-1}\right|_{V_{\I,\Delta}}$ defines a bounded $\cc_{\I}$-linear operator on $V_{\I,\Delta}$. Because of \Cref{ABC}, the resolvent of $T_{\Delta}$ at $z$ is therefore given by 
\[
R_{z}(T)\bv = \Q_{s}(T_{[\Delta]})^{-1}(\bv \overline{z} - T_{\Delta}\bv)\qquad \forall \bv\in V_{\I,\Delta}.
\] 

Altogether, we conclude $\rho(T_{\Delta}) \supset \cc_{\I}^{-} \cup \left(\cc_{\I}^{\geq}\setminus \clos{\Delta}\right) = \cc_{\I}\setminus\clos{\Delta}$ and in turn $\sigma(T_{\Delta}) \subset \clos{\Delta}$. Similarly, we see that $\sigma(T_{\Delta})\subset\clos{\Delta}$ if $\Delta\subset \cc_{\I}^{-}$. If on the other hand $\Delta\subset \rr$, then $E_{\I}(\Delta) = E(\Delta)$ so that $T_{\Delta}$ is a quaternionic linear operator with $\sigma_{S}(T_{\Delta})\subset\clos{\Delta}$. By \Cref{ABC}, we have $\sigma(T_{\Delta} )= \sigma_{\cc_{\I}}(T_{\Delta}) = \sigma_{S}(T)\subset\clos{\Delta}$. Finally, if $\Delta\in\Borel(\cc_{\I})$ is arbitrary and $z\notin\clos{\Delta}$, we can set $\Delta_{+}:= \Delta\cap\cc_{\I}^{+}$, $\Delta_{-} : = \Delta\cap\cc_{\I}^{-}$ and $\Delta_{\rr}: = \Delta\cap\rr$. Then $z$ belong to the resolvent sets of each the operators $T_{\Delta_{+}}$, $T_{\Delta_{-}}$ and $T_{\Delta_{\rr}}$ and we find
\[
R_{z}(T) = R_{z}(T_{\Delta_{+}})E_{\I}(\Delta_{+}) + R_{z}(T_{\Delta_{\rr}})E(\Delta_{\rr}) + R_{z}(T_{\Delta_{-}})E_{\I}(\Delta_{-}).
\]
We thus have $\sigma(T_{\Delta})\subset \clos{\Delta}$. Hence, $T$ is a spectral operator on $V_{R,\I}$ and  $E_{\I}$ is its ($\cc_{\I}$-complex) spectral resolution on $V_{R,\I}$.

Now assume that $T$ is a bounded quaternionic linear operator on $V_{R}$ and that for some $\I\in\SS$ there exists a $\cc_{\I}$-linear spectral resolution $E_{\I}$ for $T$ as a $\cc_{\I}$-linear operator on $V_{R,\I}$. Following Definition~6 of \cite[Chapter XV.2]{Dunford:1958}, an analytic extension of $R_{z}(T)\bv$ with $\bv\in V_{R,\I} = V_{R}$ is a holomorphic function $\bof$ defined on a set $\fdom(\bof)$ such that $(z\id_{V_{R,\I}} - T)\bof(z) = \bv $ for $z\in\fdom(\bof)$. The resolvent $\rho(\bv)$ is the domain of the unique maximal analytic extension of $R_{z}(T)\bv$ and the spectrum $\sigma(\bv)$ is the complement of $\rho(\bv)$ in $\cc_{\I}$. (We defined the quaternionic counterparts of these concepts in  \Cref{SHExt} and \Cref{SpecVec}.) Analogue to \Cref{SubSpec}, we have
\begin{equation}\label{tmp1}
E_{\I}(\Delta) V_{R,\I} = \{ \bv\in V_{R,\I} = V_{R}: \sigma(\bv)\subset\Delta\},\qquad\forall\Delta\in\Borel(\cc_{\I}).
\end{equation}

Let $\bv \in V_{R,\I}$, let $\J\in\SS$ with $\I\perp\J$ and let $\bof$ be the unique maximal analytic extension of $R_{z}(T)\bv$ defined on $\rho(\bv)$. The mapping $z\mapsto \bof\left(\overline{z}\right)\J$ is then holomorphic on $\overline{\rho(\bv)}$: for any $z\in\overline{\rho(\bv)}$, we have $\overline{z}\in\rho(\bv)$ and in turn
\[
\lim_{h\to 0} \left(\bof\left(\overline{z +h}\right)\J - \bof\left(\overline{z}\right)\J\right) h^{-1} =   \lim_{h\to 0} \left(\bof(\overline{z} +\overline{h})- \bof(\overline{z})\right) \overline{h}^{-1}\J = \bof'\left(\overline{z}\right)\J.
\]
 Since $T$ is quaternionic linear, we moreover have for $z\in\overline{\rho(\bv)}$ that
\[
\left(z\id_{V_{R,\I}} - T\right) (\bof\left(\overline{z}\right)\J) = \bof\left(\overline{z}\right)\J z - T\left(\bof\left(\overline{z}\right)\J\right) = \big(\bof\left(\overline{z}\right)\overline{ z } - T\left(\bof\left(\overline{z}\right)\right)\big)\J = \bv \J.
\]
Hence $z\mapsto\bof\left(\overline{z}\right)\J$ is an analytic extension of $R_{z}(T)(\bv\J)$ that is defined on $\overline{\rho(\bv)}$. Consequently $\rho(\bv\J) \supset \overline{\rho(\bv)}$ and in turn $\sigma(\bv\J)\subset\overline{\sigma(\bv)}$. If $\tilde{\bof}$ is the maximal analytic extension of $R_{z}(T)(\bv\J)$, then 
similar arguments show that $z\mapsto \tilde{\bof}\left(\overline{z}\right)(-\J)$ is an analytic extension of $R_{z}(T)\bv$. Since this function is defined on $\overline{\rho(\bv\J)}$, we find $\rho(\bv)\supset \overline{\rho(\bv\J)}$ and in turn $\sigma(\bv)\subset\overline{\sigma(\bv\J)}$. Altogether, we obtain $\sigma(\bv) = \overline{\sigma(\bv\J)}$ and $\tilde{\bof}(z) = \bof\left(\overline{z}\right)\J$. From \eqref{tmp1} we deduce
\begin{equation}\label{KaJuk1}
\begin{split}
 \ran E_{\I}\left(\overline{\Delta}\right) =& \left\{ \bv\in V_{R,\I} = V_{R}: \sigma(\bv)\subset\overline{\Delta}\right\} \\
 =& \left\{ \bv\J \in V_{R,\I} = V_{R}: \sigma(\bv)\subset \Delta\right\} = \left(\ran E_{\I}(\Delta) \right)\J.
\end{split}
\end{equation}

In order to construct the quaternionic spectral resolution of $T$, we define now
\[
E(\Delta) := E_{\I}(\Delta\cap\cc_{\I}),\qquad \forall\Delta\in\sBorel(\hh). 
\]
Obviously this operator is a bounded $\cc_{\I}$-linear projection on $V_{R} = V_{R,\I}$. We show now that it is also quaternionic linear. Due to the axial symmetry of $\Delta$, the identity \eqref{KaJuk1} implies  
\[
(\ran E(\Delta))\J = \left(\ran E_{\I} (\Delta\cap\cc_{\I}) \right)\J = \ran E_{\I}\left(\overline{\Delta\cap\cc_{\I}}\right) = \ran E_{\I}\left(\Delta\cap\cc_{\I}\right) = \ran E(\Delta).
\]
Similarly we find 
\begin{align*}
(\ker E(\Delta))\J =& (\ker E_{\I}(\Delta\cap\cc_{\I}))\J =   (\ran E_{\I}(\cc_{\I}\setminus\Delta))\J \\
=& \ran E_{\I}\left(\overline{\cc_{\I}\setminus\Delta}\right) =  \ran E_{\I}\left(\cc_{\I}\setminus\Delta\right) \\
=& \ker E_{\I}(\Delta\cap\cc_{\I}) = \ker E(\Delta). 
\end{align*}
If we write $\bv\in V_{R}$ as $\bv = \bv_{0} + \bv_{1}$ with $\bv_{0}\in\ker E(\Delta)$ and $\bv_{1}\in\ran E(\Delta)$, we thus have
\[
E(\Delta)(\bv\J) = E(\Delta)(\bv_{0}\J) + E(\Delta)(\bv_{1}\J) = \bv_{1}\J = (E(\Delta)\bv)\J.
\]
Writing $a\in\hh$ as $a= a_{1} + \J a_{2} $ with $a_1,a_2\in\cc_{\I}$  find due to the $\cc_{\I}$-linearity of $E(\Delta)$ that even
\[
E(\Delta)(\bv a) = (E(\Delta)\bv) a_{1} + (E(\Delta)\bv\J) a_{2} = (E(\Delta)\bv)a_{1} + (E(\Delta)\bv)\J a_{2} = (E(\Delta)\bv)a.
\]
Hence, the set function  $ \Delta \mapsto E(\Delta)$ defined in \eqref{RumSti} takes values that are bounded quaternionic linear projections on $V_{R}$. It is immediate that it moreover satisfies \cref{SM1,SM2,SM3,SM4} in \Cref{SpecMeas} because $E_{\I}$ is a spectral measure on $V_{R,\I}$ and hence has the respective properties. Consequently, $E$ is a quaternionic spectral measure. Since $E_{\I}$ commutes with $T$, also $E$ commutes with $T$. From \Cref{ABC}  and the fact that $\sigma(T|_{\ran E_{\I}(\Delta_{\I})})\subset \clos{\Delta_{\I}}$ for $\Delta_{\I}\in\Borel(\cc_{\I})$, we deduce for $T_{\Delta} = T|_{\ran E(\Delta)} = T|_{\ran E_{\I}(\Delta\cap\cc_{\I})}$ that
\[
\sigma_{S}(T_{\Delta}) =\left[ \sigma_{\cc_{\I}}(T_{\Delta}) \right] \subset  \left[ \clos{\Delta\cap\cc_{\I}} \right] =  \clos{\left[\Delta\cap\cc_{\I} \right]} = \clos{\Delta}.
\]
 Therefore $E$ is a spectral resolution for $T$. 
 
Let us now set $V_{0} = \ran E_{\I}(\rr)$ as well as $V_{+} := \ran E_{\I}\left(\cc_{\I}^{+}\right)$ and $V_{-} := \ran E_{\I}\left(\cc_{\I}^{-}\right)$. Then $V_{R,\I} = V_{0} \oplus V_{+}\oplus V_{-}$ is a decomposition of $V_{R}$ into closed $\cc_{\I}$-linear subspaces. The space $V_{0} = \ran E_{\I}(\rr) = \ran E(\rr)$ is even a quaternionic right linear subspace of $V_{R}$ because $E(\rr)$ is a quaternionic right linear operator. Moreover \eqref{KaJuk1} shows that $\bv \mapsto \bv \J$ is a bijection from $V_{+}$ to $V_{-}$. By \Cref{JVSplitThm} the operator
\[
J\bv = E_{\I}\left(\cc_{\I}^{+}\right)\bv\I + E_{\I}\left(\cc_{\I}^{-}\right)\bv(-\I)
\]
 is an imaginary operator on $V_R$. Since $E_{\I}$ commutes with $T$ and $E(\Delta)$ for $\Delta\in\sBorel(\hh)$,  also $J$ commutes with $T$ and $E(\Delta)$. Moreover $\ker J =  V_0 = \ran E(\rr)$ and $\ran J = \ran E_{\I}(\cc_{\I}^{+}) \oplus \ran E_{\I}(\cc_{\I}^{-}) = \ran E(\hh\setminus\rr)$ and hence $(E,J)$ is a spectral system that commutes with $T$. Finally, we have $\sigma\left(T_{+}\right) \subset \cc_{\I}^{\geq}$ for $T_{+} = T|_{V_{+} }= T|_{\ran E_{\I}(\cc_{\I}^{+})}$ and hence  the resolvent of $R_{z}(T_{+})$ exists for any $z\in\cc_{\I}^{-}$. Similarly, the resolvent $R_{z}(T_{-})$ with $T_{-} = T|_{V_{-}} = T|_{\ran E_{\I}(\cc_{\I}^{-})}$ exists for any $z\in\cc_{\I}^{+}$. For $s_0,s_1\in\rr$ with $s_{1}>0$ we can hence set $s_{\I}= s_0 + \I s_1$ and define by
 \[
 R(s_0,s_1) := \left.\left(R_{\overline{s_{\I}}}( T_{+}) E_{+} + R_{s_{\I}}(T_{-})E_{-}\right)\right|_{V_{+}\oplus V_{-}}
 \]
with $E_{+} = E_{\I}(\cc_{\I}^{+})$ and $E_{-} = E_{\I}(\cc_{\I}^{-})$ a bounded operator on $V_{+} \oplus V_{-} = \ran E(\hh\setminus\rr)$. Since $T$ leaves $V_{+}$ and $V_{-}$ invariant,  we then have for $\bv = \bv_{+} + \bv_{-}\in V_{+}\oplus V_{-}$ that
 \begin{align*}
 &R(s_0,s_1) (s_0\id - s_1J  - T) \bv \\
 =& R(s_0,s_1) \left( \bv_{+} s_0 - J \bv_{+} s_1  -T\bv_{+} + \bv_{-} s_0 - J \bv_{-} s_1  -T\bv_{-} \right) \\
 = & R(s_0,s_1) \left(\bv_{+} \overline{s_{\I}} -T\bv_{+} \right) + R(s_0,s_1) \left( \bv_{-} s_{\I}  - T\bv_{-} \right) \\
 =& R_{\overline{s_{\I}}}(T_{+}) \left(\bv_{+} \overline{s_{\I}} -T_{+}\bv_{+} \right) + R_{s_{\I}}(T_{-}) \left( \bv_{-} s_{\I}  - T_{-}\bv_{-} \right) = \bv_{+} + \bv_{-}  = \bv.
 \end{align*}
Similarly we find that
 \begin{align*}
   & (s_0\id - s_1J  - T) R(s_0,s_1)\bv = \\
 =&(s_0\id - s_1J  - T) R_{\overline{s_{\I}}}(T_{+})\bv_{+} + (s_0\id - s_1J  - T) R_{s_{\I}}(T_{-})\bv_{-}\\
= & R_{\overline{s_{\I}}}(T_{+})\bv_{+}s_0 - J(R_{\overline{s_{\I}}}(T_{+})\bv_{+})s_1  - T R_{\overline{s_{\I}}}(T_{+})\bv_{+} \\
&+ R_{s_{\I}}(T_{-})\bv_{-}s_0 - J(R_{s_{\I}}(T_{-})\bv_{-})s_1  - TR_{s_{\I}}(T_{-})\bv_{-}\\
=& R_{\overline{s_{\I}}}(T_{+})\bv_{+}(s_0 - \I s_1)  -  R_{\overline{s_{\I}}}(T_{+})T_{+}\bv_{+} \\
&+ R_{s_{\I}}(T_{-})\bv_{-}(s_0 +\I s_1)  - R_{s_{\I}}(T_{-})T_{-}\bv_{-}\\
= & R_{\overline{s_{\I}}}(T_{+}) \left(\bv_{+}\overline{s} - T_{+}\bv_{+}\right) + R_{s_{\I}}(T_{-}) \left(\bv_{-} s - T_{-}\bv_{-}\right) = \bv_{+} + \bv_{-} = \bv. 
 \end{align*}
 Hence $ R(s_0,s_1)$ is the bounded inverse of $(s_0\id - s_1J  - T)|_{\ran E(\hh\setminus\rr)}$ and so $J$ is actually a spectral orientation for $T$. Consequently, $T$ is a quaternionic spectral operator and the relation \eqref{QWerID} holds true. 
 
\end{proof}
\begin{remark}
We want to stress that \Cref{SOpBackwards} showed a one-to-one relation between quaternionic spectral operators on $V_{R}$ and $\cc_{\I}$-complex spectral operators on $V_{R,\I}$ that are furthermore compatible with the quaternionic scalar multiplication. It did not show a one-to-one relation between quaternionic spectral operators on $V_{R}$ and $\cc_{\I}$-complex spectral operators on $V_{R,\I}$. There exist $\cc_{\I}$-complex spectral operators on $V_{R,\I}$ that are not quaternionic linear and can hence not be quaternionic spectral operators.
\end{remark}

\section{Canonical reduction and intrinsic $S$-functional calculus for quaternionic spectral operators}\label{SpecOpSect2}
As in the complex case any bounded quaternionic spectral operator $T$ can be decomposed into the sum $T = S + N$ of a scalar operator $S$ and a quasi-nilpotent operator $N$. The intrinsic $S$-functional calculus for a spectral operator can then be expressed as a Taylor series similar to the one in \cite{Colombo:2016a} that involves functions of $S$ obtained via spectral integration and powers of $N$. Analogue to the complex case in \cite{Dunford:1958}, the operator $f(T)$ is therefore already determined by the values of $f$ on $\sigma_{S}(T)$ and not only by its values on a neighborhood of $\sigma_{S}(T)$.

\begin{definition}\label{ScalOpDef}
An operator $S\in\boundOP(V_{R})$ is said to be of scalar type if it is a spectral operator and satisfies the identity
\begin{equation}\label{ScalOpEQ}
 S= \int  s\, dE_{J}(s),
\end{equation}
where $(E,J)$ is the spectral decomposition of $S$.
\end{definition}
\begin{remark}\label{ScalarSpecInt}
If we start from a spectral system $(E,J)$ and $S$ is the operator defined by \eqref{ScalOpEQ}, then $S$ is an operator of scalar type and $(E,J)$ is its spectral decomposition. This can easily be checked by direct calculations or indirectly via the following argument: by \Cref{BackCompat}, we can choose $\I\in\SS$ and find
\[
S = \int_{\hh}s \,dE_{J}(s) = \int_{\cc_{\I}} z\,dE_{\I}(z),
\]
where $E_{\I}$ is the spectral measure constructed in \eqref{CSpecMeas}. From the complex theory in \cite{Dunford:1958}, we deduce that $S$ is a spectral operator on $V_{R,\I}$ with spectral decomposition $E_{\I}$ that is furthermore quaternionic linear. By \Cref{SOpBackwards} this is equivalent to $S$ being a quaternionic spectral operator on $V_R$ with spectral decomposition $(E,J)$.

\end{remark}
\begin{lemma}\label{ScalCom}
Let $S$ be an operator of scalar type with spectral decomposition $(E,J)$. An operator $A\in\boundOP(V_R)$ commutes with $S$ if and only if it commutes with the spectral system $(E,J)$.
\end{lemma}
\begin{proof}
If $A\in\boundOP(V_{R})$ commutes with $(E,J)$ then it commutes with $S = \int_{\hh} s\, dE_{J}(s)$ because of \Cref{GenIntProp}. If on the other hand $A$ commutes with $S$, then it also commutes with $E$ by \Cref{Commut}. By \Cref{HomLem} it  commutes in turn with the operator $f(T) = \int_{\hh} f(s) \,dE(s)$  for any $f\in\bsMeas(\hh,\rr)$. If we define 
\[
S_0 := \int_{\hh}\Re(s)\,dE(s)\qquad\text{ and }\qquad S_{1} := \int_{\hh} \underline{s}\, dE_J(S) = J \int_{\hh} |\underline{s}|\,dE(s),
\]
where $\underline{s}=\I_{s}s_1$ denotes the imaginary part of a quaternion $s$, then $AS = S A$ and  $AS_0 = S_0 A$ and in turn
\[
A S_1 = A(S -S_0) = AS - AS_0 =SA -  S_0A = (S-S_0)A = S_1 A.
\]
We can now choose a sequence of pairwise disjoint  sets $\Delta_n\in\sBorel(\hh)$ such that $\sigma_{S}(T) \setminus\rr= \bigcup_{n\in\nn} \Delta_n$ and such that $\dist(\Delta_n,\rr)>0$ for any $n\in\nn$. Then $f_{n}(s)\mapsto |\underline{s}|^{-1}\chi_{\Delta_{n}}(s)$ belongs to $\bsMeas(\hh,\rr)$ for any $n\in\nn$ and in turn
\begin{align*}
AJE(\Delta_n) =& AJ \left(\int_{\hh} |\underline{s}| |\underline{s}|^{-1}\chi_{\Delta_n}(s)\,dE(s)\right)E(\Delta_{n}) \\
= &AJ \left(\int_{\hh} |\underline{s}|\,dE(s)\right) \left(\int_{\hh} |\underline{s}|^{-1}\chi_{\Delta_n}(s)\,dE(s)\right)E(\Delta_{n})\\
= & AS_1 \left(\int_{\hh} |\underline{s}|^{-1}\chi_{\Delta_n}(s)\,dE(s)\right)E(\Delta_{n}) \\
= & S_1 \left(\int_{\hh} |\underline{s}|^{-1}\chi_{\Delta_n}(s)\,dE(s)\right)E(\Delta_{n})A\\
= & J \left(\int_{\hh} |\underline{s}| \,dE(s)\right) \left(\int_{\hh} |\underline{s}|^{-1}\chi_{\Delta_n}(s)\,dE(s)\right)E(\Delta_{n})A\\
= & J \left(\int_{\hh} |\underline{s}| |\underline{s}|^{-1}\chi_{\Delta_n}(s)\,dE(s)\right)E(\Delta_{n})A = JE(\Delta_{n})A.
\end{align*}
Since $\sigma_{S}(S)\setminus\rr \subset \bigcup_{n\in\nn} \Delta_{n}$, we have $\sum_{n=0}^{+\infty} E(\Delta_n)\bv = E(\sigma_{S}(T)\setminus\rr)\bv = {E(\hh\setminus\rr) \bv}$ for all $\bv\in V_{R}$ by \Cref{Support}. As $J = J E(\hh\setminus\rr)$, we hence find
\begin{align*}
AJ\bv  = &AJ E(\hh\setminus\rr)\bv= \sum_{n=1}^{+\infty} AJE(\Delta_n)\bv \\
=& \sum_{n=1}^{+\infty}JE(\Delta_{n})A\bv = JE(\hh\setminus\rr)A\bv = JA\bv,
\end{align*}
which finishes the proof. 

\end{proof}

\begin{definition}
An operator $N\in\boundOP(V_{R})$ is called quasi-nilpotent if 
\begin{equation}\label{QNCond}
\lim_{n\to\infty}\|N^n\|^{\frac1n} = 0.
\end{equation}
\end{definition}
The following corollaries are immediate consequences of Gelfands formula
\[
r(T) = \lim_{n\to+\infty}\left\| T^n\right\|^{\frac{1}{n}},
\]
 for the spectral radius $r(T) = \max_{s\in\sigma_{S}(T)}|s|$ of $T$ and of \cref{Funf} of \Cref{NewSCalcProp}.
\begin{corollary}\label{NilSpec}
An operator $N\in\boundOP(V_{R})$ is quasi-nilpotent if and only if $\sigma_S(T) = \{0\}$.
\end{corollary}
\begin{corollary}\label{SumNilSpec}
Let $S,N\in\boundOP(V_{R})$ be commuting operators and let $N$ be quasi-nilpotent. Then $\sigma_S(S+N) = \sigma_S(S)$.
\end{corollary}
We are now ready to show the main result of this section: the canonical reduction of a spectral operator, the quaternionic analogue of Theorem~5 in \cite[Chapter XV.4.3]{Dunford:1958}.
\begin{theorem}\label{DecompTHM}
An operator $T\in\boundOP(V_{R})$ is a spectral operator if and only if it is the sum $T = S+N$ of a bounded operator $S$ of scalar type and a quasi-nilpotent operator $N$ that commutes with $S$.  Furthermore, this decomposition is unique and $T$ and $S$ have the same $S$-spectrum and the same spectral decomposition $(E, J)$.
\end{theorem}
\begin{proof}
Let us first show that any operator $T\in\boundOP(V_{R})$ that is the sum $T=S+N$ of an operator $S$ of scalar type and a quasi-nilpotent operator $N$ commuting with $S$ is a spectral operator. If $(E,J)$ is the spectral decomposition of $S$, then \Cref{ScalCom} implies $E(\Delta)N = NE(\Delta)$ for all $\Delta\in\sBorel(\hh)$ and $J N = N J$.  Since $T = S+ N$, we find that also $T$ commutes with $(E,J)$.

Let now $\Delta\in\sBorel(\hh)$. Then $T_{\Delta}  = S_{\Delta} + N_{\Delta}$, where as usual the subscript $\Delta$ denotes the restriction of an operator to $V_{\Delta} = E(\Delta)V_{R}$. Since $N_{\Delta}$  inherits the property of being quasi-nilpotent from $N$ and commutes with $S_{\Delta}$, we deduce from \Cref{SumNilSpec} that
\[\sigma_{S}(T_{\Delta}) = \sigma_{S}(S_{\Delta} + N_{\Delta}) = \sigma_{S}(S_{\Delta}) \subset \clos{\Delta}.\]
Thus $(E,J)$ satisfies \cref{SOi} and \cref{SOii} of \Cref{SpecOP}. It remains to show that also \cref{SOiii} holds true. Therefore let $V_{0} = \ran E(\hh\setminus\rr)$ and set $T_{0} = T|_{V_0}$,  $S_{0}=S|_{V_0}$, $N_0 = N|_{V_{0}}$ and $J_0 = J|_{V_0}$ and choose  $s_0,s_1\in\rr$ with $s_1>0$. Since $(E,J)$ is the spectral resolution of $S$,  the operator $s_{0} \id_{V_0} - s_1 J_0 - S_0$ has a bounded inverse  $R(s_0,s_1) = (s_{0} \id_{V_0} - s_1 J_0 - S_0)^{-1} \in \boundOP(V_0)$. The operator $N_{0}$ is quasi-nilpotent because $N$ is quasi-nilpotent and hence it satisfies \eqref{QNCond}. The root test thus shows the convergence of the series $\sum_{n=0}^{+\infty} N_{0}^{n}R(s_0,s_1)^{n+1}$ in $\boundOP(V_0)$. Since $T_0$, $N_0$, $S_0$ and $J_0$ commute mutually, we have
\begin{align*}
&(s_0\id_{V_0} - s_1J_0 - T_{0}) \sum_{n=0}^{+\infty} N_0^n R(s_0,s_1)^{n+1} \\
=& \sum_{n=0}^{+\infty}N_0^n R(s_0,s_1)^{n+1}(s_0\id_{V_0} - s_1J_0 - S_{0}-N_{0})\\
= &\sum_{n=0}^{+\infty}N_0^n R(s_0,s_1)^{n+1}(s_0\id_{V_0} - s_1J_0 - S_{0}) - \sum_{n=0}^{+\infty}N_0^n R(s_0,s_1)^{n+1}N_0 \\
=& \sum_{n=0}^{+\infty}N_0^n R(s_0,s_1)^{n} - \sum_{n=0}^{+\infty}N_0^{n+1} R(s_0,s_1)^{n+1} = \id_{V_0}.
\end{align*}
We find that $s_0\id_0 - s_1J_0 - T_{0}$ has a bounded inverse for $s_0,s_1\in\rr$ with $s_1>0$ such that $J$ is a spectral orientation for $T$. Hence, $T$ is a spectral operator and $T$ and $S$ have the same spectral decomposition  $(E,J)$. 

Since the spectral decomposition of $T$ is uniquely determined, $S= \int_{\hh} s\, dE_{J}(s)$ and in turn also $N=T-S$ are uniquely determined. Moreover, \Cref{SumNilSpec} implies that $\sigma_{S}(T) = \sigma_{S}(S)$.

Now assume that $T$ is a spectral operator and let $(E,J)$ be its spectral decomposition. We set
\[S := \int_{\hh} s\, dE_{J}(s)\qquad \text{and}\qquad  N:= T-S.\]
By \Cref{ScalarSpecInt} the operator $S$ is of scalar type and its spectral decomposition is $(E,J)$. Since $T$ commutes with $(E,J)$, it commutes with $S$ by \Cref{ScalCom}. Consequently,  $N = T - S$ also commutes with $S$ and with $T$. What remains to show is that $N$ is quasi-nilpotent. In view of Corollary~\ref{NilSpec}, it is sufficient to show that $\sigma_S(N)$ is for any $\varepsilon>0$ contained in the open ball $B_{\varepsilon}(0)$ of radius $\varepsilon$ centered at $0$ .

For arbitrary  $\varepsilon >0$, we choose $\alpha>0$ such that $0<(1+ C_{E,J})\alpha < \varepsilon$, where $C_{E,J}>0$ is the constant in~\eqref{NormEst2}. We decompose $\sigma_S(T)$ into the union of disjoint axially symmetric Borel sets $\Delta_1,\ldots,\Delta_n\in\sBorel(\hh)$ such that for each $\ell\in\{1,\ldots,n\}$ the set $\Delta_\ell$ is contained in a closed axially symmetric set, whose intersection with any complex halfplane is a half-disk of diameter $\alpha$. More precisely, we assume that there exist points $s_{1},\ldots, s_{n}\in\hh$ such that for all $\ell = 1,\ldots,n$
\[
\Delta_\ell \subset B_{\alpha}^+([s_{\ell}]) = \{ p\in\hh: \dist(p,[s_{\ell}]) \leq \alpha \ \text{and}\  p_1 \geq s_{\ell,1}\}.
\]
Observe that we have either $s_\ell \in\rr$ or $B_{\alpha}^+([s_{\ell}]) \cap\rr=\emptyset$. 

We set $V_{\Delta_\ell} = E(\Delta_\ell)V_{R}$. As $T$ and $S$ commute with $E(\Delta_{\ell})$, also $N=T-S$ does and so $N V_{\Delta_\ell} \subset V_{\Delta_\ell}$. Hence, $N_{\Delta_\ell} = N|_{V_{\Delta_\ell}}\in\boundOP(V_{\Delta_\ell})$. If $s$ belongs to $\rho_S(N_{\Delta_\ell})$ for all $\ell\in\{1,\ldots,n\}$, we can set 
\[\Q(s)^{-1} := \sum_{\ell=1}^n \Q_{s}(N_{\Delta_{\ell}})^{-1}E(\Delta_{\ell}),\]
 where 
 \[
 \Q_{s}(N_{\Delta_{\ell}})^{-1} =\big( N_{\Delta_\ell}^2 - 2s_0 N_{\Delta_{\ell}}+|s|^2\id_{V_{\Delta_{\ell}}}\big)^{-1}\in\boundOP(V_{\Delta_{\ell}})
 \]
 is the pseudoresolvent of $N_{\Delta_{\ell}}$ as $s$. The operator $\Q(s)^{-1}$ commutes with $E(\Delta_{\ell})$ for any $\ell\in\{1,\ldots,n\}$ such that
\begin{gather*}(N^2 - 2s_0 N + |s|^2\id_{V_{R}}) \Q(s)^{-1} \\
= \sum_{\ell=1}^n (N_{\Delta_{\ell}}^2 - 2s_0 N_{\Delta_{\ell}} + |s|^2\id_{V_{\Delta_{\ell}}})\Q_{s}(N_{\Delta_{\ell}})^{-1}E(\Delta_\ell) =\sum_{\ell=1}^nE(\Delta_\ell)=\id_{V_R}
\end{gather*}
and
\begin{align*} 
&\Q(s)^{-1}(N^2 - 2s_0 N + |s|^2\id_{V_{R}}) =\\
=& \sum_{\ell=1}^{n}\Q_{s}(N_{\Delta_{\ell}})^{-1} E(\Delta_{\ell})(N^2 - 2s_0 N + |s|^2\id_{V_R}) \\
= & \sum_{\ell=1}^{n}\Q_{s}(N_{\Delta_{\ell}})^{-1}(N_{\Delta_{\ell}}^2 - 2s_0 N_{\Delta_{\ell}} + |s|^2\id_{V_{\Delta_{\ell}}})E(\Delta_{\ell})\\
= & \sum_{\ell=1}^n E(\Delta_{\ell}) = \id_{V_{R}}.
\end{align*}
Therefore $s\in\rho_S(N)$ such that $  \bigcap_{\ell = 1}^n \rho_S(N_{\Delta_\ell})  \subset \rho_S(N)$ and in turn $ \sigma_S(N)\subset \bigcup_{\ell=1}^n\sigma_S(N_{\Delta_{\ell}})$. It is hence sufficient to show that $\sigma_S(N_{\Delta_\ell})\subset B_{\varepsilon}(0)$ for all $\ell = 1,\ldots,n$.

We distinguish two cases: if $s_{\ell}\in \rr$, then we write
\[N_{\Delta_{\ell}} = (T_{\Delta_{\ell}}- s_{\ell}\id_{V_{\Delta_{\ell}}}) + (s_{\ell}\id_{V_{\Delta_{\ell}}} - S_{\Delta_{\ell}}). \]
Since $s_{\ell}\in\rr$, we have for $p\in\hh$ that
\begin{align*}
 &\Q_{p}(T_{\Delta_{\ell}} - s_{\ell}\id_{V_{\Delta_{\ell}}})=\\
  = &(T_{\Delta_{\ell}}^2 - 2s_{\ell}T_{\Delta_{\ell}} + s_{\ell}^2 \id_{V_{\Delta_{\ell}}} - 2p_0(T_{\Delta_{\ell}} - s_{\ell}\id_{V_{\Delta_{\ell}}}) + (p_0^2 + p_1^2)\id_{V_{\Delta_{\ell}}}\\
 = &T_{\Delta_{\ell}}^2 - 2 (p_0 - s_\ell)T_{\Delta_{\ell}} + \left((p_0 - s_{\ell})^2 + p_1^2\right)\id_{V_{\Delta_{\ell}}} = \Q_{p-s_{\ell}}(T_{\Delta_{\ell}})
 \end{align*}
 and thus
 \begin{equation}\label{KLaKuz}
 \begin{split}
 \sigma_{S}(T_{\Delta_{\ell}}-s_{\ell}\id_{V_{\ell}}) = &\{p - s_{\ell} \in\hh: p \in \sigma_{S}(T_{\Delta_{\ell}})\}\\
  \subset&   \{p - s_{\ell} \in\hh: p \in B_{\alpha}^+(s_\ell)\} = B_{\alpha}(0).
 \end{split}
  \end{equation}
 Moreover, the function $f(s) = (s_{\ell} - s)\chi_{\Delta_{\ell}}(s)$ is an intrinsic slice function because $s_{\ell}\in\rr$. As it is bounded, its integral with respect to $(E,J)$ is defined and
 \[
  s_{\ell}\id_{V_{\Delta_{\ell}}} - S_{\Delta_{\ell}} = \left.\left( \int_{\hh} (s_{\ell} - s  )\chi_{\Delta_{\ell}}(s) \, dE_{J}(s)\right)\right|_{V_{\Delta_{\ell}}}.
 \]
 We thus have
 \begin{equation}\label{JA2} 
 \| s_{\ell} \id_{V_{\Delta_{\ell}}}- S_{\Delta_{\ell}}\| \leq C_{E,J} \|  (s_{\ell} - s) \chi_{\Delta_{\ell}}(s)\|_{\infty} \leq C_{E,J} \alpha
 \end{equation}
 because $\Delta_{\ell}\subset B_{\alpha}([s_{\ell}])^+ = \clos{B_{\alpha}(s_{\ell})}$. Since the operators $T_{\Delta_{\ell}}- s_{\ell}\id_{V_{\Delta_{\ell}}}$ and $s_{\ell}\id_{V_{\Delta_{\ell}}} - S_{\Delta_{\ell}}$ commute, we conclude from \Cref{Funf} in \Cref{NewSCalcProp} together with \eqref{KLaKuz} and \eqref{JA2} that
 \begin{gather*}
 \sigma_{S}(T_{\Delta{\ell}}) =  \sigma_{S}\left((T_{\Delta_{\ell}}- s_{\ell}\id_{V_{\Delta_{\ell}}}) + (s_{\ell}\id_{V_{\Delta_{\ell}}} - S_{\Delta_{\ell}})\right)\\
  \subset \left\{s\in\hh: \dist\left(s,\sigma_{S}\big(T_{\Delta_{\ell}}- s_{\ell}\id_{V_{\Delta_{\ell}}}\big)\right) \leq C_{E,J}\alpha\right\} \subset B_{\alpha(1+C_{E,J})}(0) \subset B_{\varepsilon}(0). 
 \end{gather*}
 
 If $s_{\ell}\notin\rr$, then let us write
\begin{equation}\label{AA!11}
N_{\Delta_{\ell}} = (T_{\Delta_{\ell}}- s_{\ell}\id_{V_{\Delta_{\ell}}} - s_{\ell,1}J_{\Delta_{\ell}}) + (s_{\ell}\id_{V_{\Delta_{\ell}}} + s_{\ell,1}J_{\Delta_{\ell}} - S_{\Delta_{\ell}}) 
\end{equation} 
with $J_{\Delta_{\ell} }= J|_{V_{\Delta_{\ell}}}$. Since  $E(\Delta_{\ell})$ and $J$ commute, $J_{\Delta_{\ell}}$ is an imaginary operator on $V_{\Delta_{\ell}}$ and it moreover commutes with $T_{\Delta_{\ell}}$.  Since $-J_{\Delta_{\ell}}^2 = -J^2 |_{V_{\Delta}} = E(\hh\setminus\rr)|_{V_{\Delta_{\ell}}} = \id_{V_{\Delta_{\ell}}}$ as $\Delta_{\ell}\subset\hh\setminus\rr$, we find for $s = s_0  + \I_{s}s_1\in\hh$ with $s_1\geq 0$ that
 \begin{equation}\label{NeoN}
\begin{split}
& \big(s_0\id_{V_{\Delta_{\ell}}} + s_1 J_{\Delta_{\ell}} - T_{{\Delta}_{\ell}}\big) \big(s_0\id_{V_{\Delta_{\ell}}} - s_1 J_{\Delta_{\ell}} - T_{{\Delta}_{\ell}}\big) =\\
=& s_0^2 - s_1^2J_{\Delta_{\ell}}^2 - 2 s_0 T_{\Delta_{\ell}} + T_{\Delta_{\ell}}^2  =\Q_{s}(T_{\Delta_{\ell}}). 
\end{split}
 \end{equation}
Because of condition \cref{SOiii} in \Cref{SpecOP}, the operator $(s_0\id - s_1 J - T)|_{\ran E(\hh\setminus\rr)}$ is invertible if $s_1>0$. Since this operator commutes with $E(\Delta_{\ell})$, the restriction of its inverse to $V_{\Delta_{\ell}}$ is the inverse of $(s_0\id_{V_{\Delta_{\ell}}} - s_1 J_{\Delta_{\ell}} - T_{{\Delta}_{\ell}})  $ in $\boundOP(V_{\Delta_{\ell}})$. Hence, if $s_1>0$, then $(s_0\id_{V_{\Delta_{\ell}}} - s_1 J_{\Delta_{\ell}} - T_{{\Delta}_{\ell}})^{-1}\in\boundOP(V_{\Delta_{\ell}})$ and we conclude from \eqref{NeoN} that
\begin{equation}\label{HuJi}
\big(s_0\id_{V_{\Delta_{\ell}}} + s_1 J_{\Delta_{\ell}} - T_{{\Delta}_{\ell}}\big)^{-1} \in\boundOP(V_{\Delta_{\ell}})\quad\Longleftrightarrow\quad \Q_{s}(T_{\Delta_{\ell}})^{-1}\in\boundOP(V_{\Delta_{\ell}}).
\end{equation}
 If on the other hand $s_1 =0$, then both factors on the left-hand side of \eqref{NeoN} agree and so \eqref{HuJi} holds true also in this case.
Hence,  $s\in\rho_S(T_{\Delta_{\ell}})$ if and only if $(s_0\id_{V_{\Delta_{\ell}}} + s_1J_{\Delta_{\ell}} - T)$ has an inverse in $\boundOP(V_{\Delta_{\ell}})$. Since 
\[
\sigma_{S}(T_{\Delta_{\ell}}) \subset \overline{\Delta_{\ell}} \subset B_{\alpha}^{+}([s_{\ell}]) \subset \{s = s_0 + \I_{s}s_{1}\in\hh: s_1\geq s_{\ell,1}\},
\]
the operator $s_0\id_{V_{\Delta_{\ell}}} + s_1J_{\Delta_{\ell}} - T_{\Delta_{\ell}}$ is in particular invertible for any $s\in\hh$ with $0\leq s_1 < s_{\ell,1}$. As $J_{\Delta_{\ell}}$ is a spectral orientation for $T_{\Delta_{\ell}}$, this operator is also invertible if $s_1 <0$ and hence we even find
\begin{equation}\label{SLAQ}
(s_0\id_{V_{\Delta_{\ell}}} + s_1J_{\Delta_{\ell}} - T_{\Delta_{\ell}})^{-1}\in\boundOP(V_{\Delta_{\ell}})\quad \forall s_0,s_1\in\rr:  s_1 < s_{\ell,1}.
\end{equation}

We can use these observations to deduce a spectral mapping property: a straight forward computation using the facts that $T_{\Delta_{\ell}}$ and $J_{\Delta_{\ell}}$ commute and that $J_{\Delta_{\ell}}^2 = -\id_{V_{\Delta_{\ell}}}$ shows
 \begin{equation}
 \begin{split}\label{KKAH}
  &\Q_{s}(T_{\Delta_{\ell}} - s_{\ell,0}\id_{V_{\Delta_{\ell}}}- s_{\ell,1}J_{\Delta_{\ell}})\\
   =& \left((s_0+s_{\ell,0})\id_{V_{\Delta_{\ell}}} + (s_1 + s_{\ell,1})J_{\Delta_{\ell}} - T_{{\Delta}_{\ell}}\right)\\
   &\cdot \left((s_0+s_{\ell,0})\id_{V_{\Delta_{\ell}}} + ( s_{\ell,1}-s_1) J_{\Delta_{\ell}} - T_{{\Delta}_{\ell}}\right).
\end{split}
\end{equation}
If $s_{1}>0$ then the second factor is invertible because of \eqref{SLAQ}. Hence $s\in\rho_S(T_{\Delta_{\ell}} - s_{\ell,0}\id_{V_{\Delta_{\ell}}}- s_{\ell,1}J_{\Delta_{\ell}})$ if and only if the first factor in \eqref{KKAH} is also invertible, i.e. if and only if
\begin{equation}\label{Zumpfl}
 \big((s_0+s_{\ell,0})\id_{V_{\Delta_{\ell}}} + (s_1 + s_{\ell,1})J_{\Delta_{\ell}} - T_{{\Delta}_{\ell}}\big)^{-1}\in\boundOP(V_{\Delta_{\ell}})
 \end{equation}
 exists. If on the other hand $s_{1} = 0$, then both factors in \eqref{KKAH} agree. Hence, also in this case, $s$ belongs to $\rho_S(T_{\Delta_{\ell}} - s_{\ell,0}\id_{V_{\Delta_{\ell}}}- s_{\ell,1}J_{\Delta_{\ell}})$ if and only if the operator in \eqref{Zumpfl} exists. By \eqref{HuJi}, the existence of \eqref{Zumpfl} is however equivalent to $s_0  + s_{\ell,0} + (s_1 + s_{\ell,1})\SS\subset \rho_{S}(T_{\Delta})$ so that
 \[
 \rho_{S}(T_{\Delta_{\ell}} - s_{\ell,0}\id_{V_{\Delta_{\ell}}}- s_{\ell,1}J_{\Delta_{\ell}}) = \{ s\in\hh: s_0 + s_{\ell,0} + (s_{1}+s_{\ell,1})\I_{s}\in \rho_{S}(T_{\Delta_{\ell}})\}
 \]
and  in turn
 \begin{gather*}
 \sigma_{S}(T_{\Delta_{\ell}} - s_{\ell,0}\id_{V_{\Delta_{\ell}}}- s_{\ell,1}J_{\Delta_{\ell}}) = \{ s\in\hh: s_0 + s_{\ell,1} + (s_{1}+s_{\ell,1})\I_{s}\in \sigma_{S}(T_{\Delta_{\ell}})\}\\
 \subset \{ s\in\hh: s_0 + s_{\ell,0} + (s_{1}+s_{\ell,1})\I_{s}\in B_{\alpha}^{+}(s_{\ell})\} = B_{\alpha}(0).
 \end{gather*}
   
For the second operator in \eqref{AA!11}, we have again 
\[
s_{\ell}\id_{V_{\Delta_{\ell}}} + s_{\ell,1}J_{\Delta_{\ell}} - S_{\Delta_{\ell}} = \left. \left(\int_{\hh} (s_{\ell,0} +i_{s}s_{\ell,1} - s)\chi_{\Delta_{\ell}}(s)\,dE_{J}(s)\right)\right|_{V_{\Delta_{\ell}}}
\]
and so
\[
\| s_{\ell}\id_{V_{\Delta_{\ell}}} + s_{\ell,1}J_{\Delta_{\ell}} - S_{\Delta_{\ell}}\| \leq C_{E,J} \|(s_{\ell,0} +i_{s}s_{\ell,1} - s)\chi_{\Delta_{\ell}}(s)\|_{\infty} \leq C_{E,J} \alpha.
\]
 Since the operators $T_{\Delta_{\ell}}- s_{\ell}\id_{V_{\Delta_{\ell}}} - s_{\ell,1}J_{\Delta_{\ell}}$ and $s_{\ell}\id_{V_{\Delta_{\ell}}} + s_{\ell,1}J_{\Delta_{\ell}} - S_{\Delta_{\ell}}$ commute, we conclude as before from \Cref{Funf} in \Cref{NewSCalcProp} that $\sigma_{S}(T_{\Delta_{\ell}})\subset B_{\alpha(1+C_{E,J})}(0) = B_{\varepsilon}(0)$. 
 
 Altogether, we obtain that $N$ is quasi-nilpotent, which concludes the proof.
 
 \end{proof}

\begin{definition}
Let $T\in\boundOP(V_{R})$ be a spectral operator and decompose $T=S+N$ as in \Cref{DecompTHM}. The scalar operator $S$ is called the scalar part of $T$ and the quasi-nilpotent operator $N$ is called the radical part of $T$. 
\end{definition}
\begin{remark}
Let $T\in\boundOP(V_{R})$ be a spectral operator. The canonical decomposition of $T$ into its scalar part and its radical part obviously coincides for any $\I\in\SS$ with the canonical decomposition of $T$ as a $\cc_{\I}$-linear spectral operator on $V_{\I}$.
\end{remark}
The remainder of this section discusses the $S$-functional calculus for spectral operators. Similar to the complex case, one can express $f(T)$ for any intrinsic function $f$ as a formal Taylor series in the radical part $N$ of $T$. The Taylor coefficients are spectral integrals of $f$ with respect to the spectral decomposition of $T$. Hence these coefficients, and in turn also $f(T)$, do only depend on the values of $f$ on the $S$-spectrum $\sigma_{S}(T)$ of $T$ and not on the values of $f$ on an entire neighborhood of $\sigma_{S}(T)$. The operator $f(T)$ is again a spectral operator and its spectral decomposition can easily be constructed from the spectral decomposition of $T$. 
\begin{proposition}\label{SSIDProp}
Let $S\in\boundOP(V_R)$ be an operator of scalar type. If $f\in\intrin(\sigma_{S}(S))$, then
\begin{equation}\label{SSID}
f(S) = \int_{\hh} f(s)\, dE_{J}(s),
\end{equation}
where $f(S)$ is intended in the sense of the $S$-functional calculus introduced in \Cref{IntCalcSect}.
\end{proposition}
\begin{proof}
 Since $1(T) = \id = \int_{\hh} 1\, dE_{J}(s)$ and $s(S) = S = \int_{\hh} s\,dE_J(s)$, the product rule and the $\rr$-linearity of both the $S$-functional calculus and the spectral integration imply that \eqref{SSID} holds true for any intrinsic polynomial. It in turn also holds true for any intrinsic rational function in $\intrin(\sigma_{S}(S))$, i.e. for any function $r$ of the form $r(s) = p(s)q(s)^{-1}$ with intrinsic polynomials $p$ and $q$ such that $q(s) \neq 0$ for any $s\in\sigma_{S}(S)$.

Let now $f\in\intrin(\sigma_S(S))$ be arbitrary and let $U$ be a bounded axially symmetric open set such that $\sigma_S(T)\subset U$ and $\clos{U}\subset \fdom(f)$. Runge's theorem for slice hyperholomorphic functions in \cite{Colombo:2011c} implies the existence of a sequence of intrinsic rational functions $r_{n}\in\intrin\big(\clos{U}\big)$ such that $r_{n}\to f$ uniformly on $\clos{U}$. Because of \Cref{GenIntProp}, we thus have
\[
 \int_{\hh} f(s)\, dE_{J}(s)\  = \lim_{n\to+\infty} \int_{\hh} r_n(s)\,dE_{J}(s) = \lim_{n\to+\infty} r_n(S) = f(S).
\]

\end{proof}
\begin{theorem}\label{IntCalcSer}
Let $T\in\boundOP(V_R)$ be a spectral operator with spectral decomposition $(E,J)$ and let $T=S+N$ be the decomposition of $T$ into scalar and radical part. If $f\in\intrin(\sigma_{S}(T))$, then
\begin{equation}\label{JUJI}
 f(T) = \sum_{n=0}^{+\infty} N^n\frac{1}{n!}  \int _{\hh}(\sderiv^n f)(s)\,dE_{J}(s),
 \end{equation}
where $f(T)$ is intended in the sense of the the $S$-functional calculus in \Cref{IntCalcSect} and the series converges in the operator norm.  
\end{theorem}
\begin{proof}
Since $T=S+N$ with $SN=NS$ and $\sigma_{S}(N)= \{0\}$, it follows from \cref{Funf} in \Cref{NewSCalcProp} that 
\[
 f(T) = \sum_{n=0}^{+\infty}  N^n \frac{1}{n!}\left(\sderiv^nf\right)(S).
 \]
 What remains to show is that 
 \begin{equation}
(\sderiv^nf)(S) = \int_{\hh} (\sderiv^nf)(s)\, dE_{J}(s),
\end{equation}
but this follows immediately from \Cref{SSIDProp}.

\end{proof}

The operator $f(T)$ is again a spectral operator and its radical part can be easily obtained from the above series expansion.
\begin{definition}
A spectral operator $T\in\boundOP(V_R)$ is called of type $m\in\nn$ if and only if its radical part satisfies~$N^{m+1} = 0$.
\end{definition}

\begin{lemma}\label{REMThm}
A spectral operator $T\in\boundOP(V_R)$ with spectral resolution $(E,J)$ and radical part $N$ is of type $m$ if and only if
\begin{equation}\label{REM}
 f(T) = \sum_{n=0}^{m} N^n\frac{1}{n!}  \int_{\hh} (\sderiv^n f)(s)\,dE_{J}(s) \quad \forall f\in\intrin(\sigma_S(T)).
\end{equation}
In particular $T$ is a scalar operator if and only if it is of type $0$.
\end{lemma}
\begin{proof}
 If $T$ is of type $m$ then the above formula follows immediately from \Cref{IntCalcSer} and $N^{m+1} = 0$. If on the other hand \eqref{REM} holds true, then we choose $f(s) = \frac{1}{m!} s^m$ in \eqref{JUJI} and \eqref{REM} and subtract these two expressions. We obtain
 \[ 0 = N^{m+1} \int_{\hh}\, dE_{J}(s) = N^{m+1}.\]
 
\end{proof}

\begin{theorem}\label{Kakao}
Let $T\in\boundOP(V_R)$ be a spectral operator with spectral decomposition $(E,J)$. If $f\in\intrin(\sigma_S(T))$, then $f(T)$ is a spectral operator and the spectral decomposition $(\tilde{E},\tilde{J})$ of $f(T)$ is given by
\[\tilde{E}(\Delta) = E\big(f^{-1}(\Delta)\big)\quad\forall \Delta\in\sBorel(\hh)\qquad \text{and}\qquad \tilde{J} = \int_{\hh} \I_{f(s)}\, dE_{J}(s),\]
where $\I_{f(s)} = 0$ if $f(s)\in\rr$ and $\I_{f(s)} = \underline{f(s)}/|\underline{f(s)}|$ if $f(s)\in\hh\setminus\rr$. 
For any $g\in\bsIntrin(\hh)$ we have
\begin{equation}\label{SEAL1}
\int_{\hh} g(s)\,d\tilde{E}_{\tilde{J}}(s) = \int_{\hh}(g\circ f)(s)\,dE_{J}(s)
\end{equation}
and  if $S$ is the scalar part of $T$, then $f(S)$ is the scalar part of $f(T)$.
\end{theorem}
\begin{proof}
We first show that $f(S)$ is a scalar operator with spectral decomposition $(\tilde{E}, \tilde{J})$. By \Cref{SymMeas} the function $f$ is $\sBorel(\hh)$-$\sBorel(\hh)$-measurable, such that $\tilde{E}$ is a well-defined spectral measure on $\sBorel(\hh)$.

The operator $\tilde{J}$ obviously commutes with $E$. Moreover, writing $f(s) = \alpha(s) + \I_s \beta(s)$ as in \Cref{IntStruct} we have  $\I_{f(s)} = \I_{s} \sgn(\beta(s))$. If we set $\Delta_{+} = \{ s\in\hh: \beta(s)  > 0\}$, $\Delta_{-} = \{s\in\hh:\beta(s) < 0 \}$ and $\Delta_{0} = \{s\in\hh:\beta(s) = 0\}$, we therefore have
\[
\tilde{J} = J E(\Delta_{+}) - J E( \Delta_{-}).
\]
As $\beta(s) = 0$ for any $s\in\rr$, we have $\rr \subset \Delta_0$ and hence $V_{+} = \ran E(\Delta_{+} )\subset \ran E(\hh\setminus\rr) = \ran J$ and similarly also $V_{-} = \ran E(\Delta_{-})\subset \ran J$. Since $J$ and $E$ commute, $V_{+}$ and $V_{-}$ are invariant subspaces of $J$ contained in $\ran J$ such that $J_{+}$ and $J_{-}$ define bounded surjective operators on $V_{+}$ resp. $V_{-}$. Moreover $\ker {J} = \ran E(\rr)$ and hence $\ker J|_{V_{+}} = 
V_{+}\cap \ker J = \{\bO\}$ and $\ker J|_{V_{-}} = V_{-}\cap \ker J = \{\bO\}$, 
such that $\ker \tilde{J} = \ran E(\Delta_0)$ and $\ran \tilde{J} = \ran E(\Delta_{+}) \oplus \ran E(\Delta_{-}) = \ran E(\Delta_{+}\cup\Delta_{-})$. 

Now observe that $f(s) \in \rr$ if and only if $\beta(s) = 0$. Hence, $f^{-1}(\rr) = \Delta_0$ and $f^{-1}(\hh\setminus\rr) = \Delta_{+}\cup\Delta_{-}$ and we find
\[
\ran  \tilde{J} = \ran E(\Delta_{+}\cup\Delta_{-}) = \ran E\left(f^{-1}(\hh\setminus\rr)\right) = \ran \tilde{E}(\hh\setminus\rr)
\]
and
\[
\ker\tilde{J} = \ran E(\Delta_0) = \ran E\left(f^{-1}(\rr)\right) = \ran \tilde{E}(\rr).
\]
Moreover, as $E(\Delta_+)E(\Delta_{-}) = E(\Delta_{-}) E(\Delta_{+}) = 0$ and $-J^2 = E(\hh\setminus\rr)$, we have
\begin{align*}
 - \tilde{J}^2 =& - J^2 E(\Delta_{+})^2  - (-J^2)E(\Delta_{-})^2 \\
 =& E(\hh\setminus\rr)E(\Delta_{+}) + E(\hh\setminus\rr)E(\Delta_{-}) \\
 = & E(\Delta_{+}\cup\Delta_{-}) = \tilde{E}(\hh\setminus\rr),
\end{align*}
where we used that $\Delta_{+}\subset \hh\setminus\rr$ and $\Delta_{-}\subset\hh\setminus\rr$ as $\rr\subset\Delta_{0}$. Hence $-\tilde{J}^2$ is the projection onto $\ran \tilde{J}$  along $\ker\tilde{J}$ and so  $\tilde{J}$ is actually an imaginary operator and $(\tilde{E},\tilde{J})$ in turn is a spectral system.

Let $g=\sum_{\ell=0}^na_{\ell}\chi_{\Delta_{\ell}}\in\bsMeas(\hh,\rr)$ be a simple function. Then $(g\circ f)(s)  = \sum_{\ell=0}^{n} a_{\ell} \chi_{f^{-1}(\Delta_{\ell})}(s)$ is also a simple function in $\bsMeas(\hh,\rr)$ and 
\begin{align*}
\int_{\hh}g(s)\,d\tilde{E}(s) = & \sum_{\ell=0}^na_{\ell}\tilde{E}(\Delta_{\ell}) = \sum_{\ell=0}^na_{\ell}E\big(f^{-1}(\Delta_{\ell})\big) = \int_{\hh}(g\circ f)(s) \,dE(s).
\end{align*}
Due to the density of simple functions in $(\bsMeas(\hh,\rr),\|.\|_{\infty})$, we hence find 
\begin{align*}
\int_{\hh}g(s)\,d\tilde{E}(s) = \int_{\hh}(g\circ f)(s) \,dE(s),\qquad\forall g\in\bsMeas(\hh,\rr).
\end{align*}
If $g\in \bsIntrin(\hh)$ then we deduce from \Cref{IntStruct} that $g(s) = \gamma(s) + \I_{s} \delta(s)$ with $\gamma,\delta\in\bsMeas(\hh,\rr)$ and $\I_{s} = \underline{s}/|\underline{s}|$ if $s\notin\rr$ and $\I_{s} = \delta(s) = 0$ if $s\in\rr$. We then have $(g\circ f)(s) = \gamma(f(s)) + \I_{f(s)} \delta(f(s))$ and find 
\begin{align*}
\int_{\hh}g(s)\,d\tilde{E}_{\tilde{J}}(s) =& \int_{\hh}\gamma(s)\,d\tilde{E}(s) + \tilde{J} \int_{\hh} \delta(s)\,d\tilde{E}(s)\\
=& \int_{\hh}(\gamma\circ f)(s)\,dE(s) + \tilde{J} \int_{\hh} (\delta\circ f)(s)\,dE(s)\\
= & \int_{\hh}(\gamma\circ f)(s)\,dE_J(s) + \int_{\hh} \I_{f(s)}\,dE_J(s)\int_{\hh}(\delta\circ f)(s)\,dE(s)\\
= & \int_{\hh}(\gamma\circ f)(s) + \I_{f(s)} (\delta\circ f)(s) \,dE_J(s) = \int_{\hh} (g\circ f)(s)\, dE_{J}(s)
\end{align*}
and hence \eqref{SEAL1} holds true. Choosing in particular $g(s) = s$, we deduce from \Cref{SSIDProp} that
\[
f(S) = \int_{\hh}f(s)\,dE_J(s) = \int_{\hh}s\,d\tilde{E}_{\tilde{J}}(s).
\]
By \Cref{ScalarSpecInt}, $f(S)$ is a scalar operator with spectral decomposition  $(\tilde{E},\tilde{J})$.

\Cref{IntCalcSer} implies $f(T) = f(S) + \Theta$ with
\[
\Theta := \sum_{n=1}^{+\infty} N^n \frac{1}{n!} (\sderiv^nf)(S).
\]
If we can show that  $\Theta$ is a quasi-nilpotent operator, then the statement of the theorem follows from \Cref{DecompTHM}. We first observe that each term in the sum is a quasi-nilpotent operator because $N^n$ and $(\sderiv^n f)(S)$ commute due to \cref{Zwei} in \Cref{NewSCalcProp} such that
\[ 0 \leq \lim_{k\to\infty} \left\| \left( N^n \frac{1}{n!} (\partial_S^nf)(S)\right)^k\right\| ^{\frac{1}{k}} \leq \left\|\frac{1}{n!} (\partial_S^nf)(S)\right\|  \left( \lim_{k\to\infty} \left\|N^{nk}\right\|^{\frac{1}{nk}} \right)^n=0.  \]
\Cref{NilSpec} thus implies $\sigma_S\left( N^n\frac{1}{n!}(\sderiv^nf(S))\right) = \{0\}$. By induction we conclude from \Cref{Funf} in \Cref{NewSCalcProp} and \Cref{NilSpec} that for each $m\in\nn$ the finite sum $\Theta_1(m) := \sum_{n=1}^{m} N^n \frac{1}{n!} (\sderiv^nf)(S)$ is quasi-nilpotent and satisfies $\sigma_S(\Theta(m)) = \{0\}$. 

Since the series $\Theta$ converges in the operator norm, for any $\varepsilon >0$ there exists $m_{\varepsilon} \in\nn$ such that $\Theta_2(m_{\varepsilon}):= \sum_{n=m_{\varepsilon} +1}^{+\infty} N^n \frac{1}{n!} (\sderiv^nf)(S)$ satisfies $\|\Theta_2(m_{\varepsilon})\| < \varepsilon$. Hence $\sigma_S(\Theta_2(m_{\varepsilon})) \subset B_{\varepsilon}(0)$ and as $\Theta = \Theta_1(m_{\varepsilon})  + \Theta_2(m_{\varepsilon})$ and $\Theta_{1}(m_{\varepsilon})$ and $\Theta_{2}(m_{\varepsilon})$ commute, we conclude from \Cref{Funf} in \Cref{NewSCalcProp} that $\sigma_S(\Theta)\subset B_{\varepsilon}(0)$. As $\varepsilon>0$ was arbitrary, we find $\sigma_{S}(\Theta) = \{0\}$. By \Cref{NilSpec},  $\Theta$ is quasi-nilpotent.

We have shown that  $f(T) = f(S)  +\Theta $, that $f(S)$ is a scalar operator with spectral decomposition $(\tilde{E},\tilde{J})$ and that $\Theta$ is quasi-nilpotent. From \Cref{DecompTHM} we therefore deduce that $f(T)$ is a spectral operator with spectral decomposition $(\tilde{E},\tilde{J})$, that $f(S)$ is its scalar part and that $\Theta$ is its radical part. This concludes the proof.

\end{proof}

\begin{corollary}
Let $T\in\boundOP(\hil)$ be a spectral operator and let $f\in\intrin(\sigma_{S}(T))$. If $T$ is of type $m\in\nn$, then $f(T)$ is of type $m$ too.
\end{corollary}
\begin{proof}
If $T=S+N$ is the decomposition of $T$ into its scalar and its radical part and $T$ is of type $m$ such that $N^{m+1} = 0$, then the radical part $\Theta$ of $f(T)$ is due to \Cref{REMThm} and \Cref{Kakao} given by
\[
 \Theta = f(T) - f(S) = \sum_{n=1}^{+\infty} N^n \frac{1}{n!} (\sderiv^n f)(S) =  \sum_{n=1}^{m} N^n \frac{1}{n!} (\sderiv^n f)(S) .
\]
Obviously also $\Theta^{m+1} = 0$.

\end{proof}

\section{Concluding Remarks}\label{Conclusions}
Several of our results, in particular \Cref{ABC}, \Cref{IntCalcThm}, \Cref{BackCompat} and \Cref{SOpBackwards}, show a deep relation between complex and quaternionic operator theory. Indeed, if we choose any imaginary unit $\I\in\SS$ and consider our quaternionic (right) vector space $V_{R}$ as a complex vector space over $\cc_{\I}$, then the quaternionic results coincide with the complex linear counterparts after suitable identifications: the spectrum $\sigma_{\cc_{\I}}(T)$ of $T$ as a $\cc_{\I}$-linear operator on $V_{R}$ equals $\sigma_{S}(T)\cap\cc_{\I}$, the $\cc_{\I}$-linear resolvent can be computed from the quaternionic pseudo-resolvent and vice versa and the quaternionic $S$-functional calculus for intrinsic slice hyperholomorphic functions coincides with the $\cc_{\I}$-linear Riesz-Dunford functional calculus. Similarly, spectral integration with respect to a quaternionic spectral system on $V_{R}$ is equivalent to spectral integration with respect to a suitably constructed $\cc_{\I}$-linear spectral measure and a quaternionic linear operator is a quaternionic spectral operator if and only if it is a spectral operator when considered as a $\cc_{\I}$-linear operator. Once again these relations show that the theory based on the $S$-spectrum and slice hyperholomorphic functions developed by Colombo, Sabadini and different co-authors is actually the right approach towards a mathematically rigorous extension of classical operator theory to the noncommutative setting. The relation with the complex theory presented in this article offers moreover the possibility of using the powerful tools of complex operator theory to study quaternionic linear operators. We refer for instance to \cite{HInftyComing}, where the complex Fourier transform was used in order to study spectral properties of the nabla operator, which is the quaternionification of the gradient and the divergence operator. 

We showed the correspondence of complex and quaternionic functional calculi only for intrinsic slice functions. This is however the best one can get: let $T$ be a quaternionic linear operator, let $\I\in\SS$ be an arbitrary imaginary unit and let $f:U\subset\cc_{\I}\to\cc_{\I}$ be a function suitable for the $\cc_{\I}$-linear version of   whichever functional calculus we want to apply. If $\bv$ is an eigenvector of $T$ associated with $\lambda\in\cc_{\I}$ and $\J\in\SS$ with $\I\perp\J$, then $\bv\J$ is an eigenvector of $T$ associated with $\overline{\lambda}$ as $T(\bv\J) = (T\bv)\J = (\bv\lambda )\J = (\bv\J)\overline{\lambda}$. The basic intuition of a functional calculus is that $f(T)$ should be obtained by applying $f$ to the eigenvalues resp. spectral values of $T$. Hence $\bv$ is an eigenvector of $f(T)$ associated with $f(\lambda)$ and $\bv\J$ is an eigenvector of $f(T)$ associated with $f\left(\overline{\lambda}\right)$. If however $f(T)$ is quaternionic linear, then again $f(T) (\bv\J) = (f(T)\bv)\J = (\bv f(\lambda))\J = (\bv\J)\overline{f(\lambda)}$. Hence, we must have $f\left(\overline{\lambda}\right) = \overline{f(\lambda)}$ in order to obtain a quaternionic linear operator. The slice extension of  $f$ is then an intrinsic slice function. For any function $f$ that does not satisfy this symmetry, the operator $f(T)$ will in general not be quaternionic linear. Hence we cannot expect any accordance with the quaternionic theory for such functions.

The class of intrinsic slice functions is however sufficient to recover all the spectral information about an operator: projections onto invariant subspaces obtained via the $S$-functional calculus are generated by characteristic functions of spectral sets, which are intrinsic slice functions. The continuous functional calculus is defined using intrinsic slice functions and even the spectral measure in the spectral theorem for normal operators can be constructed using intrinsic slice functions (cf. \Cref{PrelimSect} and \Cref{DiffSpInt} as well as \cite{Alpay:2015}, \cite{Ghiloni:2013} and \cite{Alpay:2016} for more details). Even more, any class of functions that can be used to recover spectral information via a functional calculus, necessarily consists of intrinsic slice functions. This is a consequence of the symmetry \eqref{EVSym} of the eigenvalues resp. the $S$-spectra of quaternionic linear operators. Indeed, as already mentioned above, the very fundamental intuition of a functional calculus is that $f(T)$ should be defined by action of $f$ on the spectral values of $T$. If in particular $\bv$ is an eigenvalue of $T$ associated with the eigenvalue $s$, then $\bv$ should be an eigenvector of $f(T)$ associated with the eigenvalue $f(s)$, i.e.
\begin{equation}\label{Intuition}
T\bv = \bv s\qquad\text{implies}\qquad f(T)\bv = \bv f(s).
\end{equation}
Let now $\I\in\SS$ be an arbitrary imaginary unit, let $s_{\I} = s_0 + \I s_1$ and let $h\in\hh\setminus\{0\}$ be such that $s_{\I} = h^{-1}s h$ as in \Cref{TurnCor}.
Then $\bv h$ is an eigenvalue of $T$ associated with $s_{\I}$ as $T(\bv h) = (T\bv)h = \bv s h = (\bv h) (h^{-1}sh)  = (\bv h) s_{\I}$. If \eqref{Intuition} holds true, then $f(T)\bv h = (\bv h) f(s_{\I})$. However, we also have
\[
f(T)(\bv h) = (f(T)\bv)h = \bv f(s) h = (\bv h) (h^{-1}f(s) h)
\]
and hence necessarily 
\begin{equation}\label{RoMp}
f(h^{-1}sh) = f(s_{\I}) = h^{-1}f(s)h.
\end{equation}
If we choose $h_{\I} = \I_{s}$, we find $f(s) = \I_{s}^{-1} f(s) \I_{s}$ and in turn $\I_{s}f(s) = f(s) \I_{s}$. A quaternion commutes with the imaginary unit $\I_{s}$ if and only if it belongs to $\cc_{\I_{s}}$. Hence, $f(s)  = \alpha + \I_{s}\beta$ with $\alpha,\beta\in\rr$. Now let again $\I\in\SS$ be arbitrary, set $s_{\I} = s_{0} + \I s_{1}$ and choose $h\in\hh\setminus\{0\}$ such that $\I = h^{-1}\I_{s} h$ and in turn  $s_{\I} =  s_0 + \I s_1 = s_0 + h^{-1} \I_{s} h s_1 = h^{-1}sh$. Equation \eqref{RoMp} implies then
\[
f(s_{\I}) = f(h^{-1}sh) = h^{-1}f(s)h = \alpha + h^{-1}\I_{s} h \beta = \alpha + \I \beta
\]
so that $f$ is an intrinsic slice function. 

The symmetry of the set of right eigenvalues and the relation \eqref{RightEVRel} between eigenvectors associated with eigenvalues in the same eigensphere therefore require that the class of functions to which a quaternionic functional applies consists of intrinsic slice functions, if this functional calculus respects the fundamental intuition of a functional calculus given in \eqref{Intuition}. This observation is reflected by resp. explains several phenomena in the existing theory.
\begin{itemize}
\item The $S$-functional calculus for left and for right slice hyperholomorphic functions are not always consistent. If $f$ is left and right slice hyperholomorphic, then both functional calculi will in general not yield the same operator unless $f$ is intrinsic. Any such function is intrinsic up to addition with a locally constant function, but for locally constant functions the two functional calculi disagree. This is due to the fact that invariant subspaces of right linear operators are in general only right linear subspaces and not closed under multiplication with scalars from the left \cite{Gantner:2017}.

\item Even if the functional calculus is defined for a class of functions that does not only contain intrinsic slice functions, the product rule only holds for the intrinsic slice functions in this class \cite{Colombo:2011, Alpay:2015, Ghiloni:2013, Alpay:2016c, Alpay:2017a}. An exception is \cite[Theorem~7.8]{Ghiloni:2013}, where the continuous functional calculus for $\cc_{\I}$-slice functions satisfies a product rule. This functional calculus is however not quaternionic in nature. Instead, it considers the operator as the quaternionic linear extension of a $\cc_{\I}$-linear operator on a suitably chosen $\cc_{\I}$-complex component space. The functional calculus can then be interpreted as applying the continuous functional calculus to the $\cc_{\I}$-linear operator on the component space and extending the obtained operator to the entire quaternionic linear space.

\item Similarly, the spectral mapping theorem only holds for those functions in the admissible function class that are intrinsic slice functions, cf. \cite{Colombo:2011, Ghiloni:2013, Alpay:2016c}. Again \cite[Theorem~7.8]{Ghiloni:2013} obtains a generalized spectral mapping property, that holds because this functional calculus is actually a functional calculus on a suitably chosen complex component space.
\end{itemize}

Finally, the above observations are another strong argument why we believe that the theory of spectral integration using spectral systems as presented in \Cref{SpecIntSect} is the correct approach to spectral integration that properly reflects the underlying intuition. In the introduction of \cite{Ghiloni:2017}, Ghiloni, Moretti and Perotti argue that the approach of spectral integration in \cite{Alpay:2016} is complex in nature as it only allows to integrate $\cc_{\I}$-valued functions defined on $\cc_{\I}^{\geq}$ for some $\I\in\SS$. They argue that their  approach using iqPVMs  on the other hand is quaternionic in nature as it allows to integrate functions that are defined on a complex halfplane and take arbitrary values in the quaternions.

 We believe that it is rather the other way around. It is the approach to spectral integration using spectral systems that is quaternionic in nature although they only allow to integrate intrinsic slice functions. This integration has however a clear interpretation in terms of the right linear structure on the space, which is the structure that is important for right linear operator. Moreover, unlike the iqPVM in \cite[Theorem A]{Ghiloni:2017}, the spectral system of a normal operator is fully determined by the operator. Finally, extending the class of integrable functions towards non-intrinsic slice functions does not seem to bring any benefit as any measurable functional calculus for such functions violates the very basic intuition \eqref{Intuition} of a functional calculus as we pointed out above.

\end{document}